\documentclass{article}

\usepackage[letterpaper, left=1.25in, right=1.25in, top=1in, bottom=1in]{geometry}

\usepackage[english]{babel}
\usepackage[autostyle]{csquotes}  
\usepackage[utf8]{inputenc}       
\usepackage[T1]{fontenc}          

\usepackage{amsmath, amsfonts, amssymb, mathtools}
\usepackage{bm}                   
\usepackage{mathrsfs}             
\usepackage{dsfont}               

\usepackage{graphicx}
\usepackage[section]{placeins}    
\usepackage{float}                
\usepackage{subcaption}           
\usepackage{threeparttable}       
\usepackage{multirow}             
\usepackage{colortbl}             
\usepackage[dvipsnames]{xcolor}   

\usepackage{algorithm}
\usepackage[noend]{algpseudocode}

\usepackage{cite}                 

\usepackage{setspace}             
\usepackage{enumitem}             
\usepackage{indentfirst}          
\usepackage{authblk}              

\newcommand{\bx}{\mathbf{x}}
\newcommand{\bb}{\mathbf{b}}
\newcommand{\bA}{\mathbf{A}}

\newcommand{\bR}{\mathbf{R}}
\newcommand{\bn}{\mathbf{n}}

\newcommand{\D}{\mathrm{d}}
\newcommand{\p}{\partial}

\newcommand{\be}{\begin{equation}}
\newcommand{\ee}{\end{equation}}
\newcommand{\bea}{\begin{eqnarray}}
\newcommand{\eea}{\end{eqnarray}}
\newcommand{\ben}{\begin{eqnarray}}
\newcommand{\een}{\end{eqnarray}}

\title{Energy Dissipation Rate Guided Adaptive Sampling for Physics-Informed Neural Networks: Resolving Surface-Bulk Dynamics in Allen-Cahn Systems}
\author[1]{Chunyan Li}
\author[2]{Wenkai Yu}
\author[3]{Qi Wang\thanks{Corresponding author: qwang@math.sc.edu}}
\affil[1]{Department of Mathematics, The Pennsylvania State University, University Park, PA, 16802, USA}
\affil[2]{Department of Mathematics, The Hong Kong University of Science and Technology Clear Water Bay, Kowloon, Hong Kong, P. R. China}
\affil[3]{Department of Mathematics, University of South Carolina, Columbia, SC, 29208, USA}
\date{}

\begin{document}

\maketitle
\begin{abstract}
We introduce the Energy Dissipation Rate guided Adaptive Sampling (EDRAS) strategy, a novel method that substantially enhances the performance of Physics-Informed Neural Networks (PINNs) in solving thermodynamically consistent partial differential equations (PDEs) over arbitrary domains. EDRAS leverages the local energy dissipation rate density as a guiding metric to identify and adaptively re-sample critical collocation points from both the interior and boundary of the computational domain. This dynamical sampling approach improves the accuracy of residual-based PINNs by aligning the training process with the underlying physical structure of the system.
In this study, we demonstrate the effectiveness of EDRAS using the Allen–Cahn phase field model in irregular geometries, achieving up to a sixfold reduction in the relative mean square error compared to traditional residual-based adaptive refinement (RAR) methods. Moreover, we compare EDRAS with other residual-based adaptive sampling approaches and show that EDRAS is not only  computationally more efficient but also more likely to identify high-impact collocation points. Through numerical solutions of the Allen–Cahn equation with both static (Neumann) and dynamic boundary conditions in 2D disk- and ellipse-shaped domains solved using PINN coupled with EDRAS, we gain significant insights into how dynamic boundary conditions influence bulk phase evolution and thermodynamic behavior. The proposed approach offers an effective, physically informed enhancement to PINN frameworks for solving thermodynamically consistent models, making PINN a robust and versatile computational tool for investigating complex thermodynamic processes in arbitrary geometries.
 
\end{abstract}
\textbf{Keywords}: {Physics-informed neural network (PINN), Energy Dissipation Rate-guided Adaptive Sampling (EDRAS) Strategy, Dynamic Boundary Conditions, Thermodynamically Consistent Models, Allen-Cahn equation.}

\section{Introduction}\label{sec: intro}
Physics-Informed Neural Networks (PINNs) offer an alternative approach to scientific computing by embedding physical laws directly into deep learning through automatic differentiation \cite{wright1999numerical, raissi2019physics}. PINNs have demonstrated remarkable potential in solving complex partial differential equation (PDE) problems across diverse fields, especially for partial differential equations in complex domains. However, PINNs struggle to accurately capture the dynamics of some partial different equations, such as Allen-Cahn equations, due to the multiscale nature of the systems and the frequency principle in neural networks \cite{xu2024}. The frequency principle describes the tendency of deep neural networks to fit target functions from low to high frequencies in tandem during the training process.  Meanwhile, the performance of neural networks can be hindered by the insufficient learning of high-frequency structure. Then, one has to make extra efforts to learn high-frequency structures to reach the desired accuracy. There are three key components when using PINNs: neural network structures integrated with physical laws, training process/optimizers and sampling strategies. One can improve the performance of PINNs in the three ways, respectively or all together. For example, inspired by the frequency principle of neural networks, a multiscale neural network is proposed to approximate the Green's function in linear PDEs \cite{hao2024multiscale}. Another multi-scale deep neural networks (MscaleDNNs) using the idea of radial scaling in frequency domain and activation functions with compact support are proposed for  the Poisson-Boltzmann equation \cite{CiCP-28-5}. In the optimization process, for instance, the training strategy based on the homotopy dynamics is proposed to learn the sharp interface of steady state solution of the Allen-Cahn equation in \cite{chen2025learn}.  The last area one can seek improvements is the sampling strategy which is the focus of this work.

The performance of PINN-based solvers critically depends on the selection of collocation points where the physical laws and constraints are enforced \cite{nabian2021efficient, mao2023, han2022residual, guo2024tcas}. Initially, non-adaptive sampling strategies such as spatially uniform grids, uniformly random sampling, Latin hypercube sampling \cite{stein1987large}, and quasi-random low-discrepancy sequences, such as the Halton sequence \cite{halton1960efficiency, faure2009generalized}, Hammersley sequence \cite{hammersley1960monte, hammersley2013monte}, and Sobol sequence \cite{sobol1967distribution, joe2008constructing} were widely used in PINNs. While straightforward to implement, these methods often fail to efficiently capture the complex dynamics of the underlying systems. A comprehensive review of these non-adaptive approaches and their impact on PINNs can be found in \cite{wu2023comprehensive}.

In contrast, adaptive sampling strategies dynamically adjust collocation points during training process to better capture critical solution features. Residual-based adaptive refinement (RAR) is one such method that leverages residual information to enhance the  performance of PINNs for PDEs   \cite{lu2021deepxde}. More specifically, the collocation points in regions with large residuals are added into the training set adaptively during the training process.  Later on, several distribution-based adaptive sampling methods are proposed. For example, importance sampling strategies are proposed to select collocation points based on a distribution proportional to the loss function or its approximation \cite{nabian2021efficient}.  Evolutionary sampling strategies evolve collocation points  iteratively, retaining high-residual points and re-sampling others to avoid trivial solutions \cite{daw2022rethinking}.  Inspired by these distribution-based adaptive sampling methods, the same group who proposed RAR method presents yet other two methods, residual-based adaptive distribution (RAD) and the combination of RAR and RAD which is named as the residual-based adaptive refinement with distribution (RAR-D). These methods dynamically redistribute collocation points in training set based on global residual distribution during training process\cite{wu2023comprehensive}.

The RAR method has a drawback, that is, only collocation points with extremely large residuals can be added to training set but not the collocation points with relatively small residuals even though the solution at these points are still far away from the true solution. This observation is demonstrated and analyzed in more detail in the current study. This drawback can be alleviated by introducing a distribution based on residuals as used in RAD and RAR-D methods. Since the residuals are scaled by the global residual distribution so that the collocation points are distributed more balanced in the sense that the points with relatively smaller residuals but large approximation errors can also be sampled into training set. However, residual-distribution based method may require lots of computational resource to estimate the global distribution.

Novel sampling methods that can alleviate or avoid the drawbacks of RAR method and also more computationally efficient than RAD method are desired. Recently, imbalanced learning-based adaptive sampling method, named as residual-based Smote (RSmote), is proposed \cite{Yangyahong} to improve the performance of PINNs on various PDE systems. The main idea is to classify the training points into two groups: group A with larger residuals while  group B with smaller residuals. Usually, these two groups are imbalanced during the training process, then one can use Synthetic Minority Over-sampling Technique (SMOTE) to generate more points belongs to minority group A so that one can focus on the areas with larger residuals during the training process. No distribution estimation is needed in this method so that the memory usage is significantly reduced. There are also some other adaptive sampling methods such as temporal causality-based adaptive sampling methods which further refine training by dynamically adjusting sampling ratios, accounting for both PDE residuals and temporal causality within sub-domains \cite{guo2024tcas}. Similarly, gradient-based adaptive sampling incorporates solution gradients to detect discontinuities and shocks \cite{mao2023}.

Despite their success, the residual distribution based adaptive sampling methods such as RAD and RAR-D are not computationally efficient enough and not suitable for high-dimensional PDEs, since the estimation of the normalization factor $\int_{\Omega}R(y)dy$ ($R(y)$ denotes residual density) requires more and more memories and computational resources and can become computationally intractable as the dimension increases. Moreover, these adaptive methods primarily rely on residuals and gradients, overlooking critical information about the dynamical process in the physical law, especially the energy landscape—a fundamental aspect of thermodynamically consistent models. Addressing this limitation is crucial for advancing the application of PINNs to effectively solve such systems.

 In this paper, we address these limitations by introducing the Energy Dissipation Rate-guided Adaptive Sampling (EDRAS) method, a novel strategy that leverages the system's energy landscape to optimize the selection of the collocation points. EDRAS dynamically adjusts sampling density based on energy dissipation rate density, concentrating computational resources in regions of significant energy variation. No global distribution density estimation is needed. The energy dissipation rate density calculation are often more efficient than the calculation of residuals due to the fact that the residuals usually involves more derivatives or higher order derivatives compared with the energy dissipation rate density. To demonstrate EDRAS's effectiveness, we apply it to the Allen-Cahn phase field model with periodic, Neumann and dynamic boundary conditions, respectively,  examining cases in 1D and in complex 2D geometries (disk and elliptic domains).

This study carries particular significance for systems where surface-bulk coupling plays a crucial role, such as in cellular dynamics where membrane surface area dominates bulk volume, phase separation within living cells, surface dynamics-controlled mesoscopic systems, etc. \cite{espath2023continuum}. In such cases, traditional no-flux boundary conditions prove insufficient, necessitating dynamic boundary conditions that can capture the complex interplay between surface and bulk dynamics \cite{espath2023continuum}. Our approach using PINNs coupled with EDRAS offers distinct advantages for these systems, especially in arbitrary domains where traditional numerical methods face significant mesh-related challenges.

The theoretical foundations of dynamic boundary conditions have been extensively developed in recent literature. Systematic approaches to deriving thermodynamically consistent dynamic boundary conditions have emerged through the application of the generalized Onsager principle \cite{xiaobocms, chunyan2023thesis, xiaobocms}. For the Allen-Cahn equation with dynamic boundary conditions, fundamental mathematical properties including local existence and uniqueness have been established in the non-isothermal case \cite{gal2008non}. The mathematical theory extends further to Cahn-Hilliard type equations with dynamic boundary conditions, where researchers have proven global existence and uniqueness for both regular and irregular potentials, while also providing comprehensive asymptotic analysis \cite{racke2003cahn, wu2004convergence, pruss2006maximal, gal2006cahn, chill2004convergence, gilardi2009cahn}. These theoretical developments establish the mathematical well-posedness of the phase field models with dynamic boundary conditions, providing a rigorous foundation for numerical investigations of surface-bulk coupled dynamics.

Despite its simplicity, the EDRAS method has demonstrated significant utility. Extensive experiments show that applying EDRAS to the Allen-Cahn phase field equation markedly improves the performance of PINNs. We tersely summarize the main contributions of this study as follows.
\begin{itemize}
	\item \textbf{Broad applicability:} The EDRAS method can be seamlessly integrated into any neural network-based solver for thermodynamically consistent models in arbitrary domains, where Monte-Carlo approximation is used, as long as the free energy functional of the system is given.
	\item \textbf{Innovative sampling:} EDRAS introduces a novel adaptive sampling technique for PINN-based PDE solvers, avoiding the need for explicit distribution definitions or complex numerical approximations.
	\item \textbf{Computational efficiency:} EDRAS is straightforward to implement and computationally efficient, as it avoids the need to track global distributions or approximate them with additional dense sampling points. Moreover, the energy dissipation rate density guided sampling demonstrates superior computational efficiency compared to residual-based adaptive methods. This advantage arises because estimating the energy dissipation rate density requires fewer derivatives—or lower-order derivatives—than residual estimation, significantly reducing the computational overhead.
	\item \textbf{Enhancing performance of PINNs:} EDRAS serves as a powerful enhancement to PINNs for solving thermodynamically consistent systems.
	\item \textbf{Revealing insights into dynamic boundary:} This work provides fresh insights into the critical relationship between dynamic boundary conditions and bulk dynamics in thermodynamic consistent  systems.
\end{itemize}
The remainder of this paper is organized as follows: Section \ref{sec:model} presents the EDRAS method and its implementation for an Allen-Cahn model with dynamic boundary conditions in an arbitrary domain. Section \ref{sec:numerial} demonstrates the effectiveness of the EDRAS method and examines the surface-bulk coupled dynamics in the Allen-Cahn equation with dynamic boundary conditions in both 1D and 2D settings. Finally, Section \ref{sec:conclusion} summarizes the findings and discusses their implications.

\section{Energy Dissipation Rate-Guided Adaptive Sampling Method for Thermodynamically Consistent Models}\label{sec:model}
This section begins by introducing the thermodynamically consistent Allen-Cahn equation with thermodynamically consistent dynamic boundary conditions. We then provide a brief overview of the Physics-Informed Neural Networks (PINNs) applied to this surface-bulk coupled system. Subsequently, we discuss the motivation behind the Energy Dissipation Rate-Guided Adaptive Sampling Method (EDRAS) method, then, we outline the EDRAS algorithm in details.
\subsection{Thermodynamically Consistent Models with Thermodynamically Consistent Dynamic Boundary Conditions and PINNs}
We consider the free energy functional for a binary material system in domain $\Omega$ and along its boundary $\partial \Omega$, given by the following
\begin{equation}
    E = \int_\Omega [\frac{\epsilon^2}{2}|\nabla \phi|^2 + f(\phi)] \D \bx + \int_{\p\Omega} [\frac{\epsilon_s^2}{2}|\nabla_s\phi|^2 + g(\phi)] \D S,
\end{equation}
where $\epsilon$ is a parameter that controls the strength of the conformational entropy in the bulk, $\epsilon_s$ is a parameter that controls the strength of the conformational entropy on the bounding surface, $f(\phi)$ and $g(\phi)$ represent the bulk and surface free energy density, respectively, $\nabla$ and $\nabla_s$ are the gradient and surface gradient operator, respectively. A thermodynamically consistent model within $\Omega$, along with dynamic boundary conditions on $\p\Omega$, can be derived by hierarchically applying the generalized Onsager principle in tandem to the energy dissipation rate functional \cite{chunyan2023thesis, liu2021frontiers, xiaobocms,xiaoboentropy}, resulting in the following governing system of equations in the bulk and on the boundary, respectively,
\begin{equation}\label{eq:AC_2d}
    \begin{cases}
         \phi_t = -M_b \mu, \quad \bx \in \Omega, t\in (0, T],\\
         \mu = -\epsilon^2 \nabla^2\phi + f'(\phi), \\
        \phi_t = -M_s\mu_s, \quad \bx \in \p \Omega, t\in (0, T],\\
        \mu_s = -\epsilon_s^2\nabla^2_s\phi + \epsilon^2\bn \cdot \nabla\phi + g'(\phi), \\
    \end{cases}
\end{equation}
where $\mu$ and $ \mu_s$ are the bulk and surface chemical potential, $M_b\geq 0$ and $ M_s\geq 0$ are the bulk and surface mobility coefficient, respectively, and $T$ is the terminal time of the interest. Model \eqref{eq:AC_2d}  exhibits structures of the Allen-Cahn equation both in the bulk and on the bounding surface. Moreover, the model in the bulk, combined with given dynamic boundary conditions, is thermodynamically consistent, as the system obeys the second law of thermodynamics in the closed domain. Consequently, the system obeys the energy dissipation rate given by
\begin{equation}
    E'(t) = \int_\Omega -\frac{\phi_t^2}{M_b}\D\bx + \int_{\p \Omega}-\frac{\phi_t^2}{M_s}\D S = \int_\Omega -M_b
    \mu^2 \D\bx + \int_{\p \Omega}-M_s\mu_s^2\D S \leq 0.
\end{equation}

In this study, we demonstrate how to solve this initial-dynamic-boundary-value problem in an arbitrary domain using PINNs and investigate the impact of dynamic boundary conditions  on bulk dynamics in system  \eqref{eq:AC_2d}.  PINNs offer a mesh-free approach, suitable for handling initial-boundary value partial differential equation problems (IBV-PDEs) defined in complex domains, in which one solves IBV-PDEs through a least-square method by representing the solution using a deep neural network. The loss function in PINNs is often defined as the sum of mean squared errors of the residuals of the PDE system, initial conditions, and boundary conditions. Specifically, for a training set of sampled collocation points $\mathcal{J}$ consisting of  interior points  $\mathcal{J}_f$, boundary points $\mathcal{J}_b$, and initial points $\mathcal{J}_i$, the loss function is defined by
\begin{equation}\label{eq:loss}
    \mathcal{L}(\theta, \mathcal{J}) = w_f\mathcal{L}_f(\theta, \mathcal{J}_f) + w_b\mathcal{L}_b(\theta, \mathcal{J}_b) + w_i\mathcal{L}_i(\theta, \mathcal{J}_i),
\end{equation}
where $\theta$ represents the collection of the parameters of the neural networks, $w_f, w_b, w_i$ are weights for the residuals of the PDE $\mathcal{L}_f$, boundary error $\mathcal{L}_b$ and initial condition error $\mathcal{L}_i$,  respectively. Here, for the solution represented as a deep neural network, $\phi(\bx, t, \theta)$, each residual is defined, respectively, by
\begin{equation}
    \mathcal{L}_f(\theta, \mathcal{J}_f) = \frac{1}{|\mathcal{J}_f|}\sum_{(\bx, t)\in \mathcal{J}_f}|\phi_t(\bx, t; \theta)+M_b\mu(\phi(\bx, t; \theta))|^2,
\end{equation}

 \begin{equation}
     \mathcal{L}_b(\theta, \mathcal{J}_b) = \frac{1}{|\mathcal{J}_b|} \sum_{(\bx^{bc}, t)\in \mathcal{J}_b}|\phi_t(\bx^{bc}, t; \theta)+M_s\mu_s(\phi(\bx^{bc}, t; \theta))|^2,
 \end{equation}

 \begin{equation}
     \mathcal{L}_i(\theta, \mathcal{J}_i) = \frac{1}{|\mathcal{J}_i|} \sum_{(\bx, 0)\in \mathcal{J}_i}|\phi(\bx, 0; \theta)-\mathcal{I}(\phi(\bx))|^2,
 \end{equation}
where $|(\bullet)|$ indicates the number of points in set $(\bullet)$ and the initial condition is given as $\phi(\bx, 0)=\mathcal{I}(\phi(\bx))$.

The solution of the IBV-PDE problem is obtained by minimizing the loss function using a stochastic gradient descent based optimizer, such as the Adam optimizer \cite{kingma2014adam}, with randomly sampled training set $\mathcal{J}$ of collocation points. It is well-known that the selection of the sampled collocation points can significantly influence the performance of PINNs \cite{nabian2021efficient}. In \cite{lu2021deepxde, wu2023comprehensive}, an adaptive Residual based Adaptive Refinement (RAR) method was proposed to  sample the training points adaptively based on the residuals/errors to improve the performance of PINNs. However, the RAR method is not always effective for any PDEs. For example, it does not work well with the initial-boundary-value Allen-Cahn model we formulated above.

\subsection{Energy Dissipation Rate Density-Guided Adaptive Sampling Method}
We use the Allen-Cahn equation in 1D given below as an example to illustrate why the RAR method are not effective enough for thermodynamically consistent models and discuss the motivation for the EDRAS method and the EDRAS algorithms in details.
\begin{equation}
\left \{
\begin{split}
      &\phi_t = 0.0001\phi_{xx}-5\phi^3+5\phi, \quad x\in [-1, 1], \quad t\in (0, 1],\\
      &\phi(0, x) =x^2\cos(\pi x),\\
      &\phi(t, -1)=\phi(t, 1),\quad \phi_x(t, -1) = \phi_x(t, 1).
\end{split}\right.
\end{equation}
Let $\phi$ be the reference solution obtained from the finite difference (FDM) method, $\hat{\phi}$ the approximate solution obtained from the PINN method, $R(x, t)=|\hat{\phi}_t -0.0001\hat{\phi}_{xx}+5\hat{\phi}^3-5\hat{\phi}|$ the absolute value of the residual of the approximate solution at $(x, t)$ and $e(x,t)=|\phi(x,t)-\hat{\phi}|$. 

\begin{figure}[ht]
    \centering
        \includegraphics[scale=0.5]{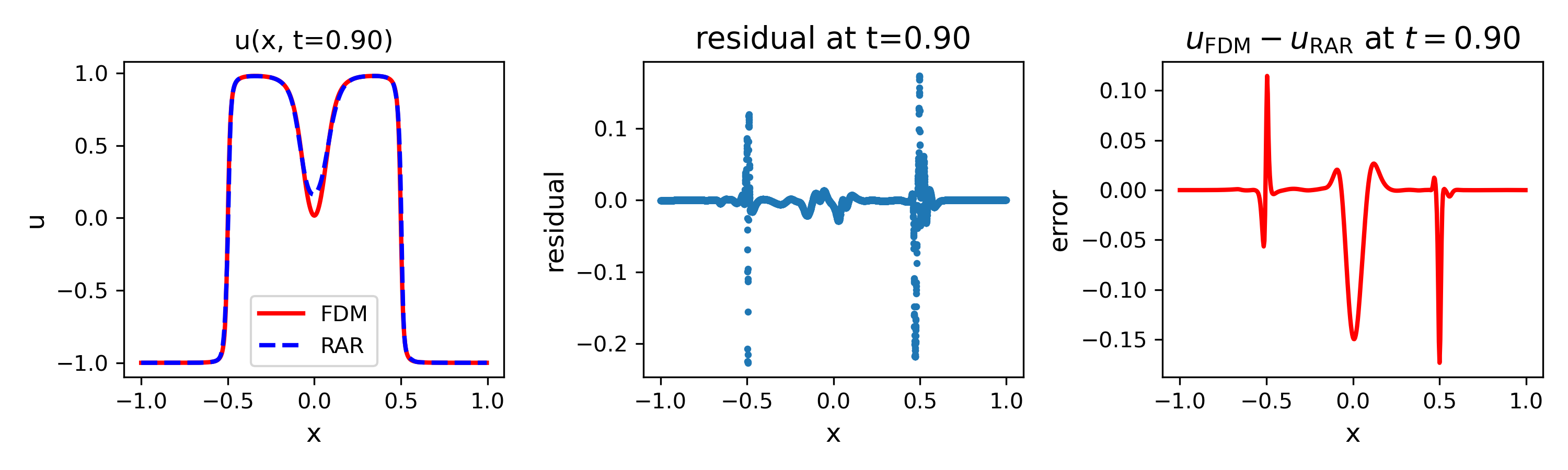}
        \caption{Left: Numerical solutions of the Allen-Cahn equation at $t=0.9$ obtained from a finite difference method (FDM) and an RAR-enhanced PINN. Middle: The residual density of the Allen-Cahn equation with the RAR-enhanced PINN solution at $t=0.9$. Right: Error (the solution difference between the FDM and the RAR-enhanced PINN solution) at $t=0.9$. It's evident that the RAR-enhanced PINN solution losses accuracy at $x=0, \pm 0.5$, respectively.}
        \label{fig:motivation}
  \end{figure}

As a motivation,  we  draw an analogy by considering a linear system, $\bA\cdot \bx=\bb$, $\bA\in \bR^{n\times n}, \bx, \bb\in \bR^n$. Solving the linear system using a residual-based method, such as GMRES, can lead to a small residual  $R=A(\bx-\Tilde{\bx})$, yet the relative error $e=\frac{||\bx-\Tilde{\bx}||}{||\bx||}$  may still be large if the condition number of $A$ is high. In the Allen-Cahn equation, the error $e(x, t)$ near $x=0$ is large despite a small residual density $R(\bx, t)$ (see Figure \ref{fig:motivation}). This suggests that eliminating the error $e(x, t)$ near $x=0$ solely through sampling interior points based on the residual error may not be sufficient enough. As it is shown in Figure \ref{fig:motivation}, we can easily classify the collocation points into two groups: in group one, the points have relatively large residuals as well as relatively large approximation errors. The group two contains the points that have relatively small residuals but relatively large approximation errors. Just as the area near $x=\pm\frac{1}{2}$ (group 1) vs. area near $x=0$ (group 2). Practically, the RAR method is more likely to add more points to group one instead of group two. To circumvent this issue, inspired by the distribution based adaptive sampling methods, one proposed residual-based adaptive distribution (RAD) and residual-based adaptive refinement with distribution (RAR-D) which is the combination of RAR and RAD. By introducing a distribution of residual normalized by $\int_{\Omega}R(x, t)dxdt$, RAD method can balance the probability of group one and group two being sampled. However, this normalization factor $\int_{\Omega}R(\bx, t)dxdt$ becomes harder and harder to estimate or approximate and even intractable as the dimension of $x$ increases.

\begin{figure}[ht]
     \centering
        \includegraphics[scale=0.3]{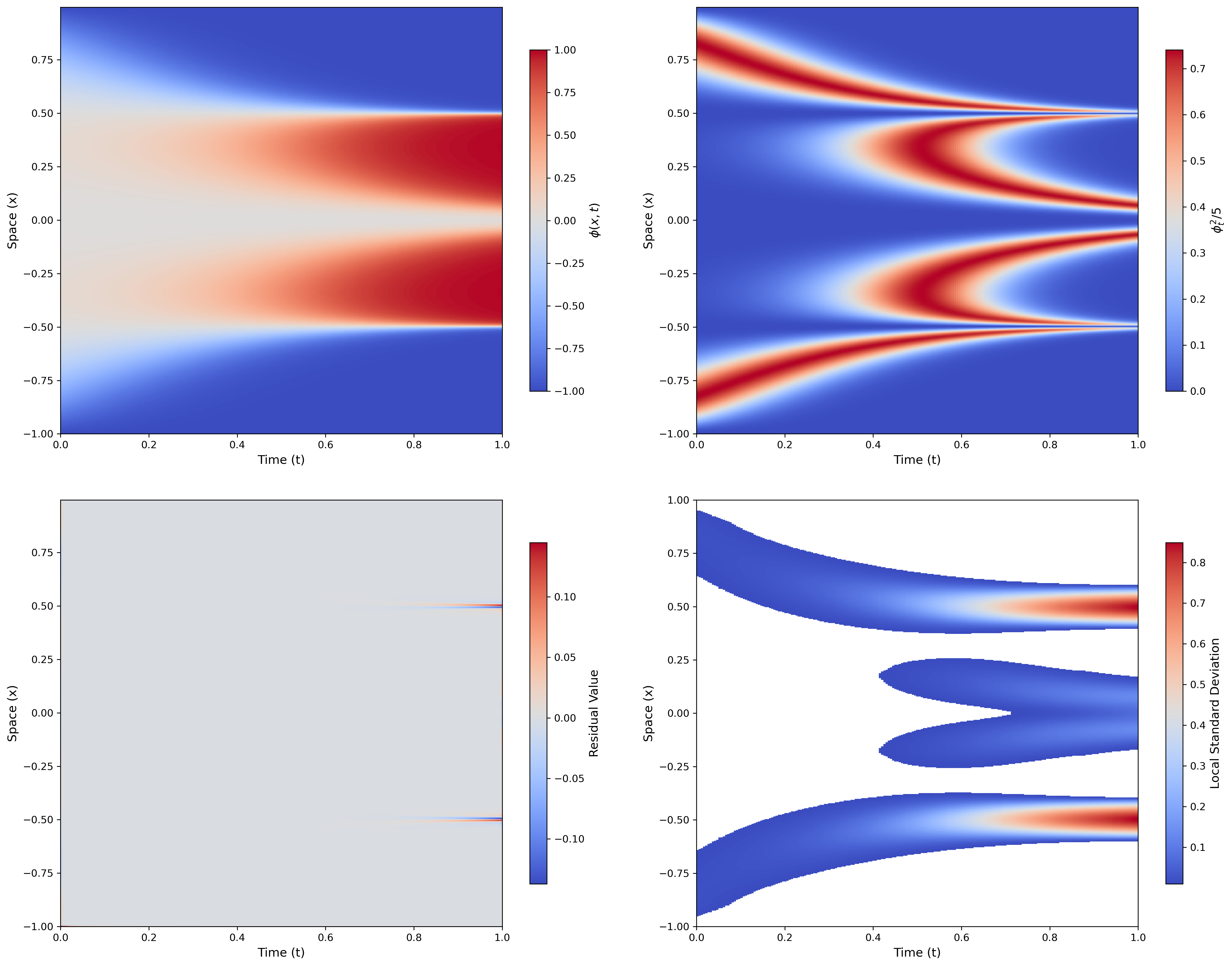}
        \caption{Top Left: The numerical solution obtained from FDM. Top Right: The energy dissipation rate density $\frac{\phi_t^2}{5}$ over space and time computed from the numerical solution using FDM.
Bottom Left: The residual of the Allen-Cahn equation over the spatial-temporal domain from the numerical solution using FDM. Bottom Right: Local standard deviation $\sigma$ of the numerical solution that exceeds threshold $\sigma>1e-2$.  There is a clear match between the large energy dissipation rate density and the large local standard deviation. The regions of high solution errors in the RAR-enhanced PINN model strongly correlate with areas exhibiting intense energy dissipation rate density fluctuations.}
    \label{fig:motivation1}
\end{figure}

These observations highlight the need for more effective sampling strategies in Physics-Informed Neural Networks (PINNs) when solving the thermodynamically consistent models such as the  Allen-Cahn model. For thermodynamically consistent PDE systems, the energy dissipation rate density serves as a key metric, quantifying spatiotemporal  dynamics of the system. Since the energy functional—a Lyapunov function for such systems—drives the solution toward local energy minima over time (see Figure \ref{fig:motivation1}), regions with high energy dissipation rate densities (e.g., near $x=0$ in the Allen-Cahn model) often correspond to areas where PINNs struggle to learn the solution, even with Residual-Based Adaptive Refinement (RAR) sampling. Moreover, it is natural to prioritize sampling in regions where the solution exhibits higher standard deviation. As illustrated in Figure \ref{fig:motivation1}, the energy dissipation rate density aligns closely with the spatiotemporal standard deviation of the solution. This connection motivates us to propose the Energy Dissipation Rate-Guided Adaptive Sampling (EDRAS) method. EDRAS uses the energy dissipation rate density as a guiding metric for adaptive sampling in PINNs, analogous to local mesh refinement in finite element methods or finite difference, where free energy and its dissipation rate act as error estimators \cite{verfurth1994posteriori, zienkiewicz2005finite}. Specifically, EDRAS targets regions (including boundaries) with high energy dissipation rate densities, where solutions exhibit significant variations. By focusing computational resources on these critical areas, EDRAS enhances the accuracy of PINNs when applied to thermodynamically consistent models as validated in numerical experiments in the next section.

Specifically, denoting $e_{\text{edrd}_f}(\mathbf{x}, t)$ as the absolute value of the energy dissipation rate density in the bulk, we utilize $e_{\text{edrd}_f}$ to identify areas with significant energy variations across a sampled spatial-temporal domain. The points corresponding to large $e_{\text{edrd}_f}$ are then incorporated into the training set.

The simplest approach is to augment the training set by adding the top $m$ points with the largest $e_{\text{edrd}_f}$ values from a dense candidate set $S$. For high-dimensional problems, where computational efficiency is crucial, one may also need to remove certain points from the training set. This can be accomplished in two ways:

\begin{itemize}
    \item Remove $m$ points with the smallest $e_{\text{edrd}_f}$ values from the current training set.
    \item Establish a threshold $\rho$ for $e_{\text{edrd}_f}(\mathbf{x}, t)$ and filter out points below this value.
\end{itemize}

However, such pruning may lead to insufficient point density in certain regions of the domain, compromising the training set's comprehensive coverage. To address this, we propose a density-aware strategy to maintain a minimum density threshold per unit area. This balanced approach prevents computational redundancy while ensuring adequate domain coverage. A detailed description of the algorithm is summarized in Algorithm 1 below. 

In practice, during the training process, let \(\mathcal{J}_f \) and \(\mathcal{J}_b \) be the training sets of interior points and boundary points, respectively. We define the following thresholds:

\begin{itemize}
    \item \textbf{Threshold of the bulk  energy dissipation rate density}:
    \[
    e_{\text{edrd}_{f0}} = \sum_{y \in \mathcal{J}_f} \frac{1}{3|\mathcal{J}_f|} e_{\text{edrd}_{f}}(x,t),
    \]
    where \( e_{\text{edrd}_{f}}(x,t) \) is the energy dissipation rate at an interior point \((x,t) \).

    \item \textbf{Threshold of the boundary  energy dissipation rate density}:
    \[
    e_{\text{edrd}_{b0}} = \sum_{y \in \mathcal{J}_b} \frac{1}{3|\mathcal{J}_b|} e_{\text{edrd}_{b}}(x,t),
    \]
    where \( e_{\text{edrd}_{b}}(x,t) \) is the energy dissipation rate at a boundary point \( (x, t) \).
\end{itemize}
Additionally, the \textbf{thresholds of the  collocation point density} in domain \( d^{(0)}_{f} \) and on boundary \( d^{(0)}_{b} \) can be determined based on the total number of points in a uniform grid, as typically used in the finite difference method.  In practices, we compute the energy dissipation rate density (both the bulk and surface) using the general flux $u_t$ rather than the generalized forces $\mu$ and $\mu_s$, as the latter involve higher-order derivatives and thus incur greater computational cost.

This methodology is simple, versatile and applicable to any thermodynamically consistent models in arbitrary domain $\Omega$. For any given $\Omega$, the interior points can be promptly acquired by initially sampling points within a larger rectangular domain, followed by a verification of their presence within the domain. The crucial part of the sampling lies in effectively sampling the boundary points, particularly in  our model, where dynamic boundary conditions are integral parts of the model. Ensuring accuracy and efficiency in this sampling practice for the boundary is paramount for optimizing the performance of PINNs. If a parameterized representation of the domain boundary is accessible, implementing the sampling approach becomes straightforward. In cases where such a representation is not available, an alternative method utilizes another feed-forward neural network, denoted as $\psi_{\Omega}(\bx; \theta)$, to approximate $H_{\Omega}$, the characteristic  function of domain $\Omega$. This function presents the shape of domain in a larger rectangular domain and replaces the sharp boundary of the domain with a continuous transition layer. Subsequently, $\psi_{\Omega}(\bx, \theta)=0.5$ serves as the boundary delineating the shape of the domain, allowing the bisection method to determine the boundary points.

\begin{algorithm}[ht]
    \label{Alg:alg1}
    \caption{EDRAS and density detection strategy for improving the distribution of collocation points for training PINNs.}
    \begin{algorithmic}[1]
        \State Given maximum iteration number $N_{iter}$, threshold for loss $\mathcal{L}_0$, threshold for density of collocation points in domain $d^{(0)}_{f}$ and on boundary $d^{(0)}_{b}$, 
          select the initial  interior points $\mathcal{J}_f$, boundary points $\mathcal{J}_b$ and initial points $\mathcal{J}_i$, and then train neural network, $u(\bx, t, \theta)$, for a limited number of iterations.
	\State Set the threshold of bulk  energy dissipation rate density $ e_{\text{edrd}_{f0}} = \frac{1}{3|\mathcal{J}_f|} \sum_{\mathbf{x} \in \mathcal{J}_f} e_{\text{edrd}_f}(\mathbf{x}, t, \theta)$, and the threshold of boundary  energy dissipation rate density
$e_{\text{edrd}_{b0}} = \frac{1}{3|\mathcal{J}_b|} \sum_{\mathbf{x}_b \in \mathcal{J}_b} e_{\text{edrd}_b}(\mathbf{x}_b, t, \theta).$ 
        \State Remove interior points with very small density of bulk energy dissipation rate density $e_{\text{edrd}_f}=\frac{\phi_t}{M_b}<e_{\text{edrd}_{f0}}$ from $\mathcal{J}_f$, remove boundary points with very small densities of boundary energy dissipation rate $e_{\text{edrd}_b}=\frac{\phi_t}{M_s}<e_{\text{edrd}_{b0}}$ from $\mathcal{J}_b$.  \{Alternatively, one can remove $m_f$ interior points with smallest $e_{\text{edrd}_f}$ and remove $m_b$ boundary points with smallest $e_{\text{edrd}_b}$ for prescribed $m_f$ and $m_b$. \} 
        \State Divide domain $\Omega$ into $N$ small subdomains $\Omega_i, i=1, ..., N$. Detect the density of training points $d_f^{(i)}$ and $d_b^{(i)}$ in each subdomain $\Omega_i$ and its boundary.
        \For{$i=1$ \textbf{to} $N$}
            \If {$d_f^{(i)}<d^{(0)}_{f}$ or $d_b^{(i)}< d^{(0)}_{b}$}
                \State Generate a dense set of interior points $S_f$ and a dense set of boundary points $S_b$.
                \State Add $d^{(0)}_f-d_f^{(i)}$ new interior points with the largest bulk energy dissipation rate density $e_{\text{edrd}_f}$ in $\mathcal{S}_f$ to $\mathcal{J}_f$ and add $d_b^{(0)}-d_b^{(i)}$ new boundary points with the largest boundary energy dissipation rate density $e_{\text{edrd}_b}$ in $\mathcal{S}_b$ to $\mathcal{J}_b$.
            \EndIf
        \EndFor
        \State Stop if $|\mathcal{L}^{n+1}-\mathcal{L}^{n}|<\mathcal{L}_0$ or $n_{iter}> N_{iter}$. Otherwise, go to step 2 and retrain the network.
    \end{algorithmic}
\end{algorithm}

\section{Results and Discussion}\label{sec:numerial}
This section is organized into three parts. First,we compare the performance of the EDRAS method with that of the Residual-based Adaptive Refinement (RAR) method while coupled with  PINNs in solving the Allen-Cahn equation in 1D with periodic boundary conditions in Subsection \ref{subsec:AC_1d}. Furthermore, we analyze and discuss limitations of RAR and RAD methods as well as advantages of EDRAS from a probabilistic perspective in Subsection \ref{subsec:AC_1d}. Second, in Subsection \ref{subsec:AC_neu}, we apply the EDRAS method to the Allen-Cahn equation in 2D with Neumann boundary conditions to demonstrate its effectiveness and examine the distinct impact of dynamic boundary conditions in comparison to Neumann boundary conditions. Finally, we further explore the influence of dynamic boundary conditions on the bulk dynamics of the system in 2D computational domains with disk and ellipse geometries in Subsection \ref{subsec:AC_dy}.

\subsection{Allen-Cahn Equation with Periodic Boundary Condition in 1D }\label{subsec:AC_1d}
We first use Allen-Cahn equation \eqref{eq:ac_1d} in 1D as an example to showcase the effectiveness of EDRAS in comparison with RAR proposed in \cite{lu2021deepxde}.
We consider the total free energy with a double well bulk potential as follows
\begin{equation}
    E(t) = \int_{\Omega} 0.00001|\nabla \phi|^2+\frac{1}{4}(\phi^2-1)^2 \D x.
\end{equation}
The model with a periodic boundary condition of period 2 reads
\begin{equation}\label{eq:ac_1d}
\left \{
\begin{split}
      &\phi_t = 0.0001\phi_{xx}-5\phi^3+5\phi, \quad x\in [-1, 1], \quad t\in (0, 1],\\
      &\phi(0, x) =x^2\cos(\pi x),\\
      &\phi(t, -1)=\phi(t, 1),\quad \phi_x(t, -1) = \phi_x(t, 1).
\end{split}\right.
\end{equation}
Then, the energy dissipation rate is given by
\begin{equation}
    \frac{dE}{dt} = \int_\Omega -5(\phi^3-\phi -0.00002 \nabla^2 \phi)^2dx = \int_\Omega -\frac{\phi^2_t}{5}\D x.
\end{equation}
We denote ${u(\bx, t;\theta)}$ as the solution obtained using PINN. Then, the loss function is defined as follows
\begin{equation}
\begin{split}
    \mathcal{L}(\theta) &= \frac{w_f}{|\mathcal{J}_f|}\sum_{(x, t)\in \mathcal{J}_f}|{u}_t(x, t;\theta)-(0.0001{u}_{xx}(x, t; \theta)-5{u}^3(x, t; \theta)+5{u}(x, t; \theta))|^2 \\
    &+ \frac{w_{b1}}{|\mathcal{J}_{b1}|}\sum_{\mathcal{J}_{b1}}|{u}(-1 ,t; \theta)-{u}(1, t; \theta)|^2
    +\frac{w_{b2}}{|\mathcal{J}_{b2}|}\sum_{\mathcal{J}_{b2}}|{u}_x(-1, t; \theta)-{u}_x(1, t; \theta)|^2 \\
    &+ \frac{w_i}{|\mathcal{J}_i|}\sum_{(0, x)\in \mathcal{J}_i}|{u}(x, 0; \theta) - x^2\cos(\pi x)|^2.
\end{split}
\end{equation}

We implement a time-marching training strategy\cite{bcpinn-mattey2022novel, Jia2020adapinn, wang2024respecting, krishnapriyan2021characterizing} using a sequence of deep neural networks (DNNs) $\{u_i(\bx, t; \theta)\}_{i=1}^n$ to approximate solutions within discrete time segments $t\in [T_{i-1}, T_{i}], i=1,..., n$. In our implementation, we divide the temporal domain into $n=6$ time segments with endpoints $T_i=0.01, 0.2, 0.4, 0.6, 0.8, 1.0$ for $i=1, \cdots, 6$, respectively. The choice of a small initial segment $T_1=0.01$ ensures the accurate learning of the solution during the initial short time period, which is critical for this dissipative system to maintain solution accuracy across the entire temporal domain.  Each DNN employs 3 hidden layers with 128 nodes per layer and the hyperbolic tangent activation function in each hidden layer, no activation function is used in the output layer.
The training process consists of several stages:
\begin{itemize}
\item Initial sampling: 1000 collocation points, 514 initial points, and 200 boundary points are randomly sampled from a uniform distribution.
\item Adaptive refinement: Every 40 epochs, 100 additional collocation points are added using either EDRAS, RAR, or their combination, until reaching 3000 additional points.
\item Optimization: Training proceeds with Adam optimizer (3000 epochs, learning rate 0.001, mini-batch size 32) followed by L-BFGS optimizer (maximum 50000 iterations).
\item Loss weighting: The loss function weights are set to $w_i=100$ for initial conditions, $w_f=1$ for the PDE residual, and $w_{b1}=1, w_{b2}=50$ for the boundary conditions in equation (\ref{eq:ac_1d}).
\end{itemize}

\begin{figure}[H]
  \centering
  \begin{minipage}[b]{0.3\textwidth}
       \includegraphics[scale=0.3]{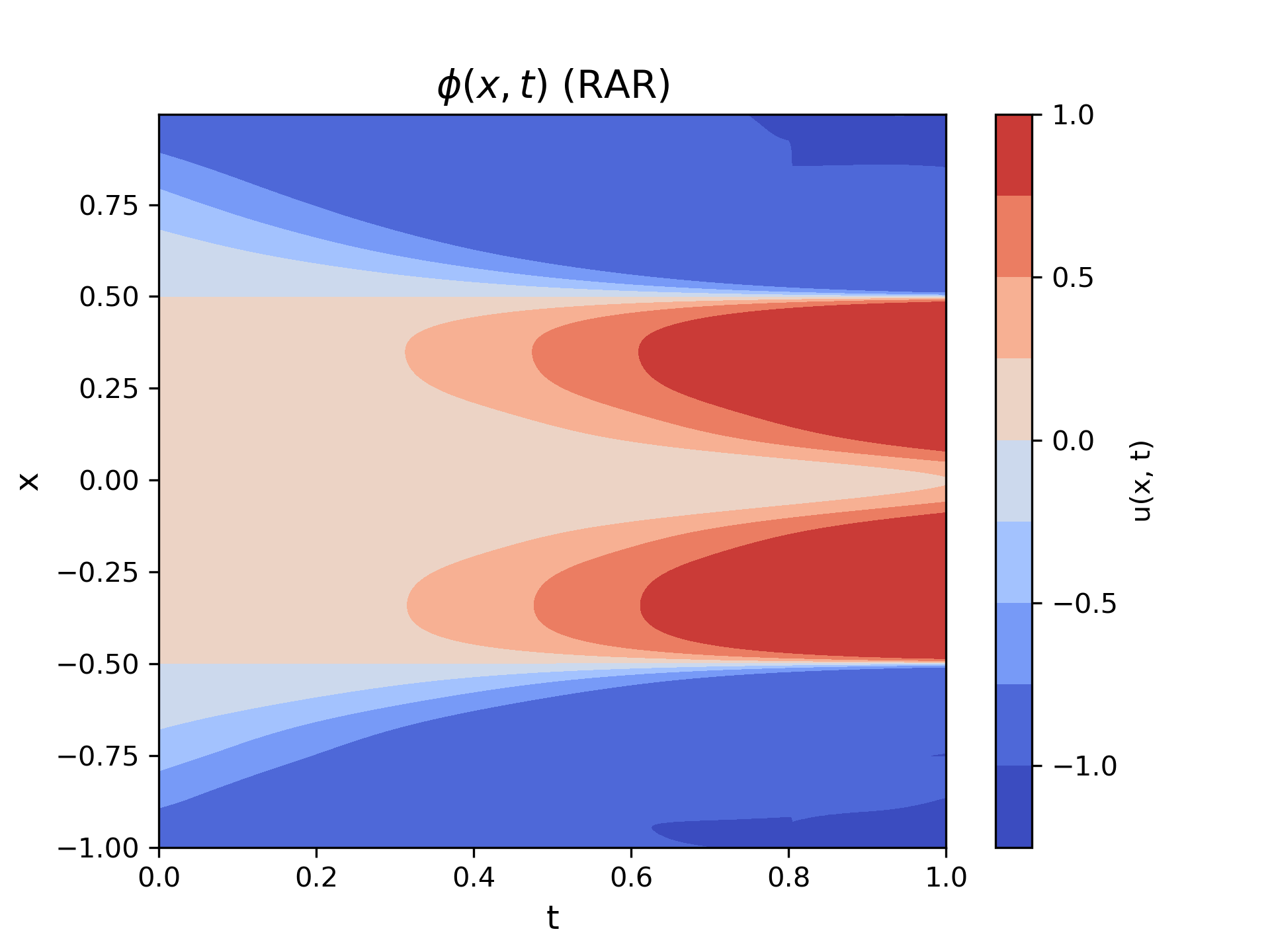}
    \end{minipage}
   \begin{minipage}[b]{0.3\textwidth}
     \includegraphics[scale=0.3]{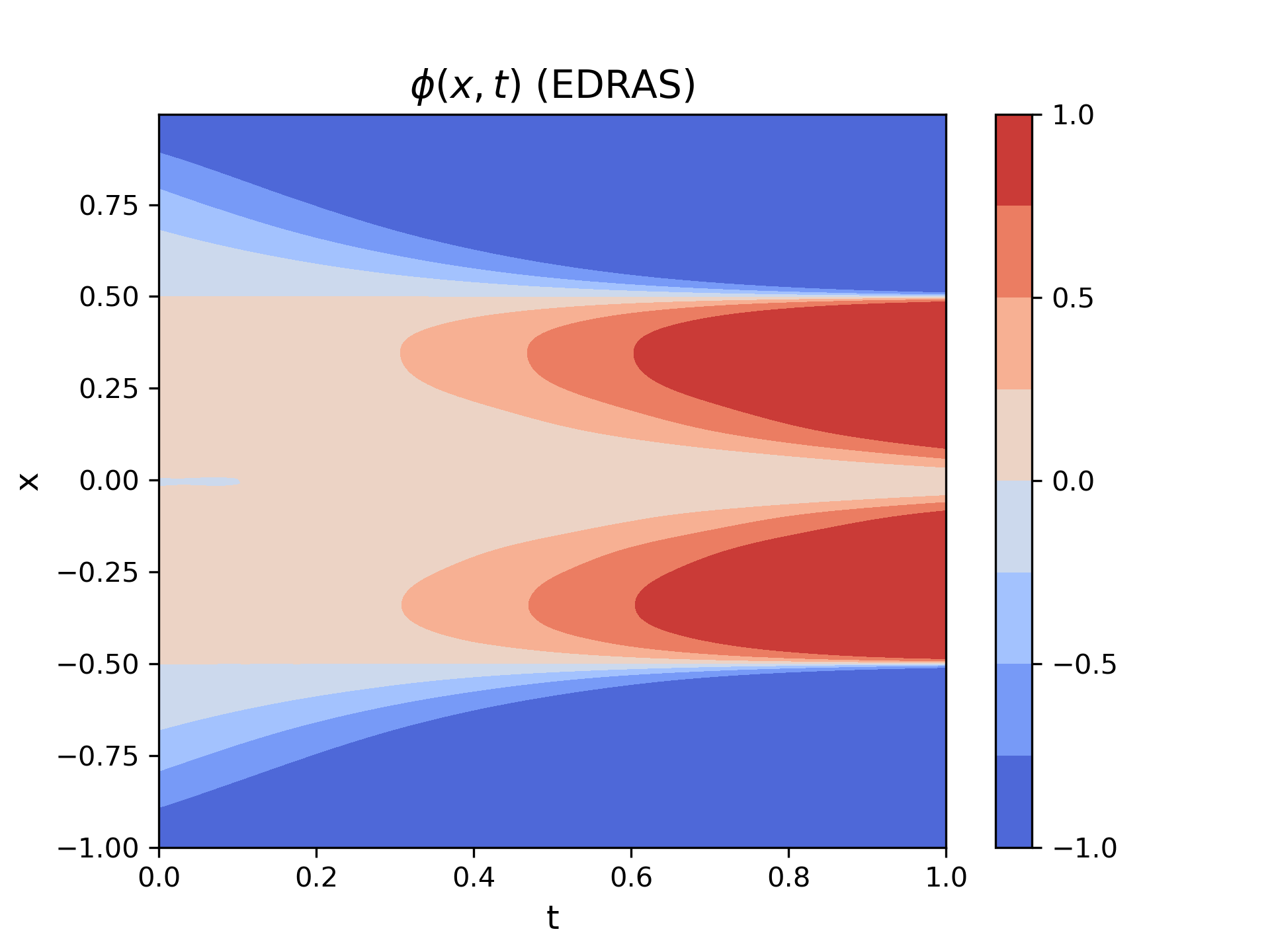}
   \end{minipage}
    \begin{minipage}[b]{0.3\textwidth}
        \includegraphics[scale=0.3]{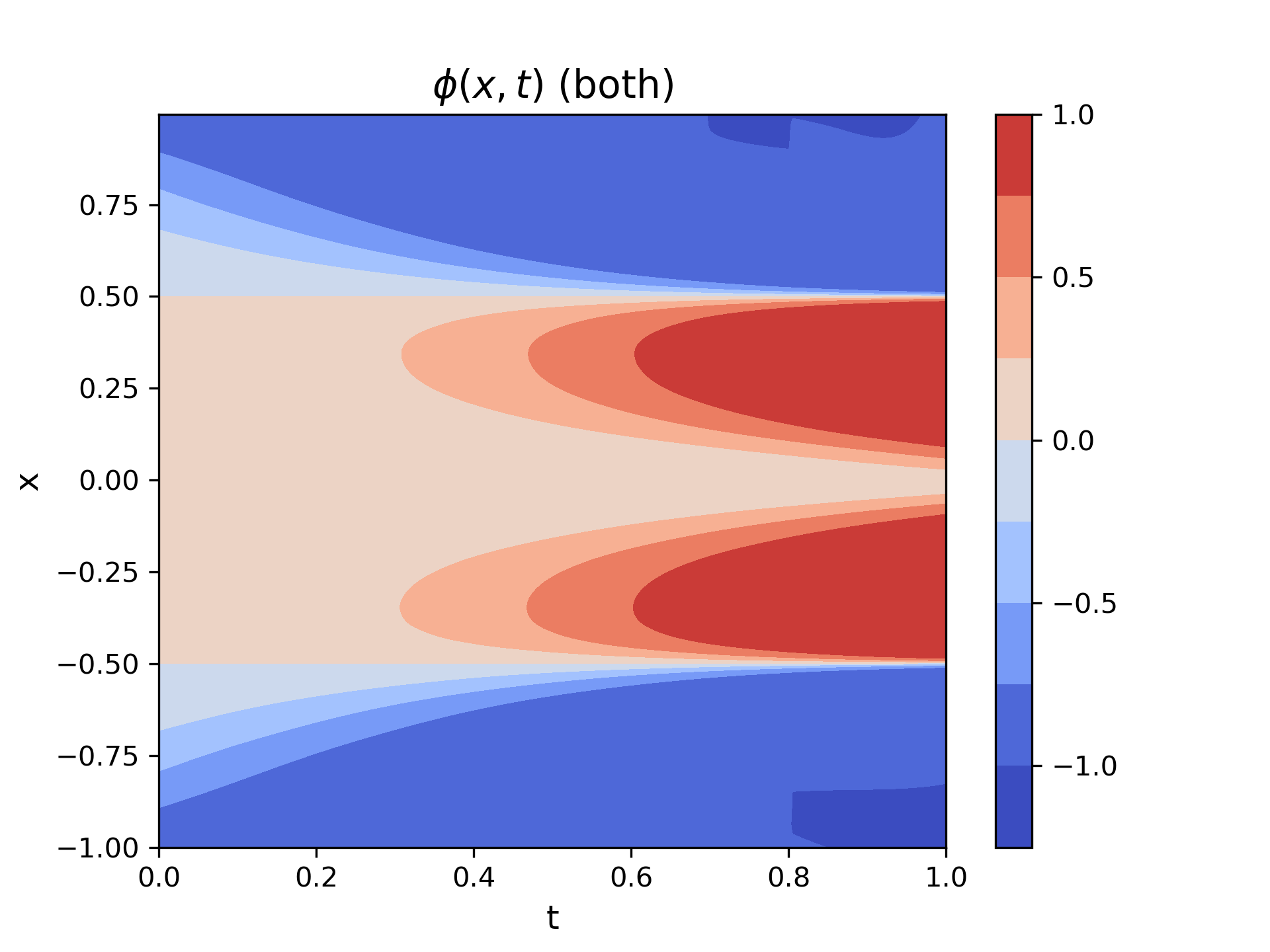}
  \end{minipage}
   \begin{minipage}[b]{0.3\textwidth}
        \includegraphics[scale=0.3]{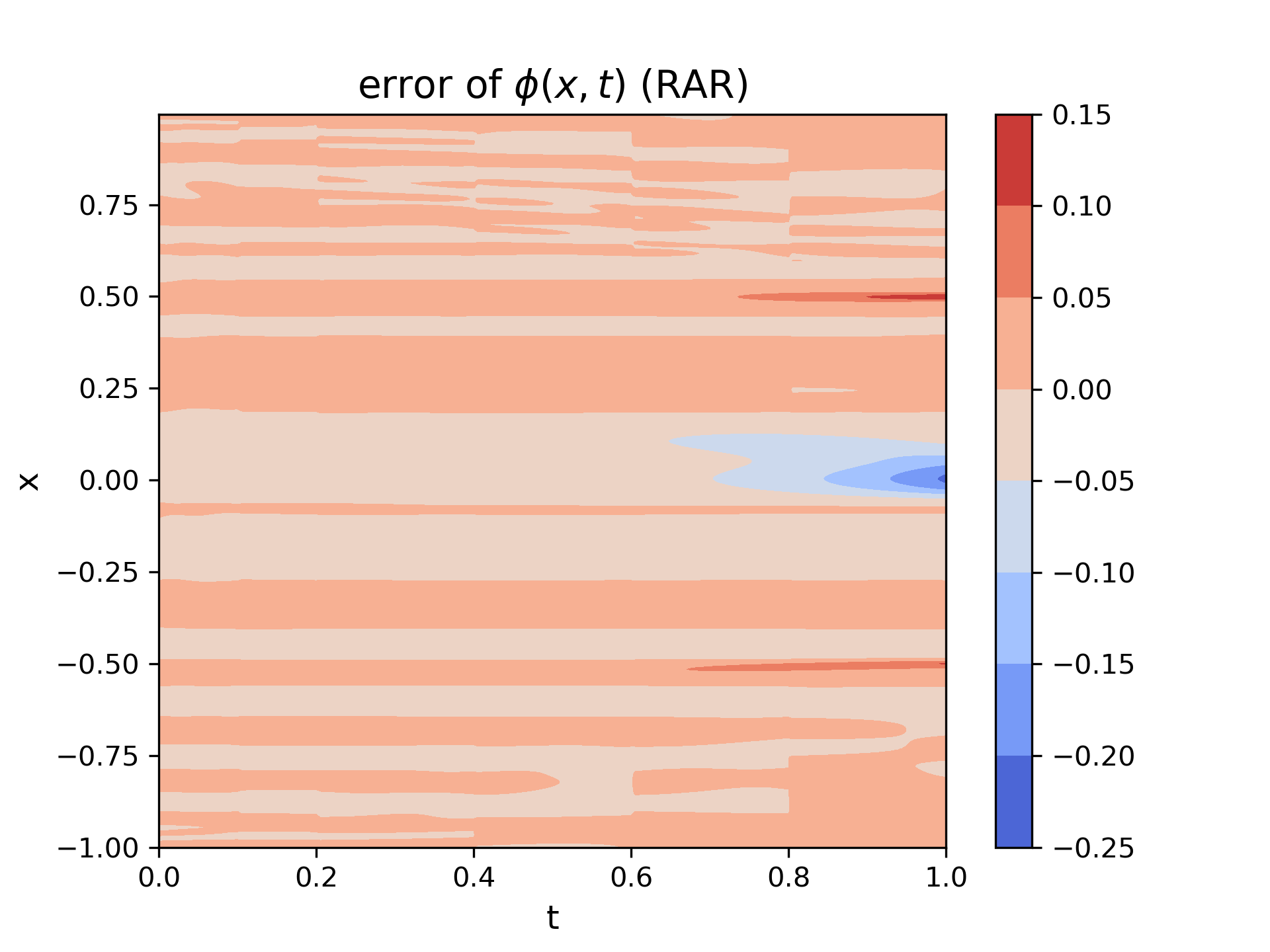}
    \end{minipage}
    \begin{minipage}[b]{0.3\textwidth}
       \includegraphics[scale=0.3]{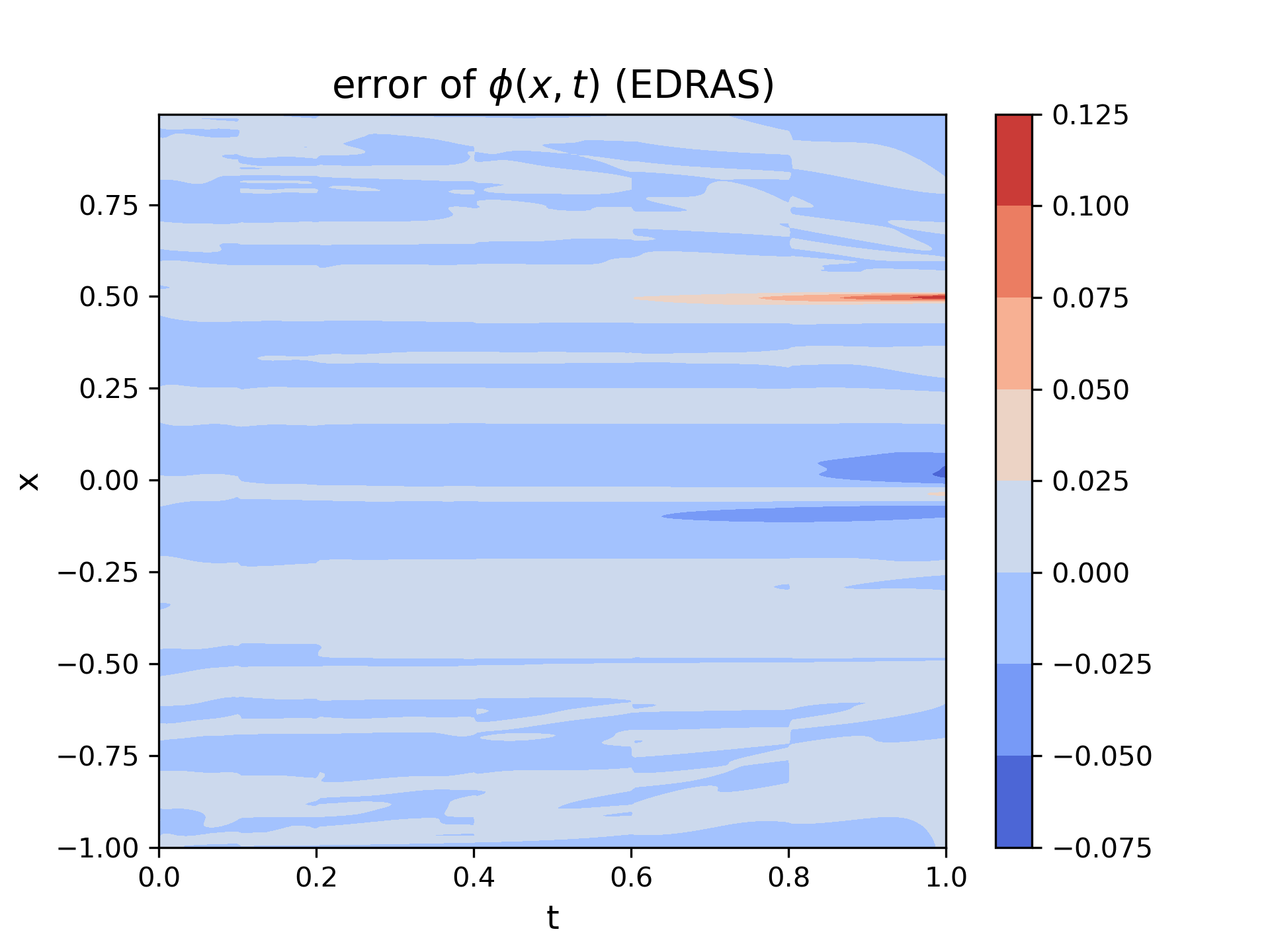}
     \end{minipage}
   \begin{minipage}[b]{0.3\textwidth}
       \includegraphics[scale=0.3]{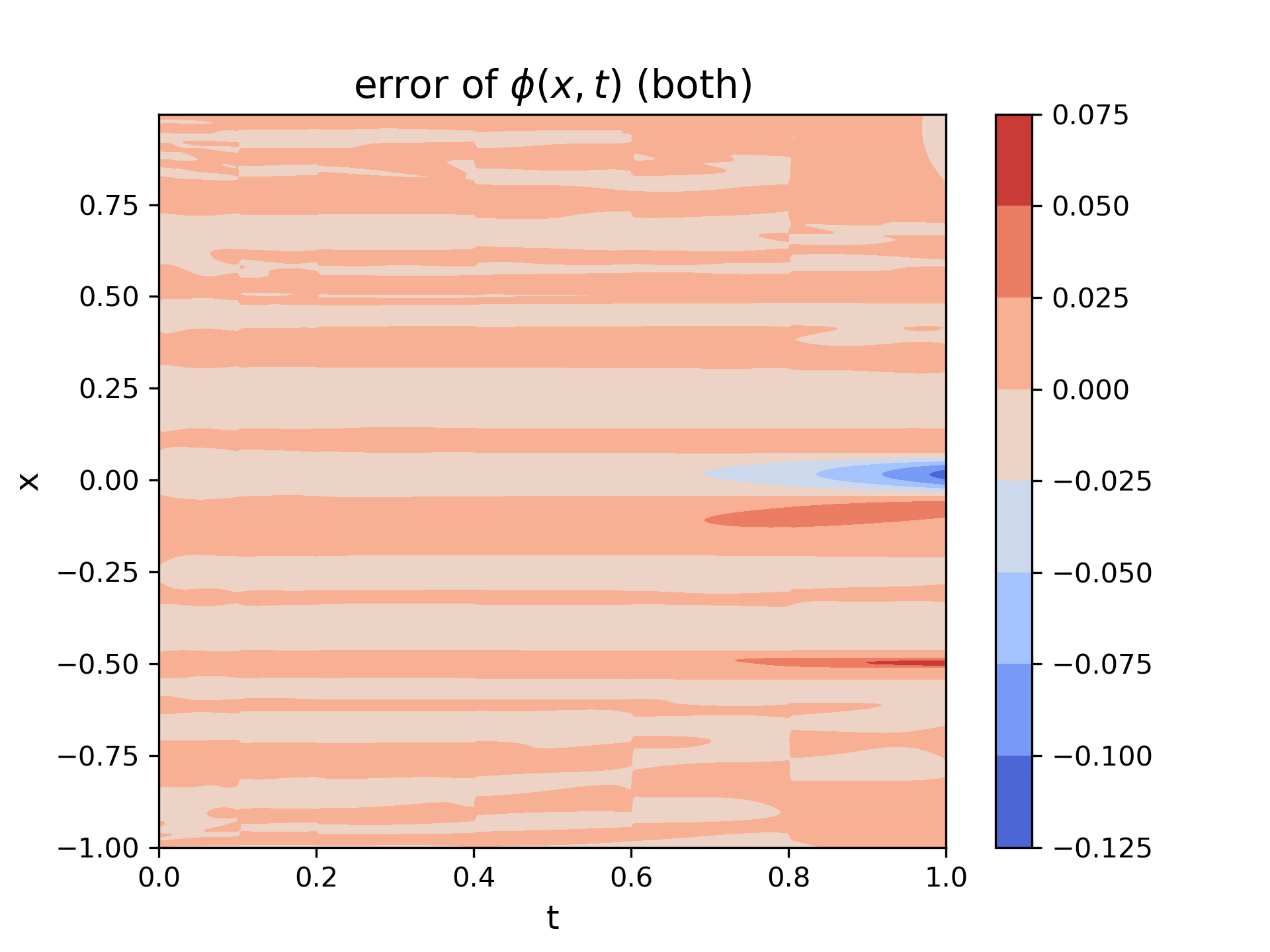}
   \end{minipage}
     \caption{Top: contour plots of solutions obtained using PINNs with RAR (left), EDRAS (middle), and RAR+EDRAS (right). Bottom: the errors between the FDM solution and the PINN solutions. The slightly lower average accuracy of the combined approach, compared to the EDRAS method alone, is due to the fact that a half of the resampled points in the combined method are chosen based on the EDRAS method and the other half based on the RAR method.}
   \label{fig:err_plot_ac1d}
\end{figure}

Our comparative analysis of EDRAS and RAR reveals EDRAS's superior performance. Figure \ref{fig:err_plot_ac1d} presents contour plots of the solutions for the Allen-Cahn equation in 1D, along with the errors between the reference solution and those obtained using the RAR method, the EDRAS method, and their combination, respectively. The reference solution is computed using the finite difference method. Quantitative metrics in Table \ref{tab:err-report-ac1d} demonstrate the effectiveness of EDRAS, achieving a sixfold reduction in relative mean squared error (MSE) and a twofold decrease in both relative mean absolute error (MAE) and relative $L^\infty$ error compared to the RAR method. These results show that as long as the EDRAS approach is incorporated, even partially, the performance of PINNs improves significantly.

\begin{table}[ht]
\centering
 \caption{Errors for EDRAS, RAR, and Combination of both for Allen-Cahn equation in 1D.}
\label{tab:err-report-ac1d}
    \begin{tabular}{cccc}
    \hline
\textbf{Metric} & \textbf{EDRAS} & \textbf{RAR} & \textbf{Combination of both} \\
    \hline
    MSE & $5.43 \times 10^{-5}$ & $3.05 \times 10^{-4}$ & $6.25 \times 10^{-5}$ \\
    Relative MSE & $1.08 \times 10^{-4}$ & $6.06 \times 10^{-4}$ & $1.24 \times 10^{-4}$ \\
    MAE & 0.0030 & 0.0064 & 0.0029 \\
    Relative Mean Abs. Error & 0.0051 & 0.0107 & 0.0049 \\
    Max for Absolute Error & 0.1127 & 0.2106 & 0.1108 \\
    Relative $L_{\infty}$ Error & 0.1127 & 0.2106 & 0.1108 \\
    \hline
    \end{tabular}
\end{table}

 \begin{figure}[ht]
     \centering
       \begin{minipage}[b]{1\textwidth}
         \includegraphics[scale=0.4]{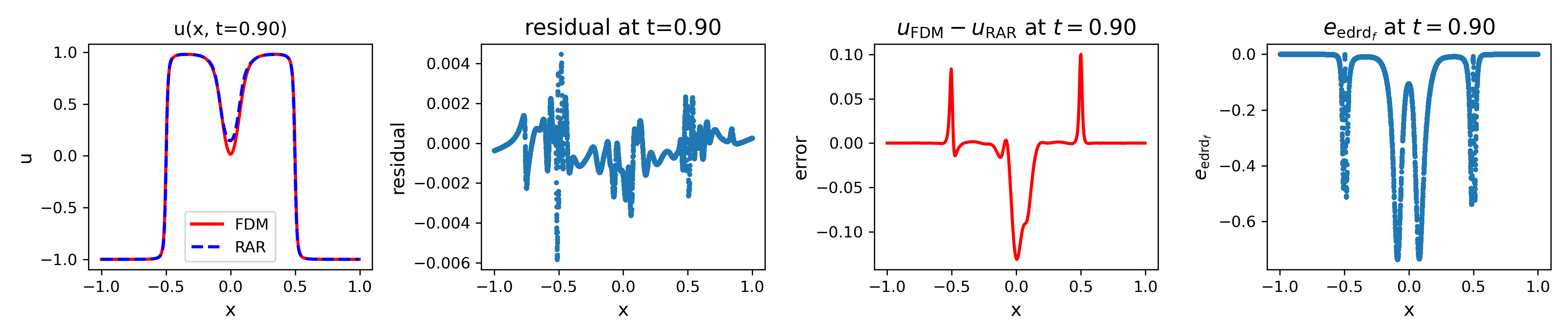}
     \end{minipage}
     \begin{minipage}[b]{1\textwidth}
         \includegraphics[scale=0.4]{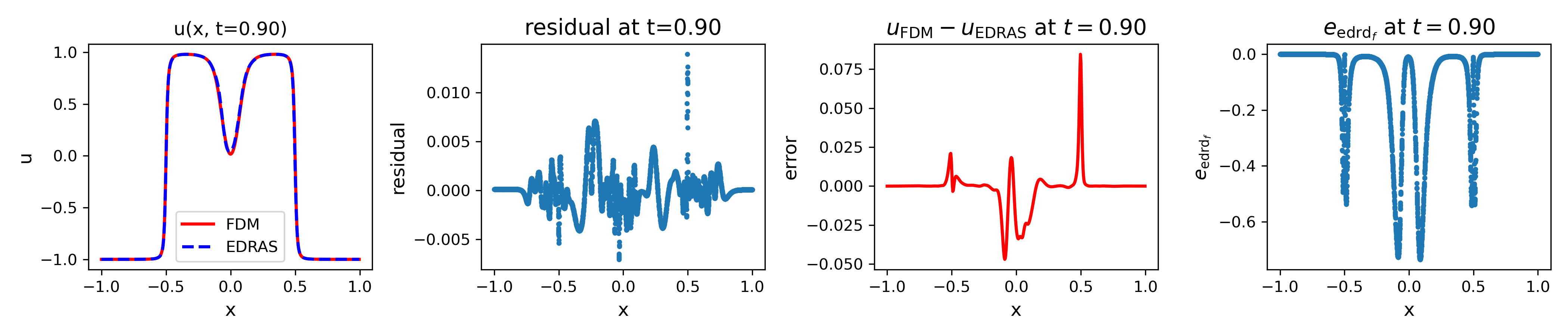}
     \end{minipage}
     \begin{minipage}[b]{1\textwidth}
         \includegraphics[scale=0.4]{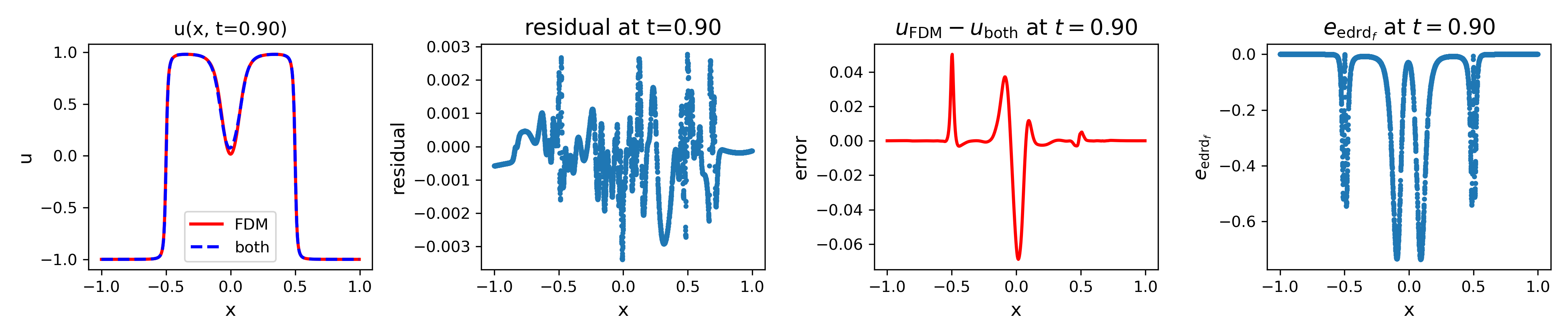}
     \end{minipage}
      \caption{Snapshots of solutions obtained from PINN  with RAR (top), EDRAS (middle) and both EDRAS+RAR (bottom) at $t=0.9$, respectively. So long as EDRAS is used, even partially, the result is improved noticeably, evidenced by the combined case.}
     \label{fig:ac1d-0.9}
\end{figure}

\begin{figure}[ht]
      \centering
       \begin{minipage}[b]{0.3\textwidth}
         \includegraphics[scale=0.2]{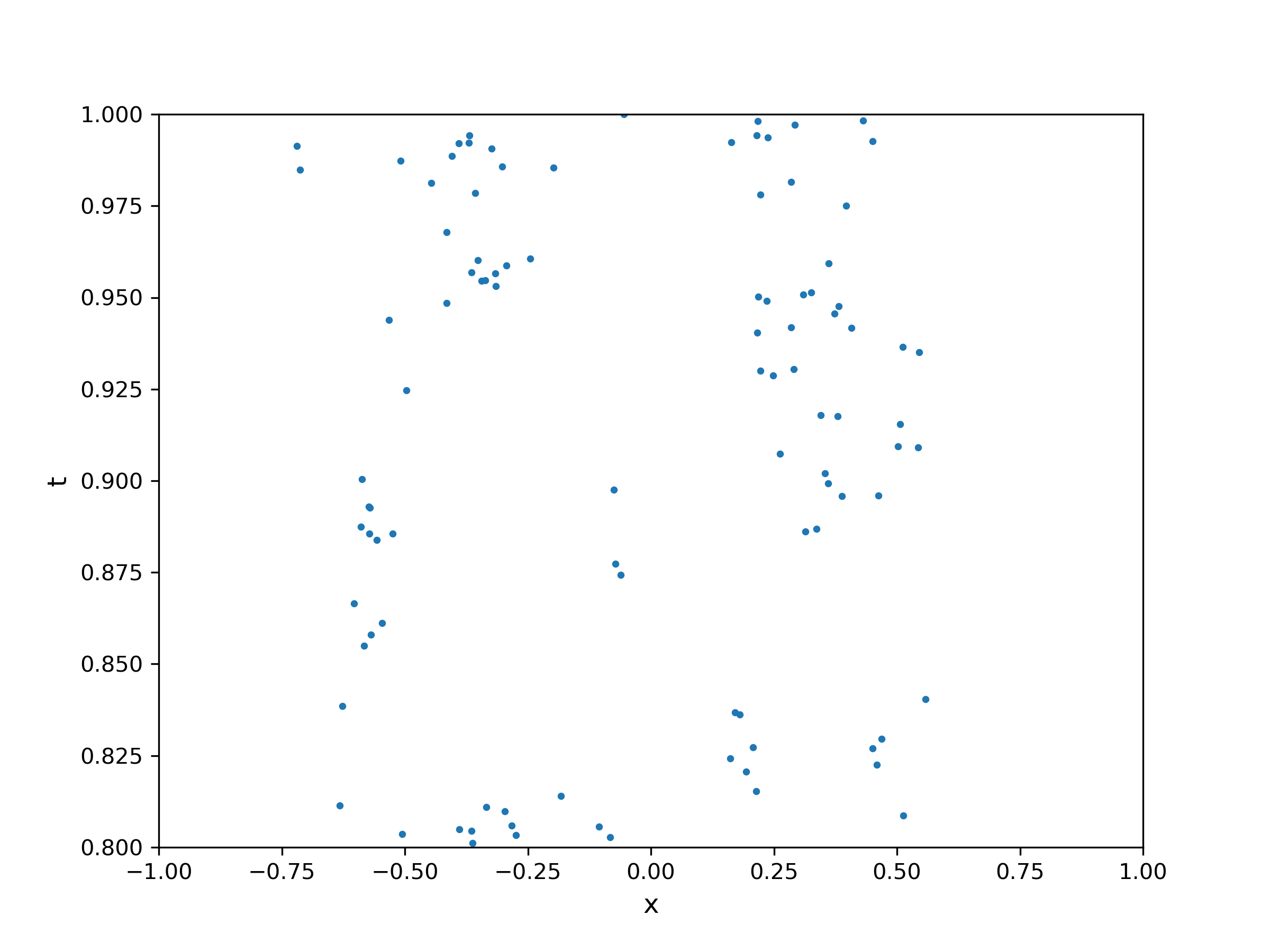}
     \end{minipage}
     \begin{minipage}[b]{0.3\textwidth}
         \includegraphics[scale=0.2]{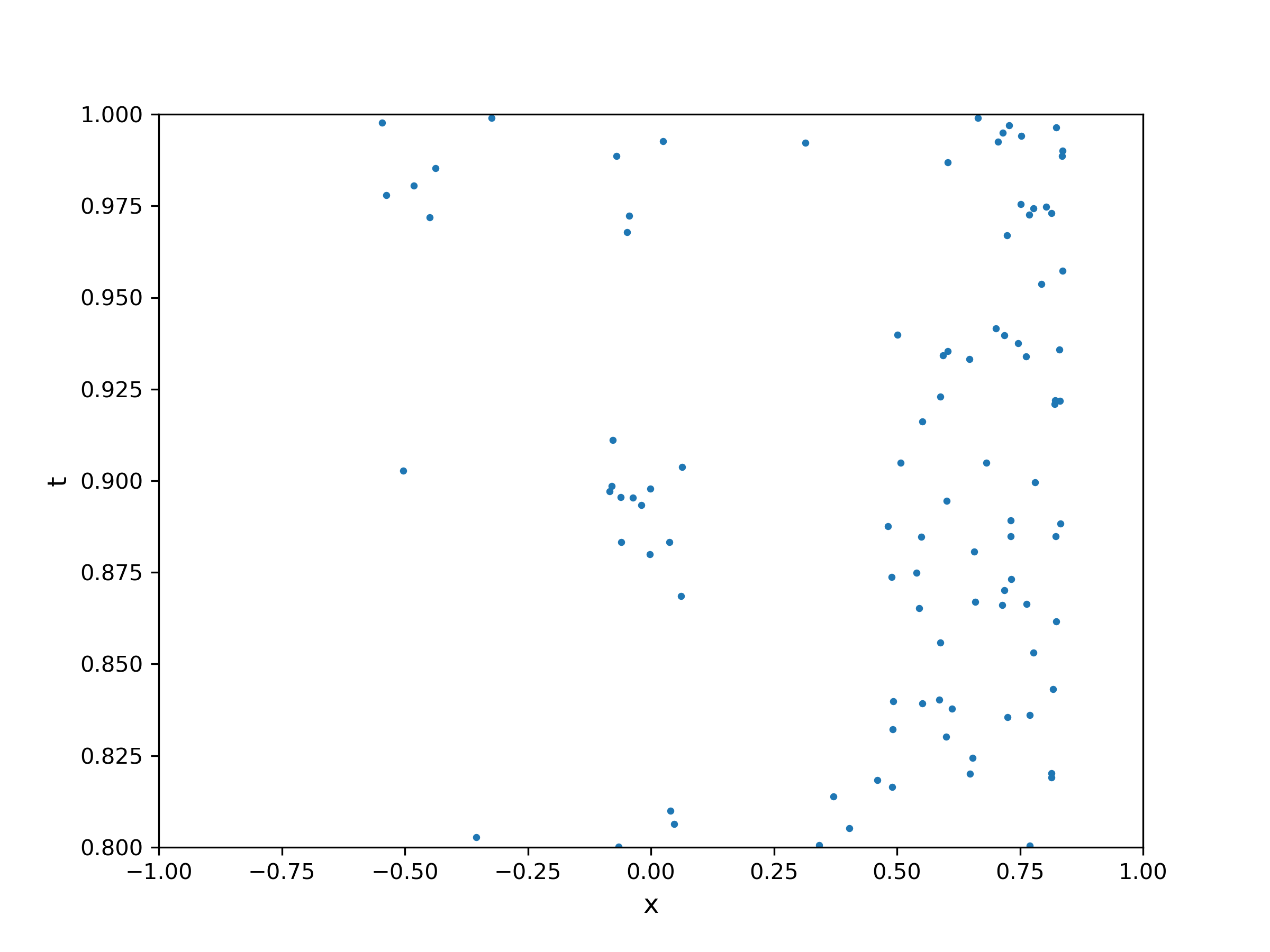}
     \end{minipage}
     \begin{minipage}[b]{0.3\textwidth}
         \includegraphics[scale=0.2]{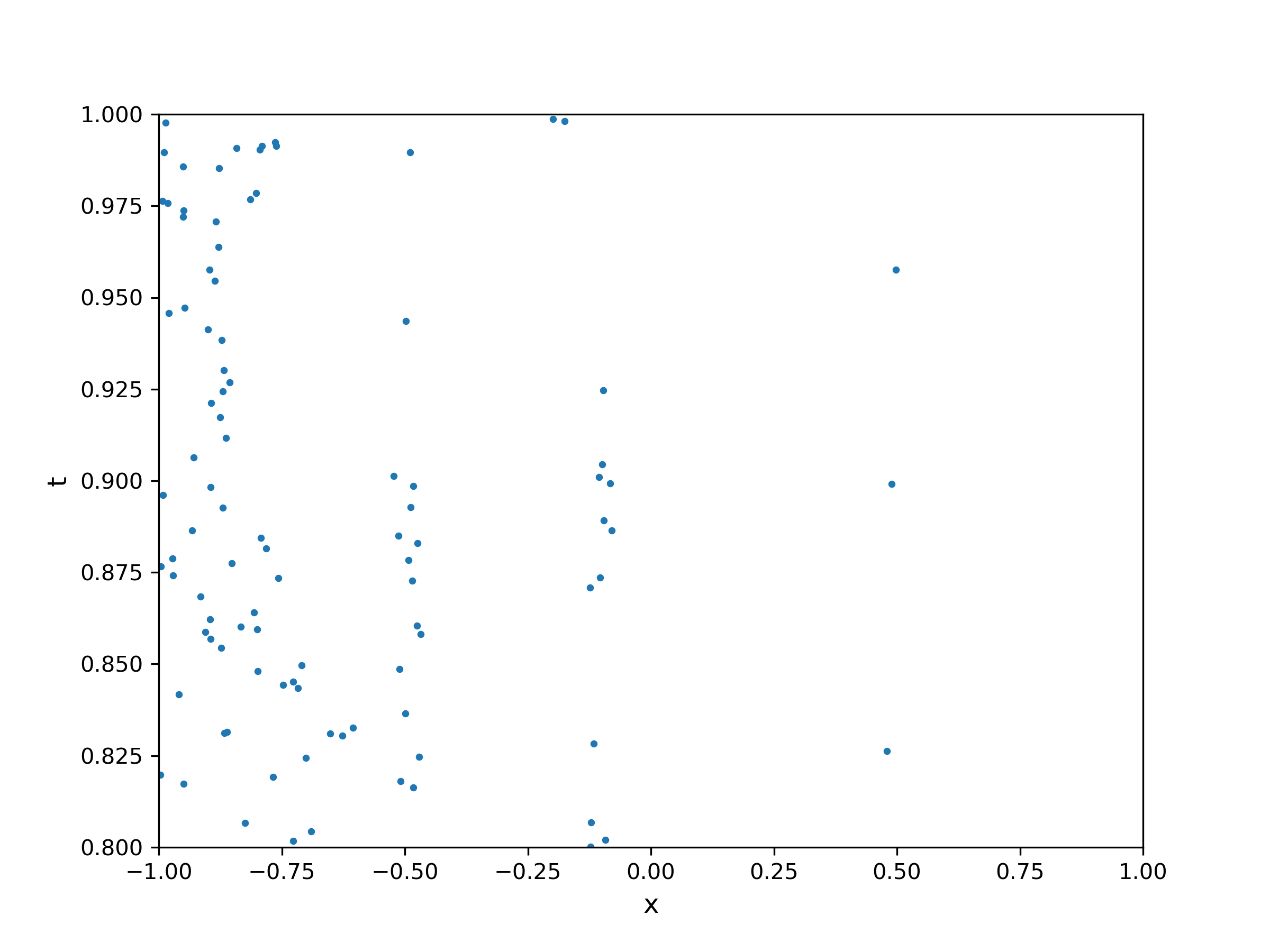}
     \end{minipage}

     \begin{minipage}[b]{0.3\textwidth}
         \includegraphics[scale=0.2]{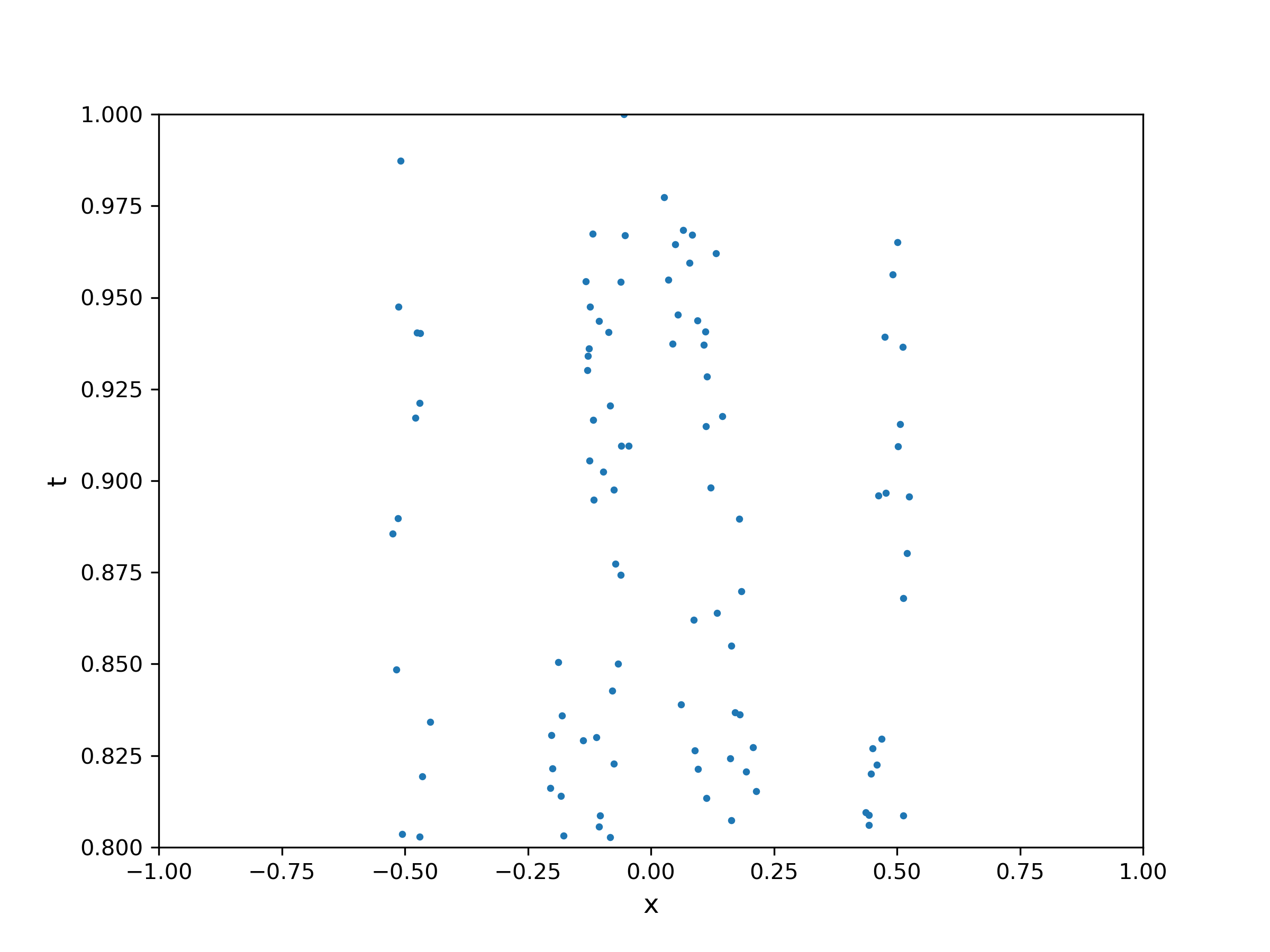}
     \end{minipage}
      \begin{minipage}[b]{0.3\textwidth}
         \includegraphics[scale=0.2]{fig/ac1d-EDRR/integrand_resample_100_0.8000_1000.png}
     \end{minipage}
     \begin{minipage}[b]{0.3\textwidth}
         \includegraphics[scale=0.2]{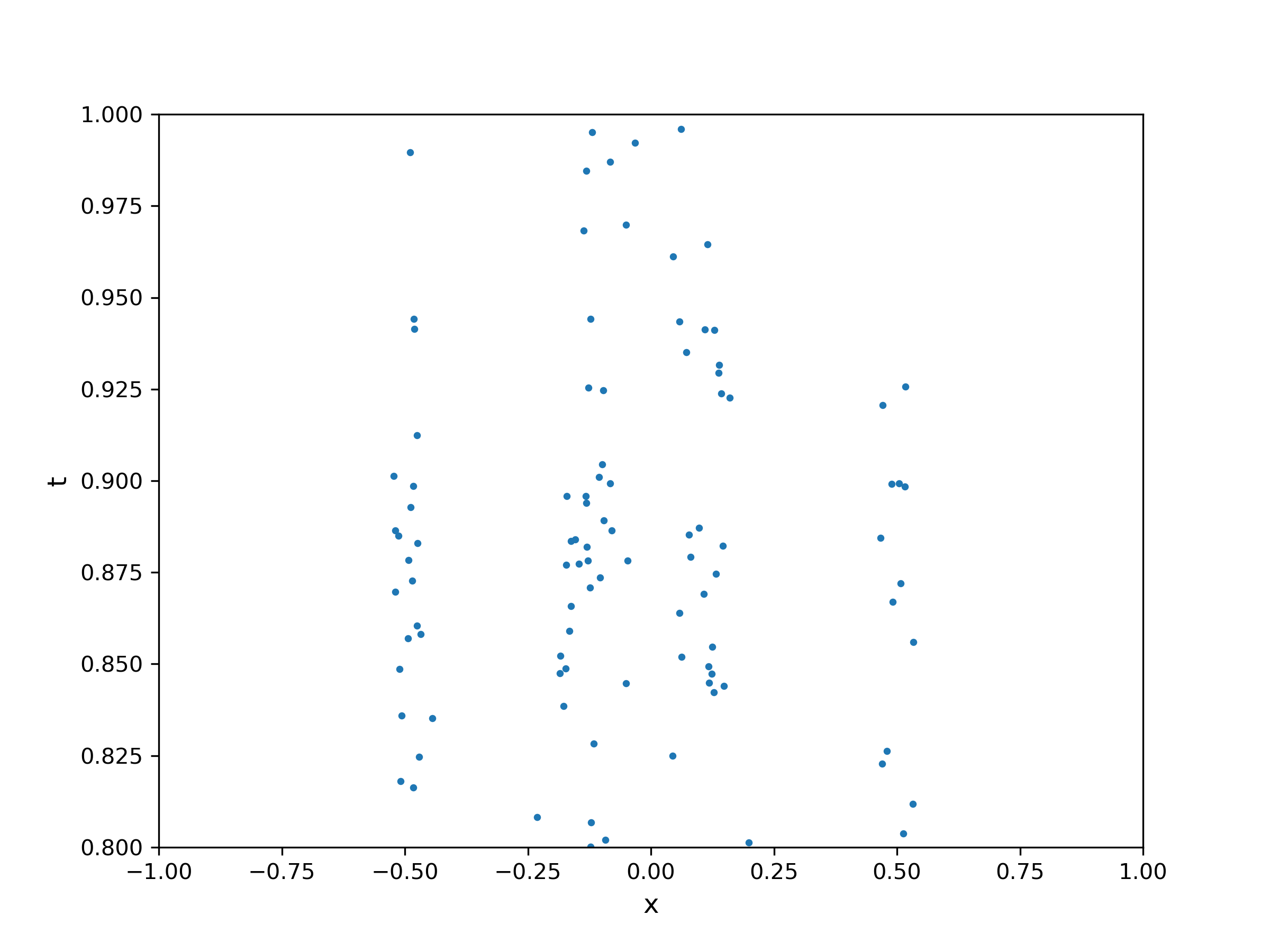}
     \end{minipage}
     \caption{Top: 3 examples of re-sampled collocation points at $[0.8, 1.0]\times[-1, 1]$ using RAR. Bottom:  3 examples of re-sampled points at $[0.8, 1.0]\times[-1, 1]$ using EDRAS. The resampled collocation points strongly correlate with energy dissipation rate fluctuations in EDRAS.  }
     \label{fig:ac1d-sample}
 \end{figure}

The temporal evolution of solutions is illustrated through snapshots taken at $t = 0.15, 0.3, 0.5, 0.7$ (in the appendix) and $t = 0.9$, respectively, shown in Figure \ref{fig:ac1d-0.9}. These results compare the RAR method (top), the EDRAS method (middle), and their combination (bottom). Our analysis reveals a strong correlation between significant errors and elevated fluctuations in the energy dissipation rate density. The EDRAS method demonstrates a more precise alignment between observed errors and energy dissipation rate density fluctuations compared to the RAR method, as illustrated in Figure \ref{fig:ac1d-0.9}.

Figure \ref{fig:ac1d-sample} compares the distribution of the resampled points in the temporal-spatial subdomain, $[0.8, 1.0] \times [-1, 1]$,  for the RAR (top row) and EDRAS (bottom row) method, respectively. The EDRAS method demonstrates a clear alignment between the resampled points and the energy dissipation rate density, with finer sampling in regions of higher local solution variance. In contrast, the RAR method concentrates points primarily in areas with steep gradients, where residuals are large, while under-sampling regions with large approximation errors but with relatively smaller residuals, such as areas near $x=0$—a key limitation of residual-based approaches. Notably, the EDRAS method’s pronounced clustering of points near $x=0$ visually underscores its advantage over RAR in adaptive sampling efficiency.

\begin{figure}[ht]
       \centering
        \includegraphics[scale=0.7]{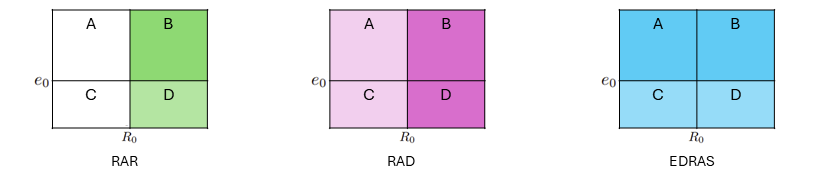}
        \caption{Schematics of the RAR (left), RAD (middle), and EDRAS (right) method, respectively. The colored regions represent areas identifiable by the adaptive sampling method, with darker shades indicating higher sampling probabilities.}
       \label{fig:prob}
\end{figure}

This observation motivates us to examine RAR and RAD further from a probability perspective.  Let $e_0>0$ be the threshold for approximation error which helps to exactly locate the important areas where more collocation points are needed. Since we don't have a reference solution available in practice, so one usually choose the residual as an alternative choice to identify the important areas. Let $R_0>0$ be a threshold value for the residual. The areas with residuals exceeding $R_0$ should sample more points into the training set due to the large residual values.  Let $y=(x,t)$ and $R(y)=R(x, t)$ for simplicity. We classify all collocation points in the sampled temporal-spatial domain into four groups, as demonstrated in Figure \ref{fig:prob},
\begin{itemize}
    \item Group A: A small residual and large approximation error. Namely, $A=\{y\in [-1,1]\times [0, 1]\|R(y)\leq R_0 \quad \& \quad e(y)>e_0\}$.
    \item Group B: A large residual and large approximation error. Namely, $B=\{y\in [-1,1]\times [0, 1]\|R(y)> R_0 \quad \& \quad e(y)> e_0\}$.
    \item Group C: A small residual and small approximation error. Namely, $C=\{y\in [-1,1]\times [0, 1]\| R(y) \leq R_{0} \quad \& \quad e(y)\leq e_0\}$.
    \item Group D: A large residual and small approximation error. Namely, $D=\{y\in [-1,1]\times [0, 1]\|R(y) > R_{0} \quad \& \quad e(y)\leq e_0\}$.
\end{itemize}

We call the collocation points in group A and group B the critical/important collocation points due to the large approximation error in the groups. Hence, during the adaptive sampling process, we expect $p(A\cup B)$, the probability of selected points that belong to $A$ and $B$, as large as possible so that we can effectively select the important points.  Next, we examine disadvantages and advantages of the three adaptive sampling methods: RAR, RAD and EDRAS, in the probability perspective.

The RAR method selects the top $m$ points with largest residuals from a large pool of $n$ dense candidates. The RAD method (residual-based adaptive distribution) samples $m$ collocation points based on the probability distribution defined by
\begin{equation}\label{eq:prob_rad}
    p(y) = \frac{R(y)}{\int_{\Omega} R(y)dy}.
\end{equation}
In practice,  $\int_{\Omega} R(y)$ is estimated by $\sum_{i=1}^n R(y_i)$ where $n >>m$.

Let $n_S$ be the number of selected points in group $S$, where $S\in \{A, B, C, D\}$, among the $m$ added points. Now we calculate the probabilities for  these four groups  with respect to the RAR or RAD method.
For the RAR method,
\begin{equation}\label{eq:rar_p}
p_{RAR}(S)=\frac{n_S}{m}.
\end{equation}
For the RAD method,
\begin{equation}\label{eq:rad_p}
 p_{RAD}(S) = \sum_{y_j\in S}p(y_j).
\end{equation}
where $y_j$ is sampled from distribution \eqref{eq:prob_rad}.

Then, based on the construction of the RAR method and RAD method, we conclude that when $m$ is not large enough, the RAR method can only select points from $B\cup D$ while the RAD method has a positive probability to select points from $A\cup C$.  Since the points in $B\cup D$ have higher residuals than those in $A\cup C$, they are more likely to be sampled. Unlike RAR, since RAD does not strictly prioritize selecting the point with the very top residual, some relatively lower-residual points in $A$ may still be sampled. Hence, RAD's sampling is less aggressive than RAR's top$-m$
 selection in $B\cup D$, as demonstrated in Figure \ref{fig:prob}. Hence, we conclude
$$p_{RAR}(B\cup D) =1, \quad p_{RAD}(B\cup D) < 1, \quad \text{as $m$ is not large enough.} $$
This explains the drawback of the RAR method (that is $p_{RAR}(A)=0$) and why the RAD method, can alleviate this issue ($p_{RAD}(A)>0$ as long as $A\neq \emptyset$).

We numerically validate the above discussion on the RAR and RAD method and also compare them with the proposed EDRAS method in terms of probability of these 4 groups on the temporal-spatial domain, $[0, 1]\times [-1, 1]$. The results in domain $[0.8, 1]\times [-1, 1]$ are shown in Figure \ref{fig:prob_view} and the results on other time intervals are given in the appendix.  We set $e_0=0.001$ and $R_0=\frac{1}{n}\sum_{i=1}^n R(y_i)$, the mean of the residuals at all points from the dense set that the important points are sampled. We select $m=100$ collocation points using these 3 methods from the temporal spatial domain, respectively. Then, we compute the corresponding probabilities for the 4 groups by \eqref{eq:rar_p} and \eqref{eq:rad_p}.  For the RAD method, we repeat this sampling process 100 times to ensure the stability and reliability of the estimation.  For the EDRAS method, we select $m=100$ collocation points with the top $m$ large energy dissipation rate density values. We implement the adaptive sampling strategy every 100 epochs during the training process and then, we calculate these probabilities of the 4 groups during each resampling process. The probabilities of the 4 groups with respect to epochs for these 3 methods are shown in Figure \ref{fig:prob_view}. We notice clearly that the numerical results are consistent with our previous analysis about the RAR method ($p_{RAR}(A)=0$) and the RAD method ($p_{RAD}(A)>0$). The numerical results also demonstrate why the EDRAS method works better than the RAR method as shown in Figure \ref{fig:prob_view} as the EDRAS method can cover both groups, $A$ and $B$, and even have a higher probability for $(A\cup B)$ compared with the RAR method as well as the RAD method.
\begin{itemize}
\item  Comparing the EDRAS method with the RAR method numerically, we find $p_{EDRAS}(A\cup B)>p_{RAR}(A\cup B).$
\item  Comparing the EDRAS method with the RAD method numerically, we find $p_{EDRAS}(A\cup B)>p_{RAD}(A\cup B).$
\end{itemize}

We observe that the RAR method focuses mainly on group A and B in practice and the RAD method also covers group C and D. The EDRAS method covers all 4 groups.  These observations align with the above analysis and discussion and demonstrate clearly the advantage of the EDRAS method compared with the RAR and RAD method. Moreover, compared with the RAD method, the EDRAS method is more computationally efficient and thus suitable for high-dimensional PDEs since the estimation of the normalization factor $\int_{\Omega}R(y)dy$ is computationally costly and  even intractable as the dimension increases.

\begin{figure}[ht]
      \centering
         \includegraphics[scale=0.3]{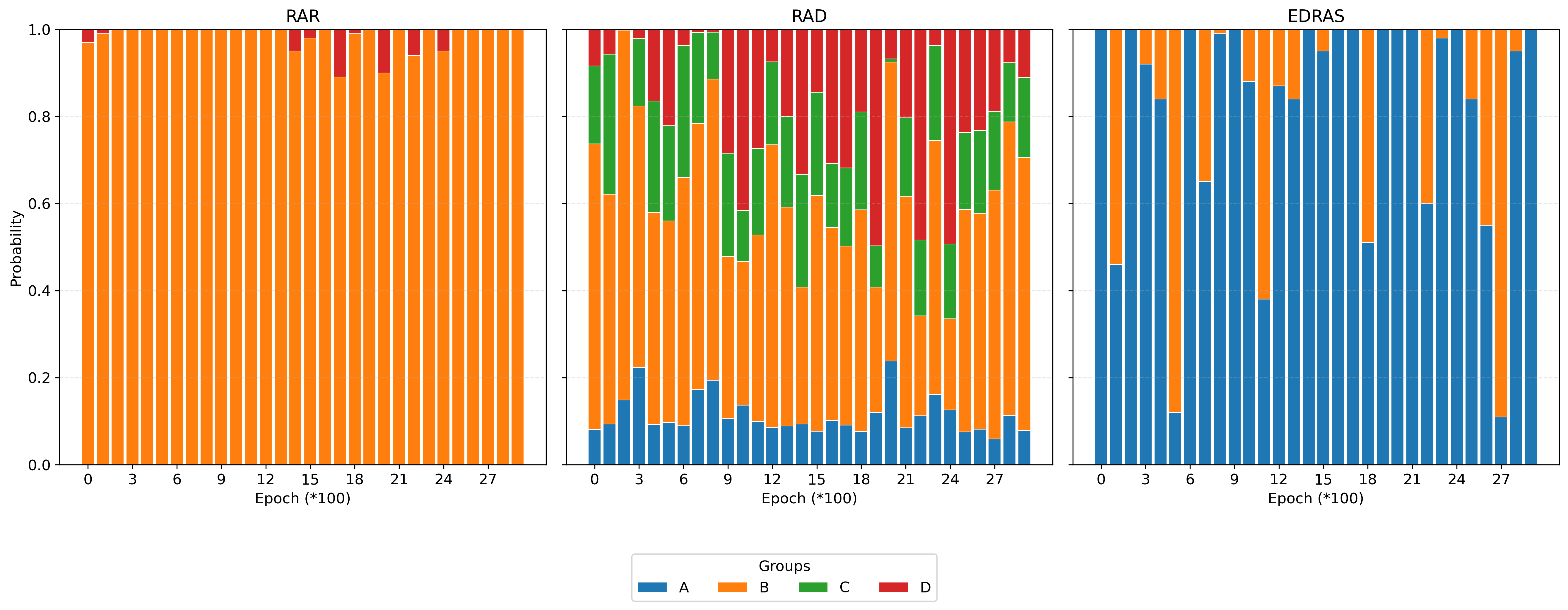}
      \caption{The sampling probabilities of the 4 groups, $(A,B,C,D)$, for the RAR, RAD and EDRAS method in domain $[0.8, 1.0]\times [-1,1]$ plotted over training epochs, respectively. We observe that RAR samples primarily based on the residual, EDRAS focuses exclusively on the error, and RAD samples from all groups. }
 \label{fig:prob_view}
\end{figure}

\subsection{Static Boundary Conditions  Versus Dynamic Boundary Conditions in 2D}\label{subsec:AC_neu}

We consider the free energy consisting of the bulk free energy and surface free energy given by
\begin{equation}\label{eq:double-well}
    E[t] = \int_{\Omega} \frac{\epsilon^2}{2}|\nabla \phi|^2+\frac{1}{4}(\phi^2-1)^2 \D\bx + \int_{\p \Omega} \frac{\epsilon_s^2}{2}|\nabla_s \phi|^2+\frac{1}{4}(\phi^2-1)^2 \D S,
\end{equation}
where $\epsilon$ and $ \epsilon_s$ measure the width of the diffuse interface layers in the bulk and on the boundary between two distinct phases, respectively. The Allen-Cahn equation with the Allen-Cahn type dynamic boundary condition is derived by the generalized Onsager principle\cite{xiaoboentropy, xiaobocms,  chunyan2023thesis} as follows
\begin{equation}\label{eq:ac_2d}
\begin{cases}
      \phi_t = -M_b(\phi^3-\phi-\epsilon^2 \nabla^2 \phi),  \quad \bx\in \Omega, \quad t\in (0, 1], \\
      \phi_t = -M_s( \phi^3-\phi-\epsilon_s^2 \nabla_s^2 \phi + \epsilon^2 \bn\cdot \nabla \phi), \quad \bx \in \p \Omega, \quad t\in (0, 1].
\end{cases}
\end{equation}
The energy dissipation rate of the system is given by
\begin{equation}
    \begin{split}
      \frac{dE}{dt}&=\int_\Omega \left(-\epsilon^2\nabla^2 \phi + \phi^3-\phi\right)\phi_t\D\bx + \int_{\p\Omega} \left(-\epsilon_s^2\nabla_s^2 \phi + \epsilon^2 \bn\cdot \nabla \phi + \phi^3-\phi \right)\phi_t\D S \\
      &= \int_\Omega -\frac{\phi^2_t}{M_b}\D\bx -\int_{\p\Omega} \frac{\phi^2_t}{M_s}\D S \leq 0.
    \end{split}
\end{equation}

In this section, we investigate the Allen-Cahn model under two distinct types of boundary conditions: dynamic and static (Neumann) boundary conditions, to highlight their fundamental differences. The dynamic boundary condition is given in \eqref{eq:ac_2d} while the static, Neumann boundary condition refers to the  zero flux at the boundary, $\bn\cdot \nabla \phi=0$, in the Allen-Cahn model without the surface free energy.
The Allen-Cahn equation with the static boundary condition is given as
\begin{equation}\label{eq:ac_2d_neumann}
\begin{cases}
      \phi_t = -M_b(\phi^3-\phi-\epsilon^2 \nabla^2 \phi),  \quad \bx\in \Omega, \quad t\in (0, 1], \\
      \bn\cdot \nabla \phi=0, \quad \bx \in \p \Omega, \quad t\in (0, 1].
\end{cases}
\end{equation}
Its energy dissipation rate is given by
\begin{equation}
      \frac{dE}{dt}=\int_\Omega \left(-\epsilon^2\nabla^2 \phi + \phi^3-\phi\right)\phi_t\D\bx = \int_\Omega -\frac{\phi^2_t}{M_b}\D\bx  \leq 0.
\end{equation}

The numerical implementation employs a time-marching strategy using a sequence of identically structured deep neural networks. Each DNN employs 3 hidden layers with 128 nodes per layer and hyperbolic tangent activation functions in each hidden layer, no activation function is used in the output layer. We adopt a uniform time step of $\Delta t=0.2$, with a refined initial segment $[0,0.01]$ to ensure accurate resolution of the initial condition—a critical requirement for solving PDEs using PINNs. The loss function weights are set to $w_r=1, w_s=1$ for residuals and surface terms, and $w_i=1000$ for initial conditions to enforce accurate initial state learning.
Our training protocol consists of:
\begin{itemize}
\item Initial sampling: 10000 initial points, 10000 residual points, and 3200 boundary points across the spatial-temporal domain.
\item Adaptive refinement: Progressive addition of residual and boundary points using EDRAS and RAR methods, either independently or in combination.
\item Optimization: Adam optimizer with learning rate 0.001 and batch size 2048.
\end{itemize}

\begin{figure}[H]
    \centering


    \begin{subfigure}[t]{1\textwidth}
        \centering
        \includegraphics[scale=0.21]{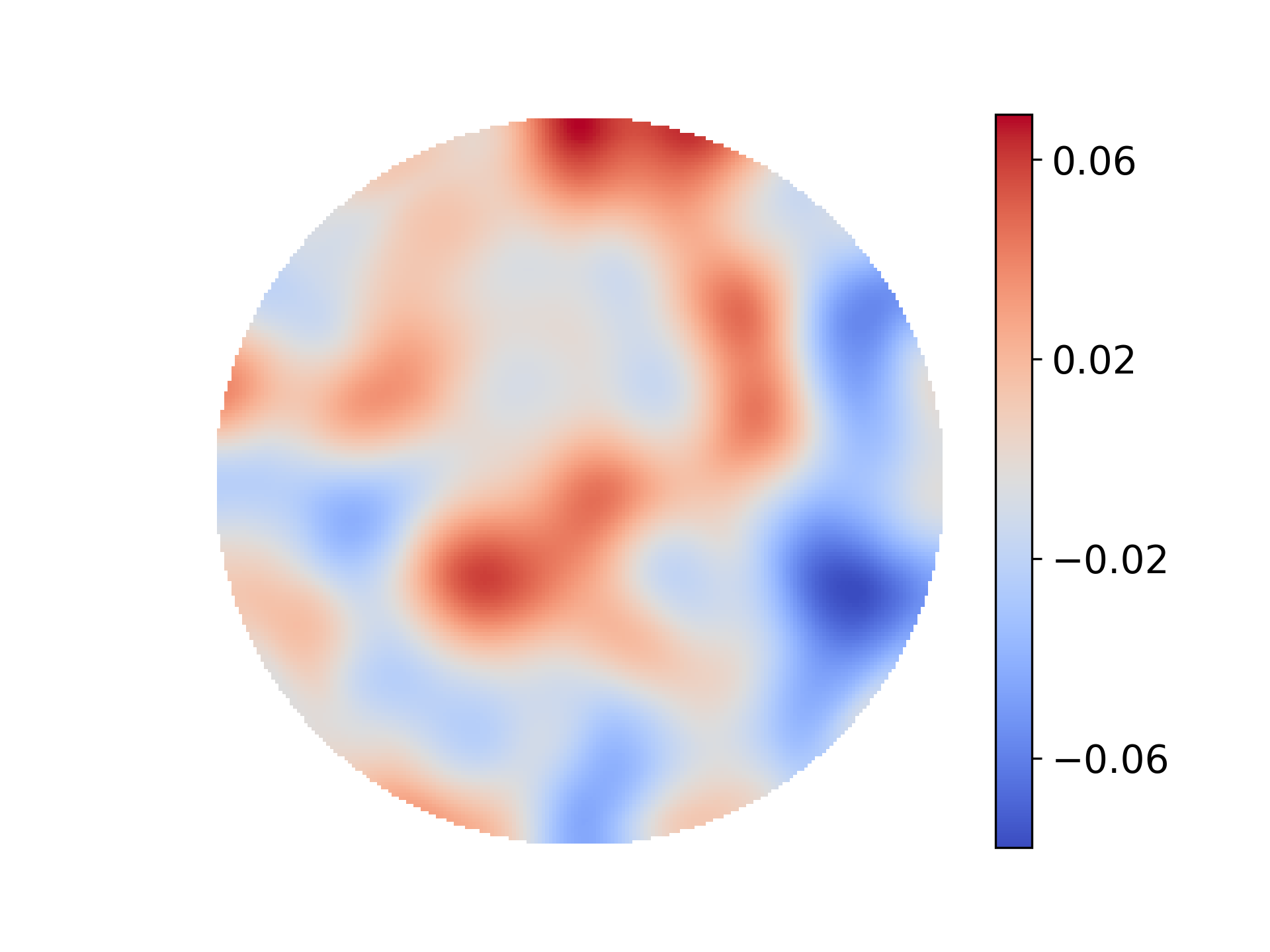}
        \includegraphics[scale=0.21]{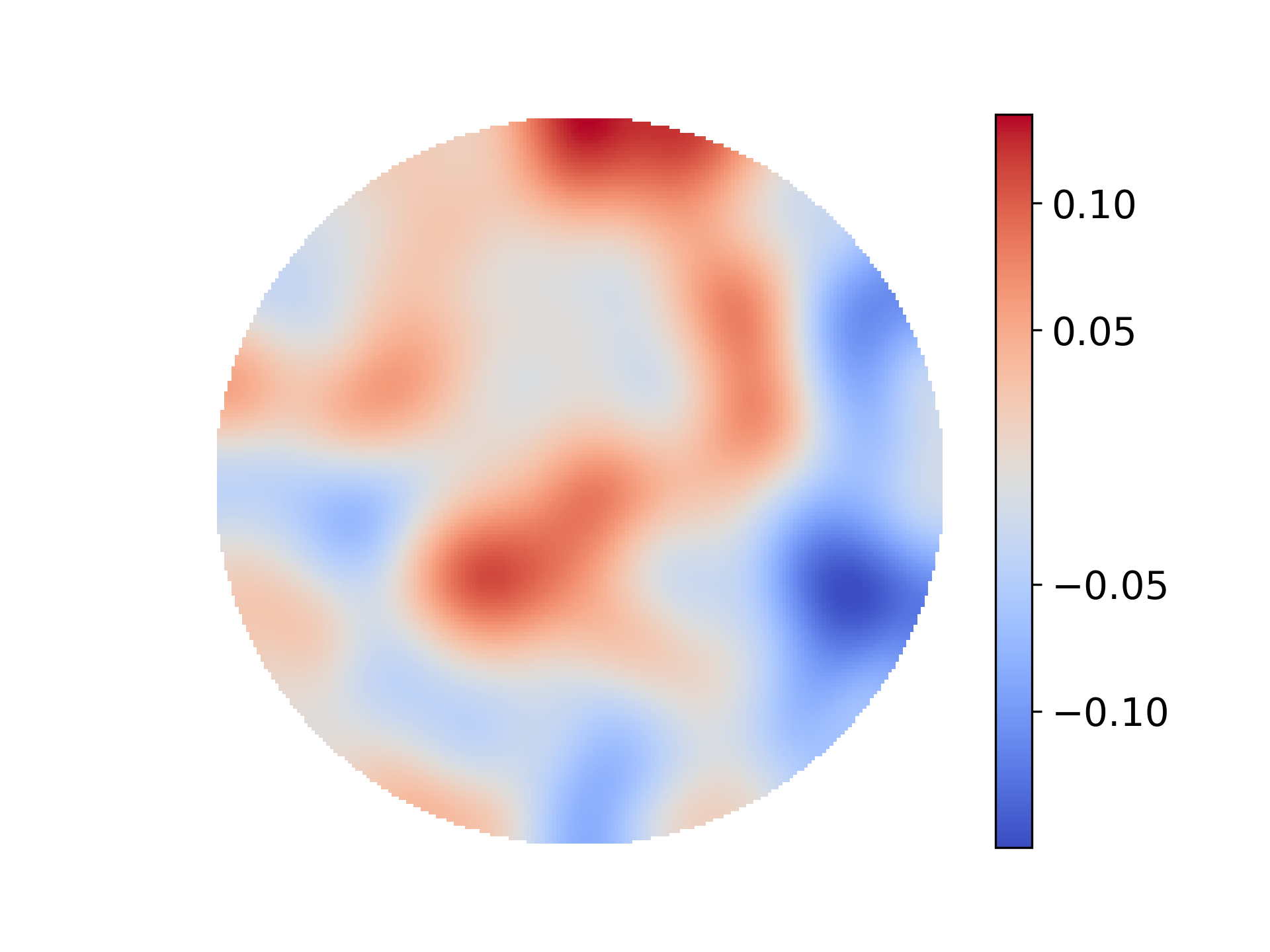}
         \includegraphics[scale=0.21]{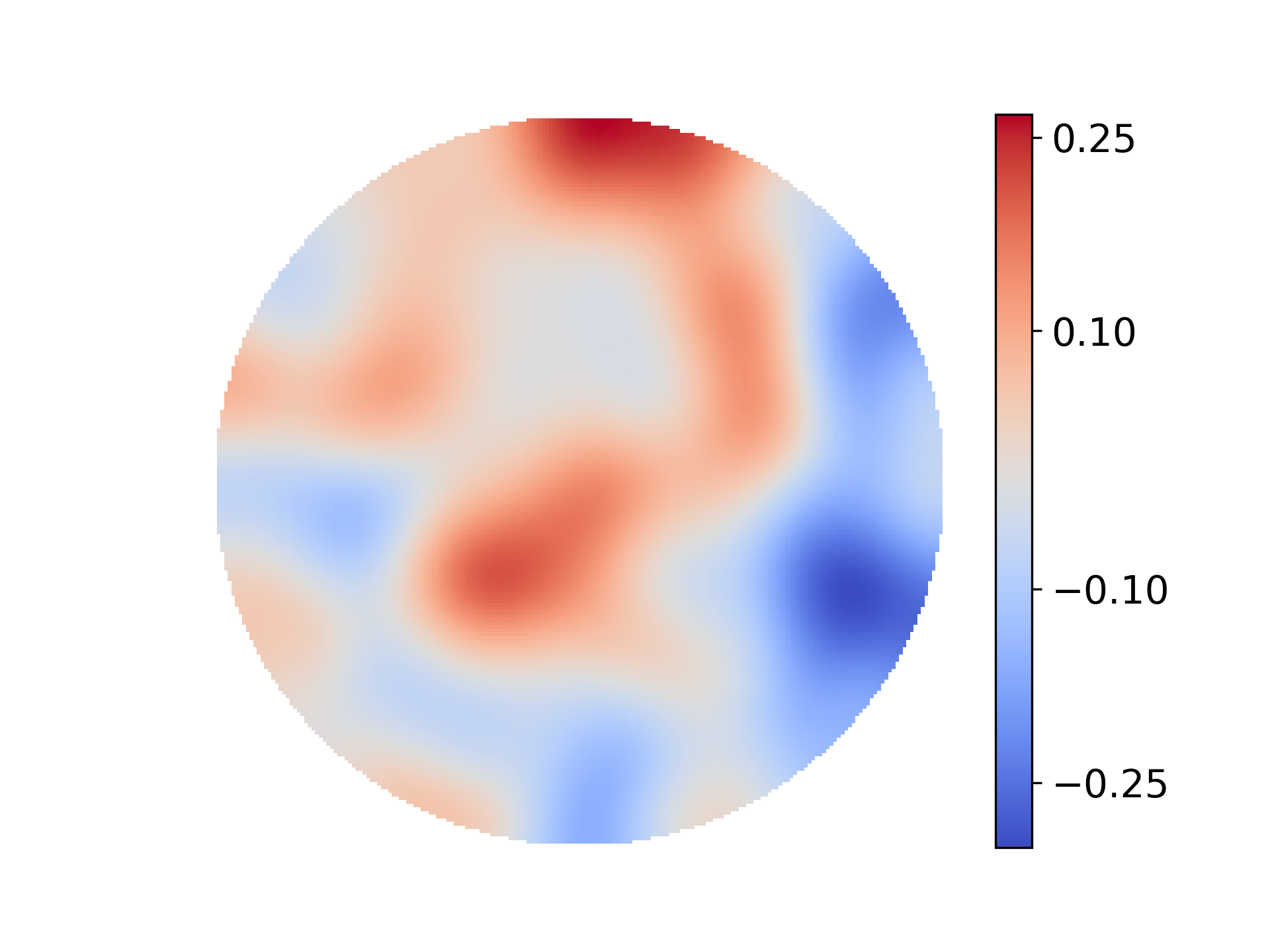}
        \includegraphics[scale=0.21]{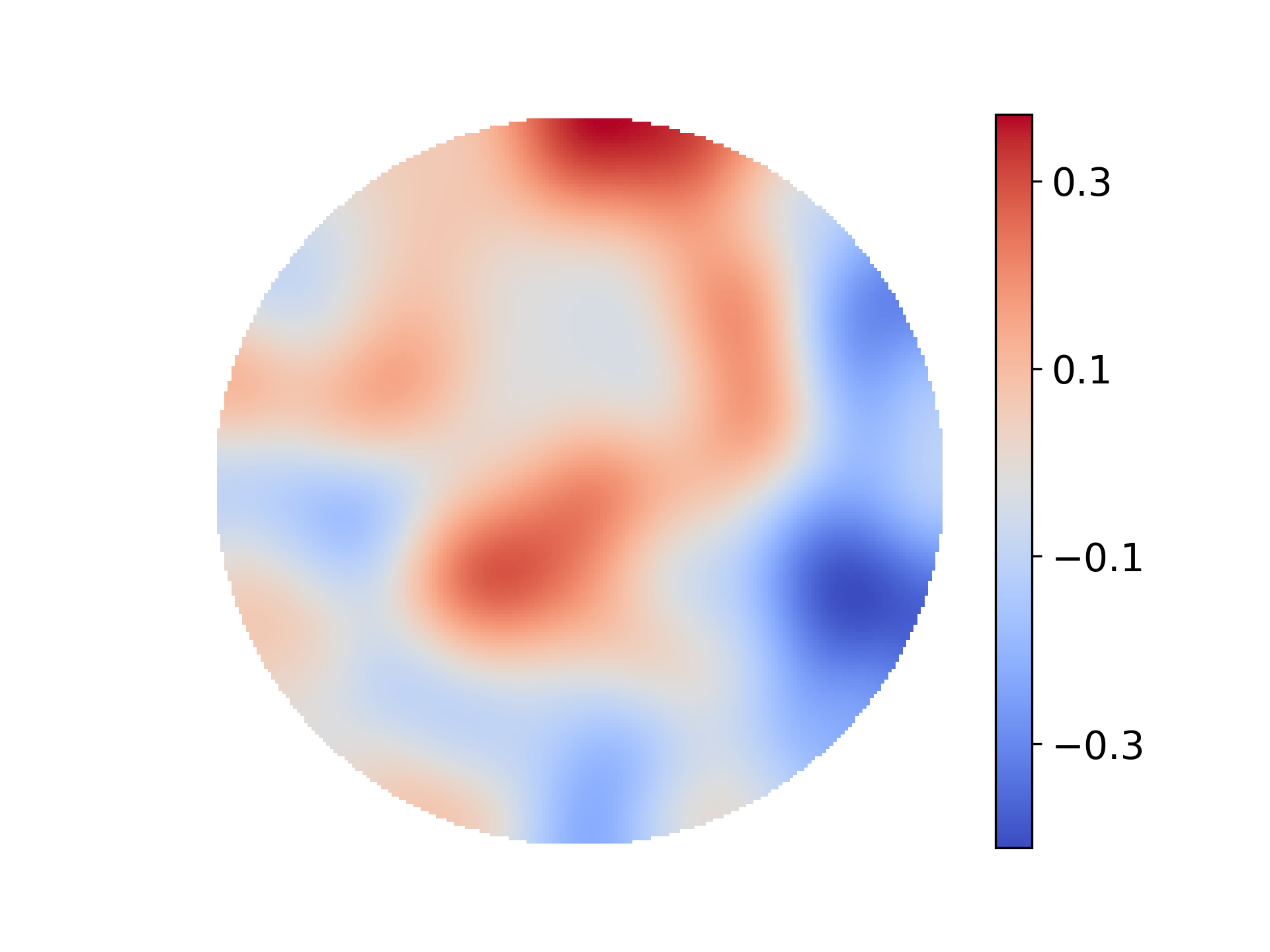}
        \caption*{EDRAS-enhanced PINN solutions for \( M_b=2 \) at \( t=0, 0.5 \) and \( t=1.0 \).}
    \end{subfigure}

 \begin{subfigure}[t]{1\textwidth}
        \centering
        \includegraphics[scale=0.21]{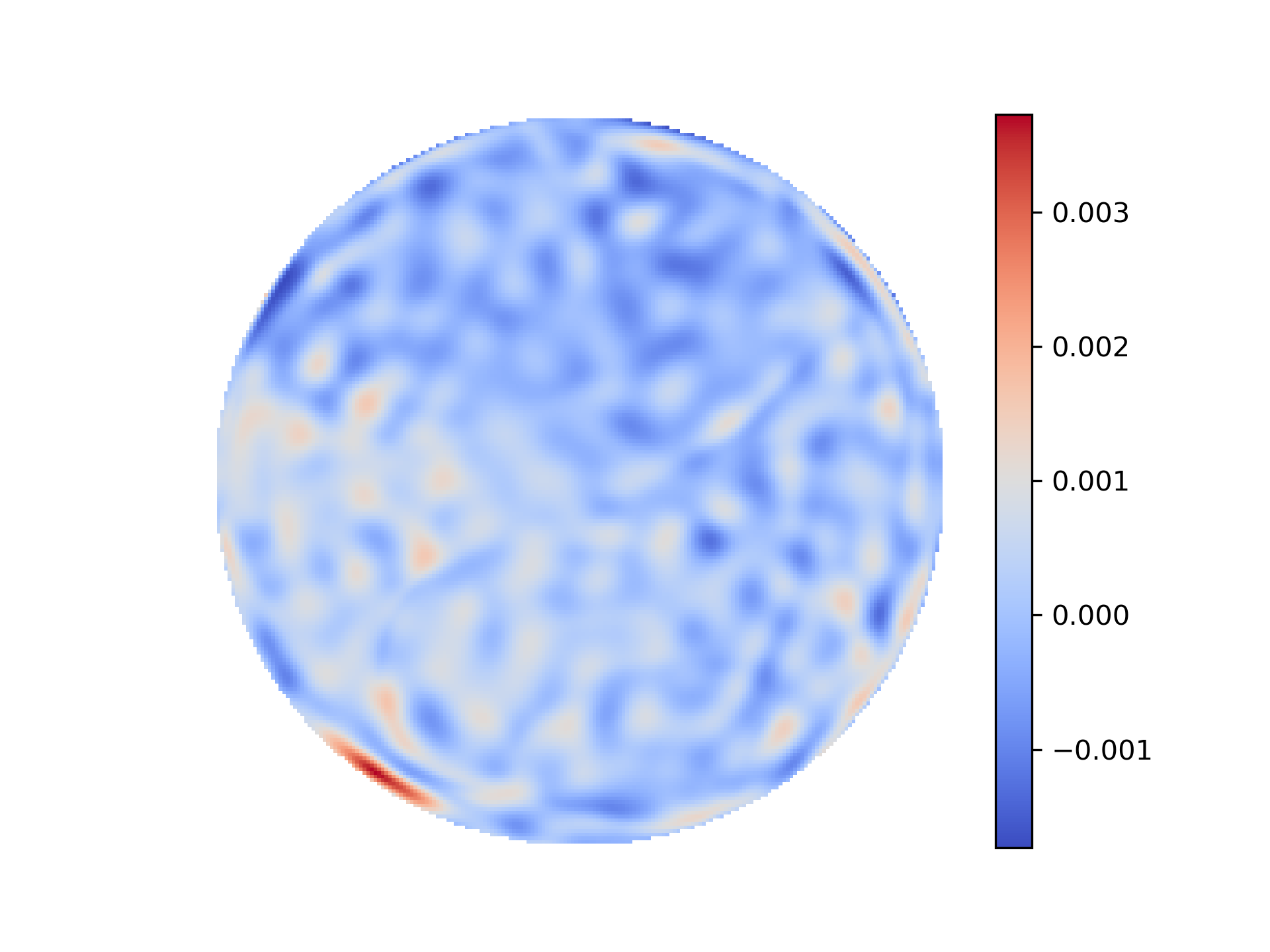}
        \includegraphics[scale=0.21]{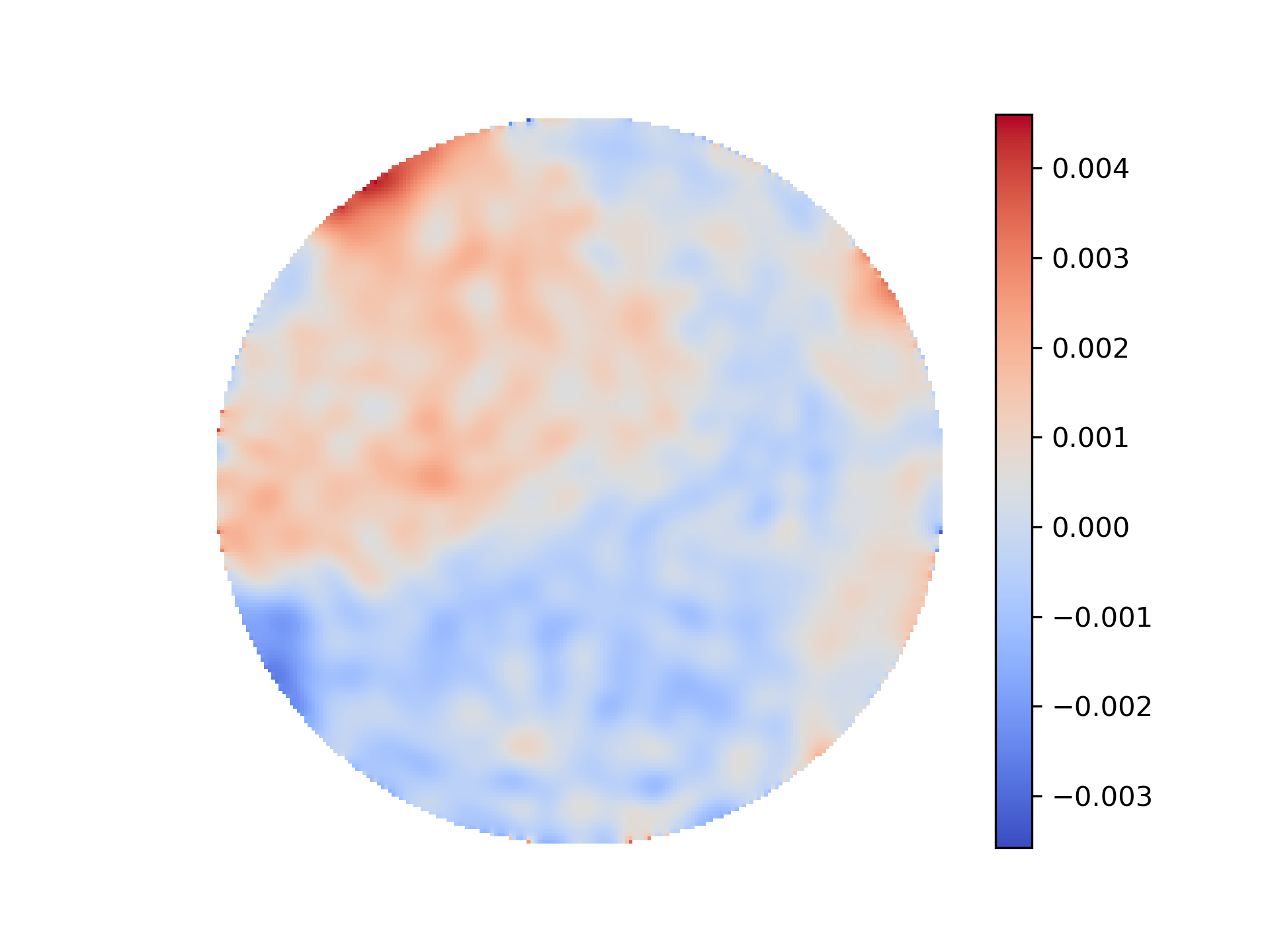}
         \includegraphics[scale=0.21]{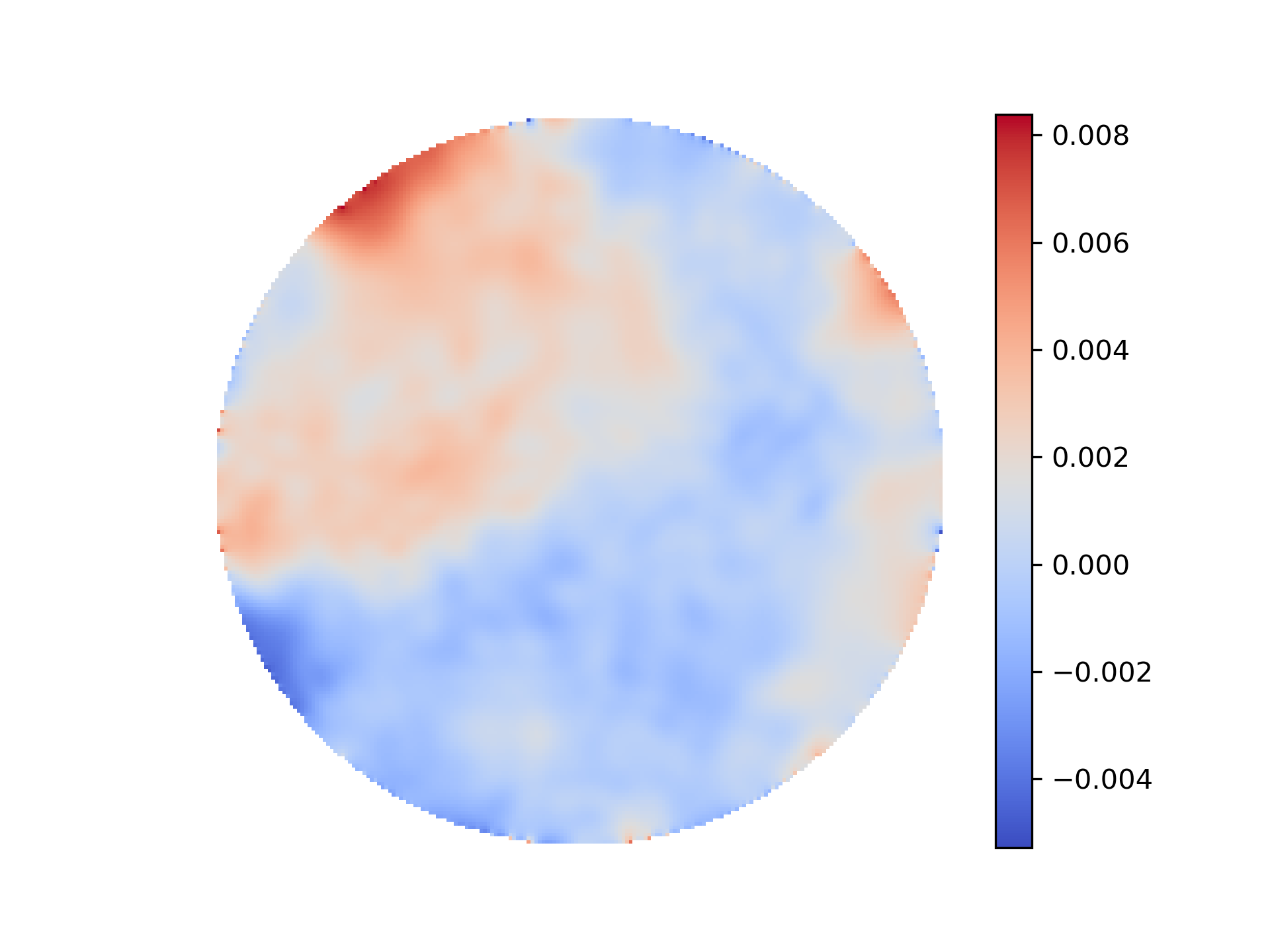}
        \includegraphics[scale=0.21]{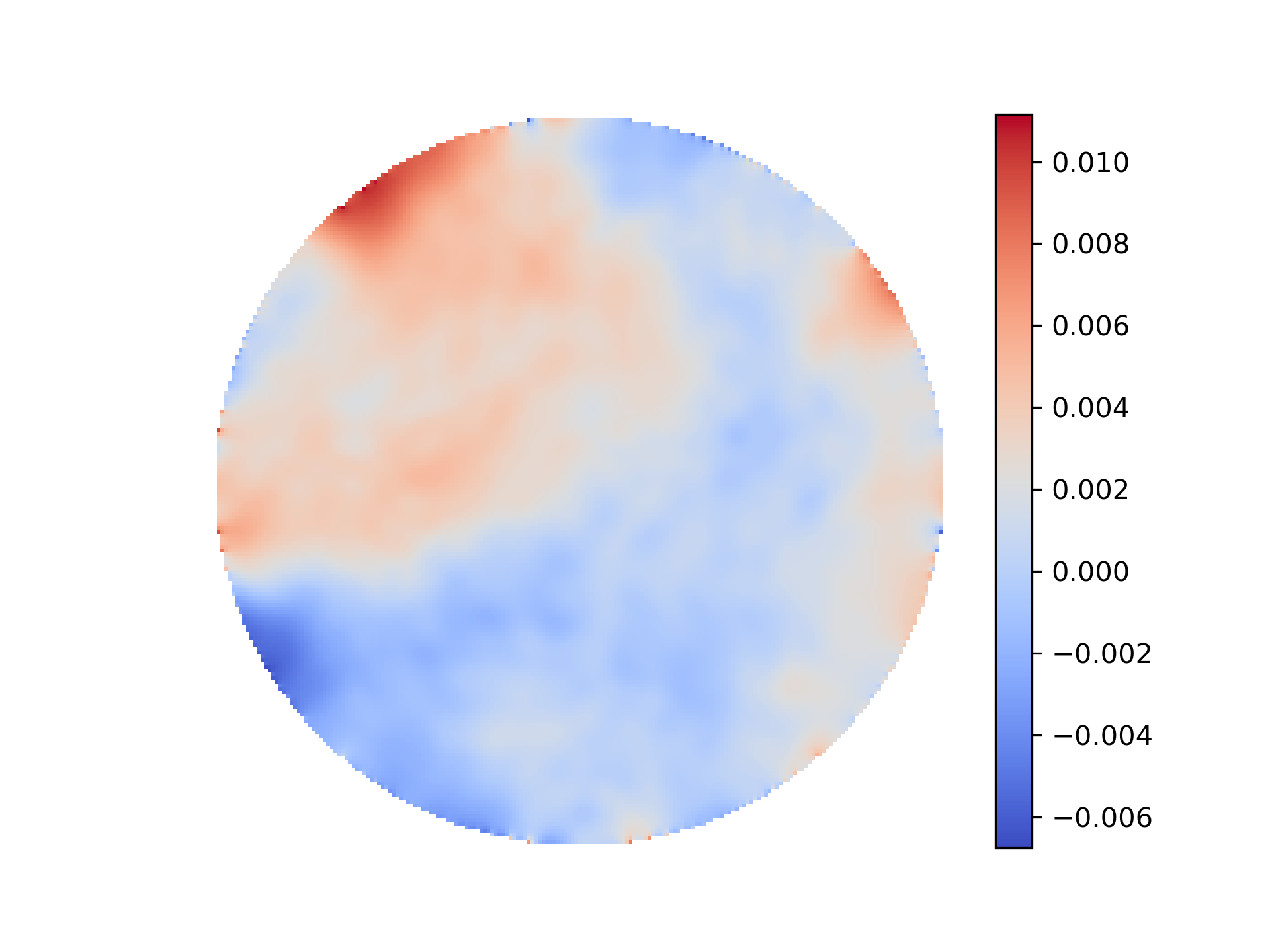}
        \caption*{Error for \( M_b=2 \) at \( t=0, 0.5 \) and \( t=1.0 \).}
    \end{subfigure}

    \begin{subfigure}[t]{1\textwidth}
        \centering
        \includegraphics[scale=0.21]{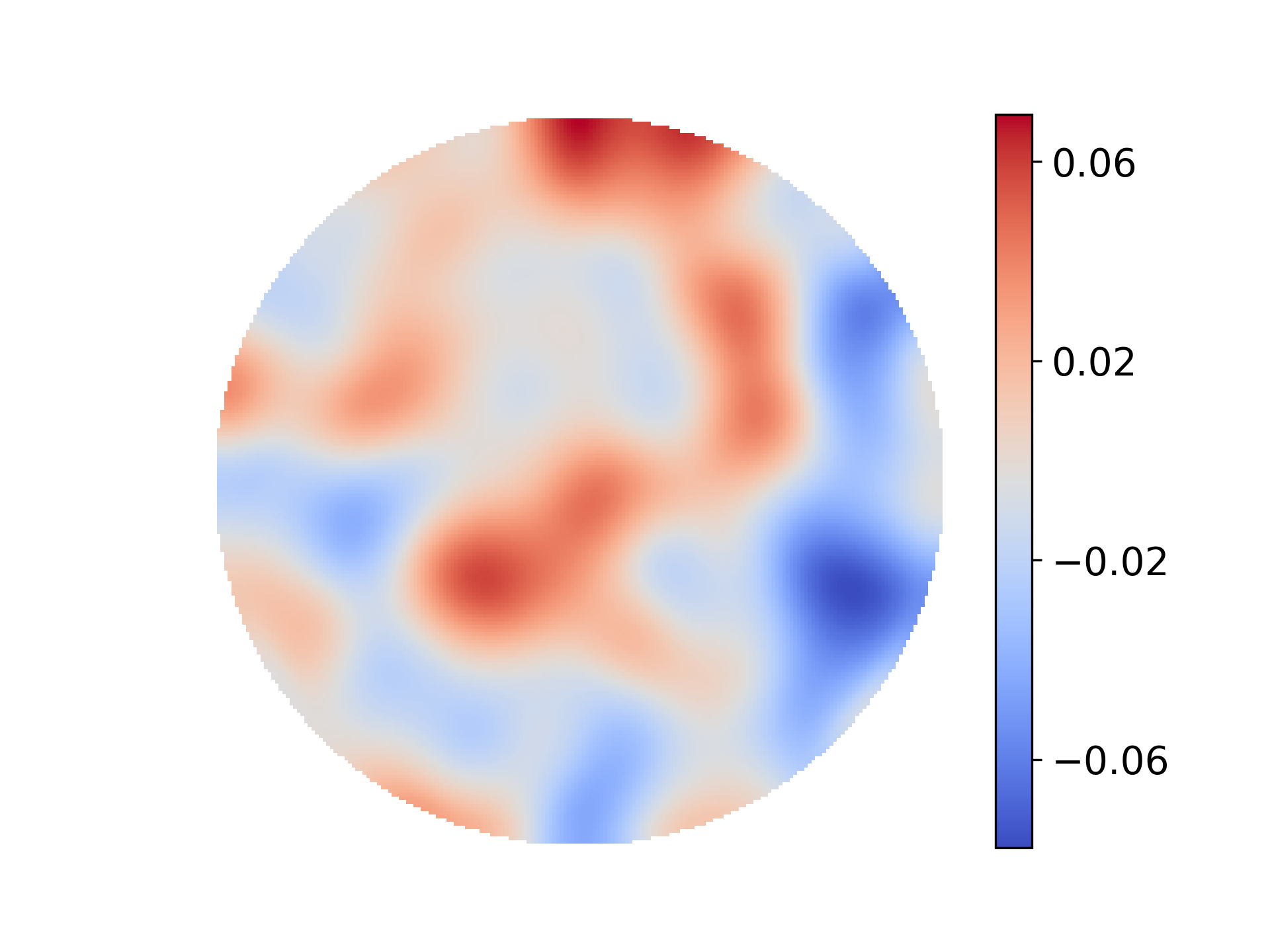}
        \includegraphics[scale=0.21]{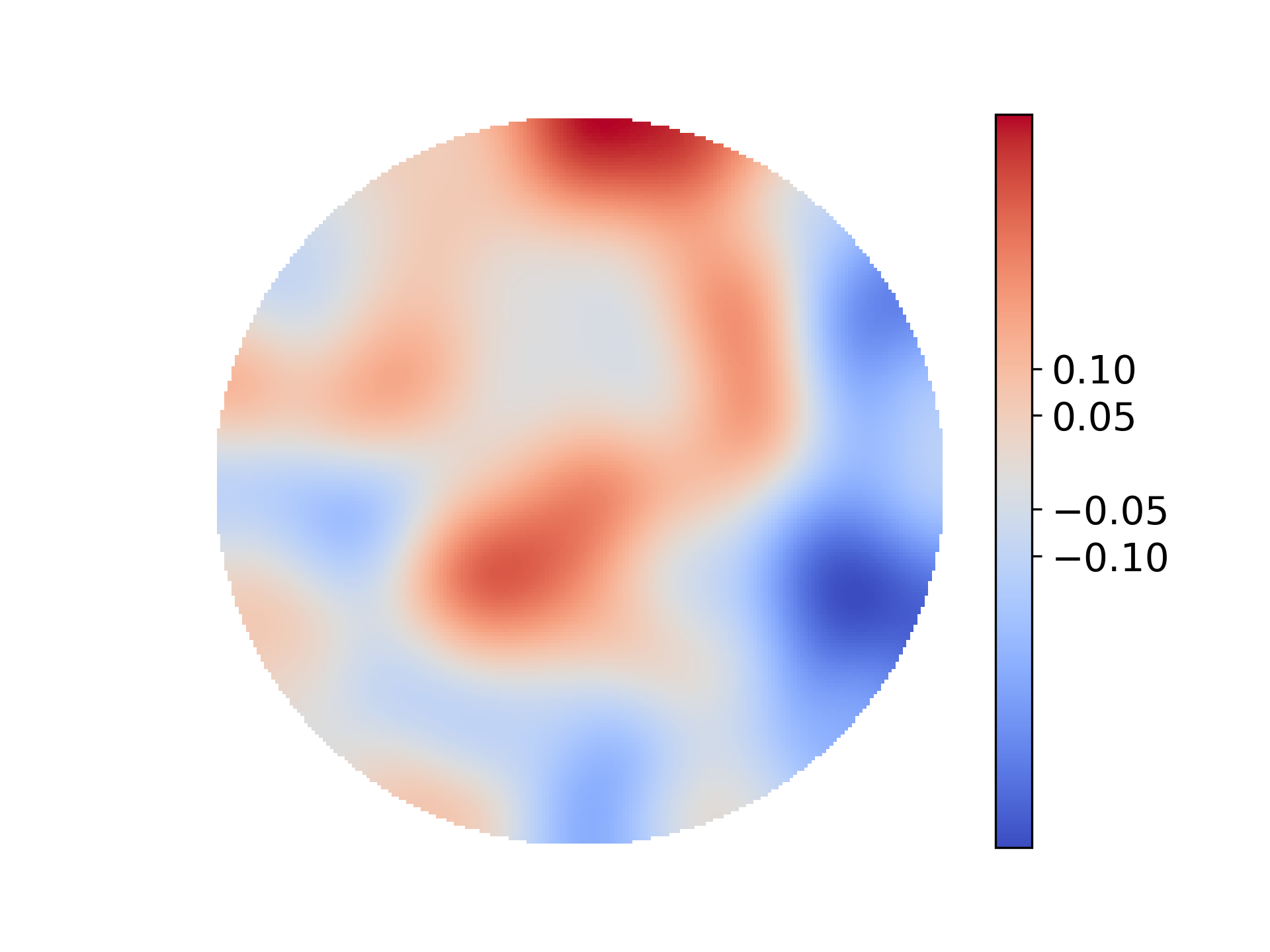}
        \includegraphics[scale=0.21]{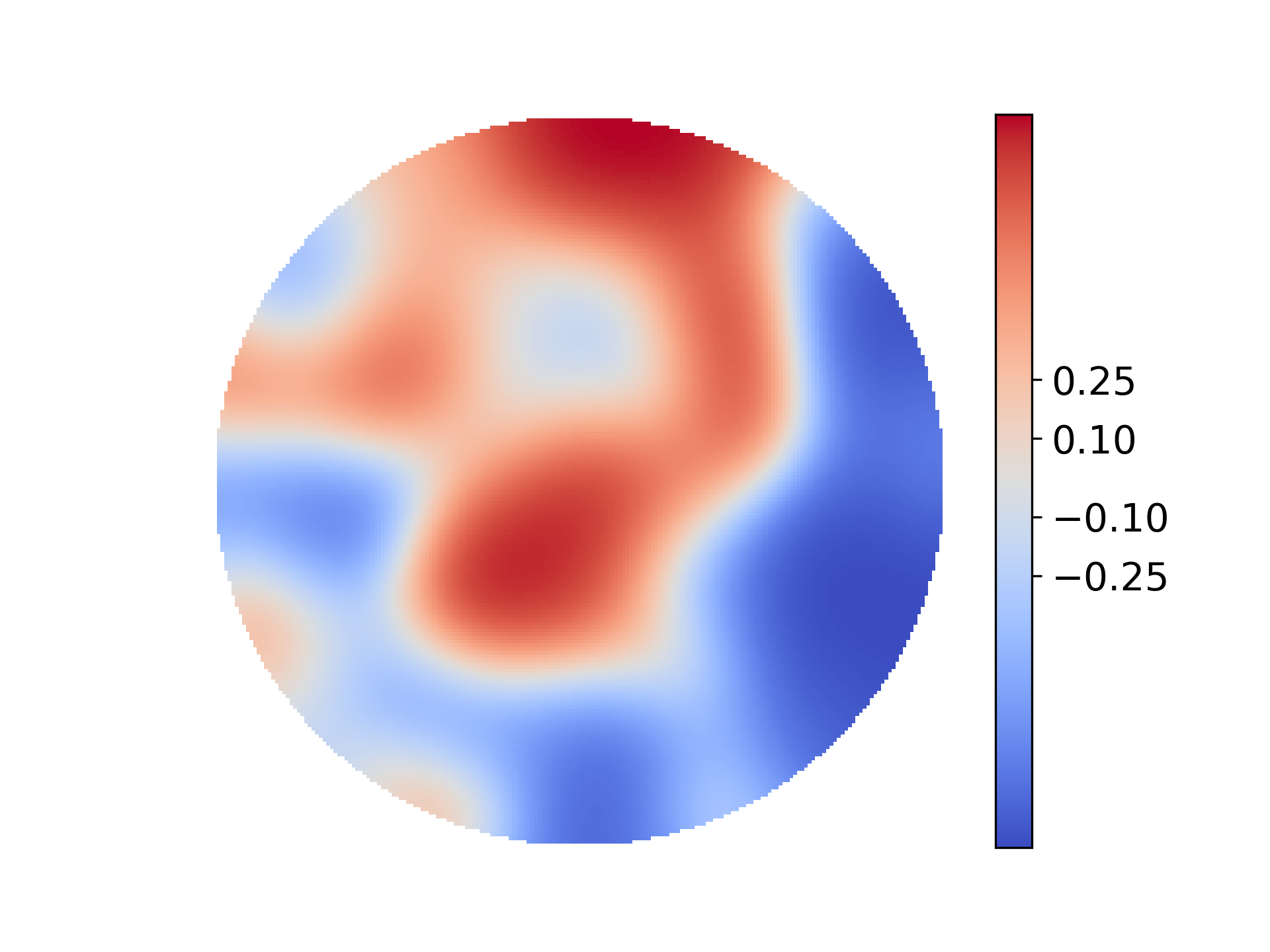}
        \includegraphics[scale=0.21]{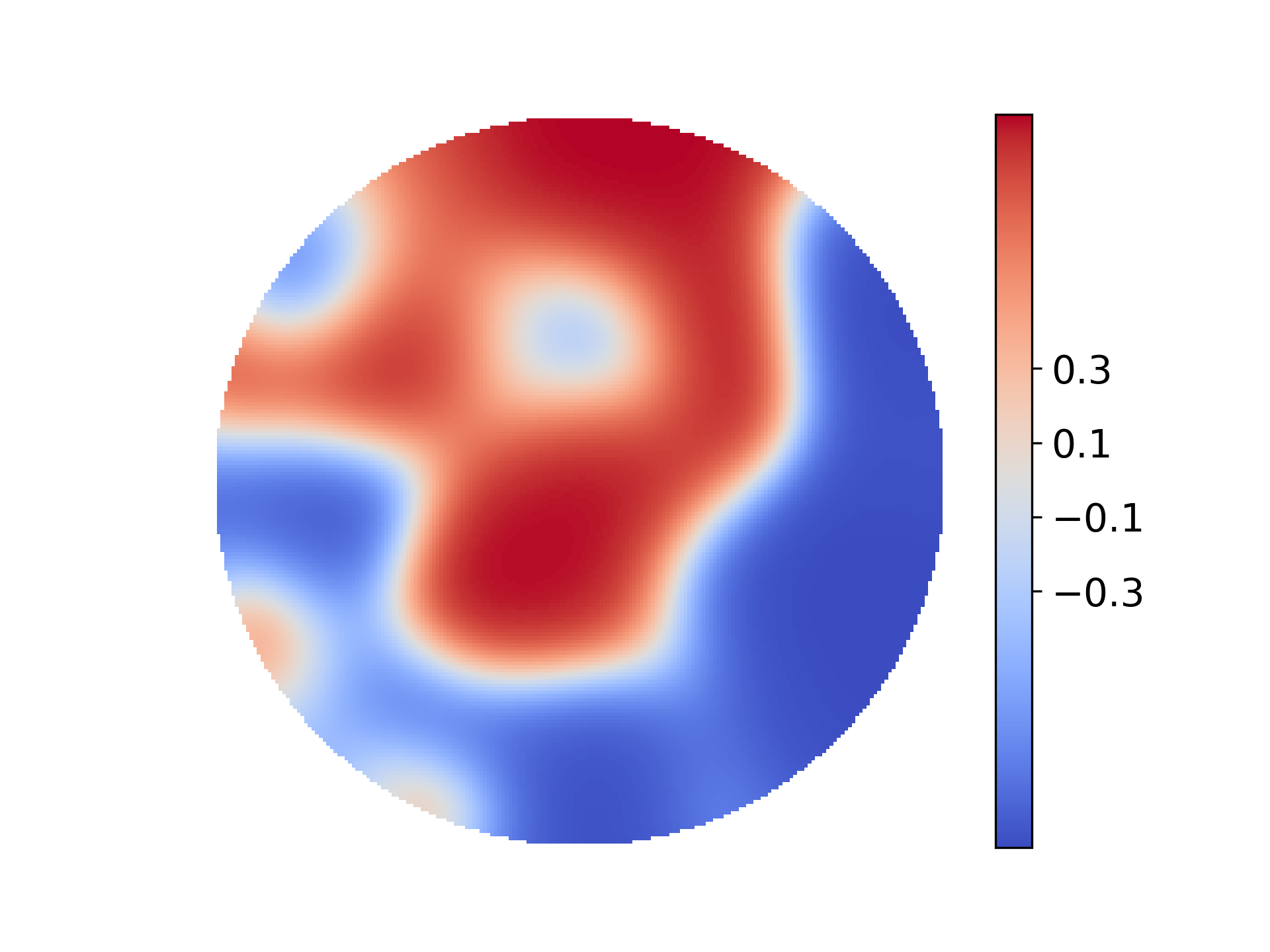}
        \caption*{EDRAS-enhanced PINN solution for \( M_b=5 \) at \( t=0, 0.5 \) and \( t=1.0 \).}
    \end{subfigure}
	
     \begin{subfigure}[t]{1\textwidth}
        \centering
        \includegraphics[scale=0.21]{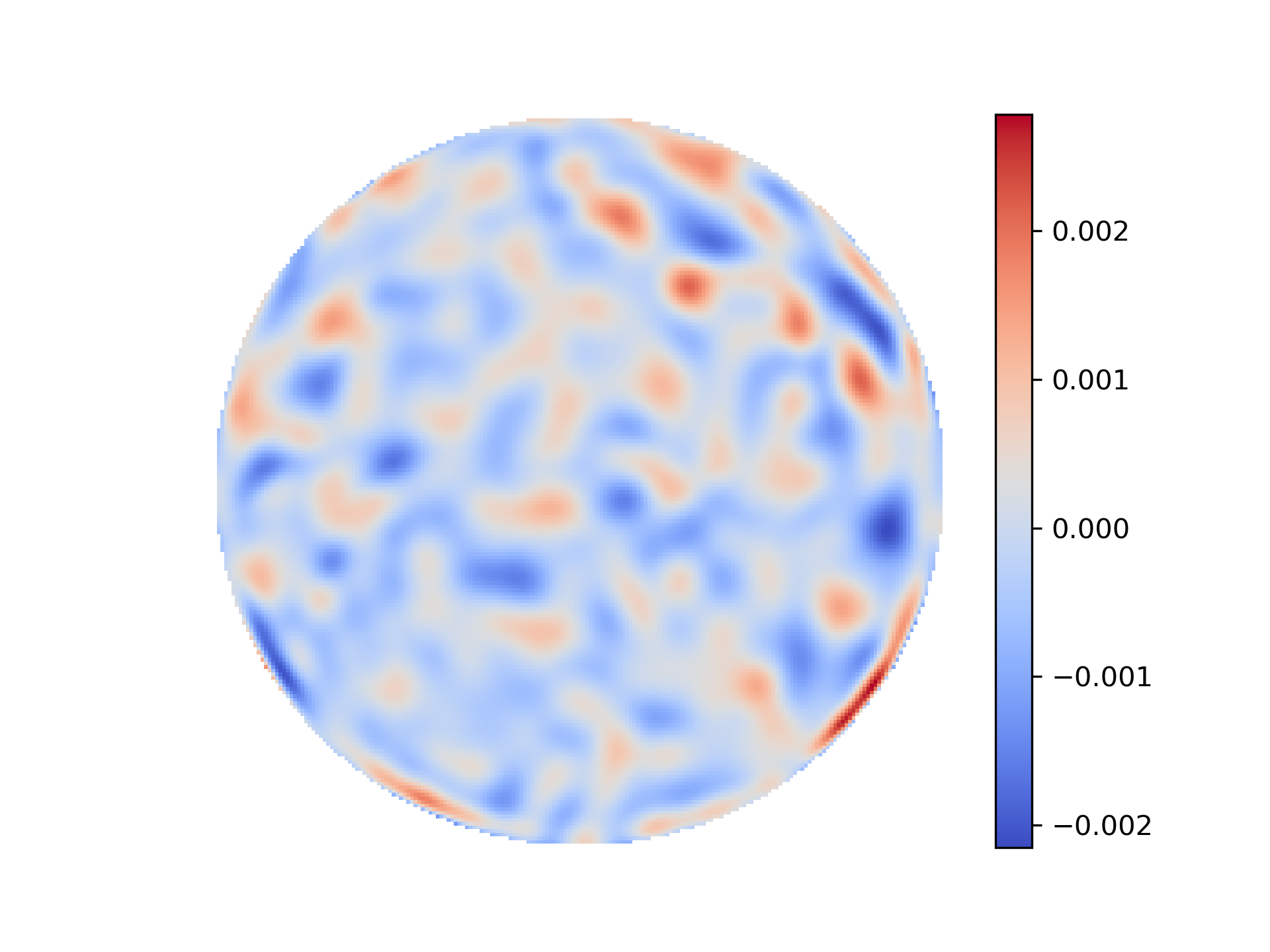}
        \includegraphics[scale=0.21]{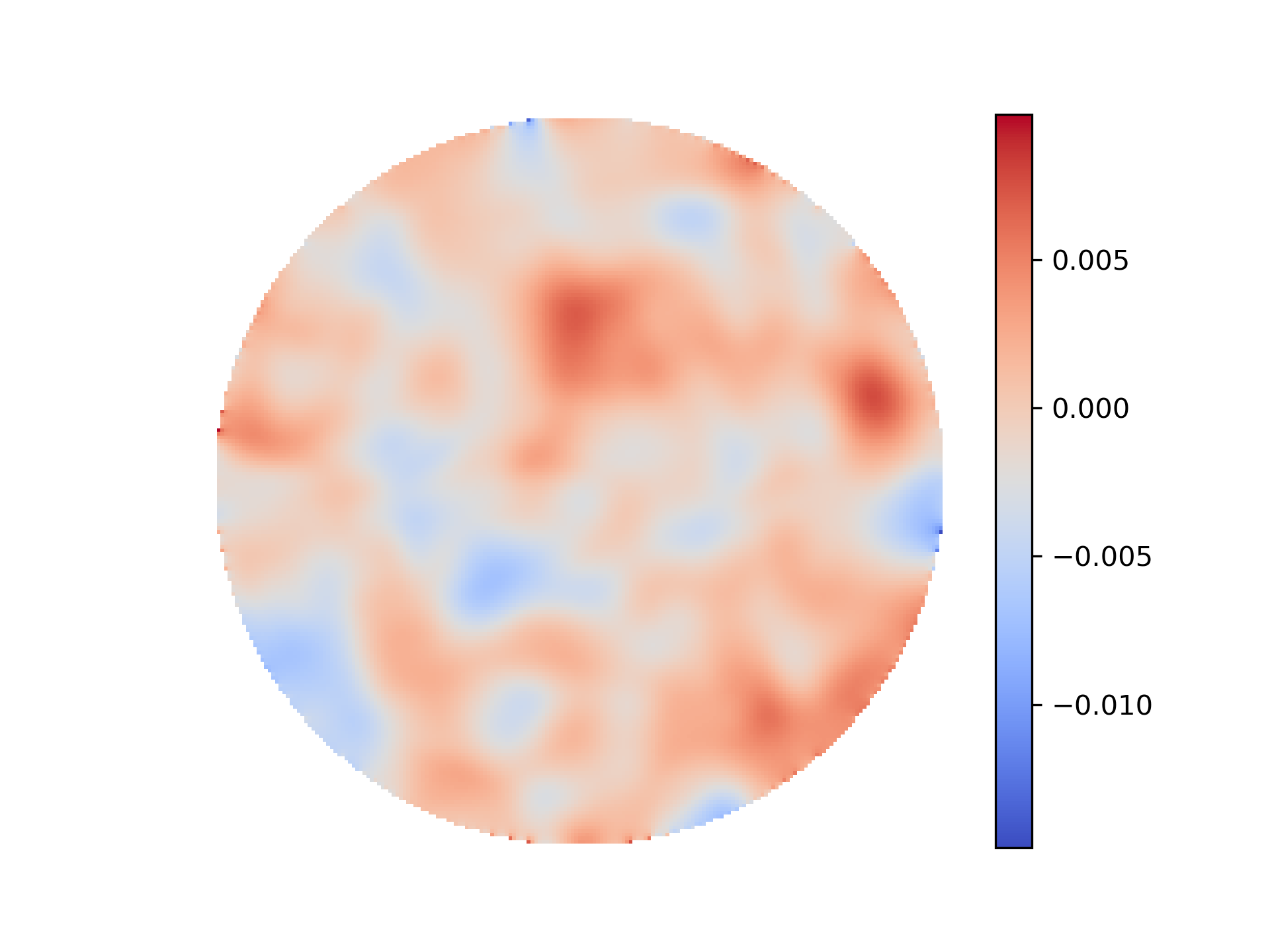}
        \includegraphics[scale=0.21]{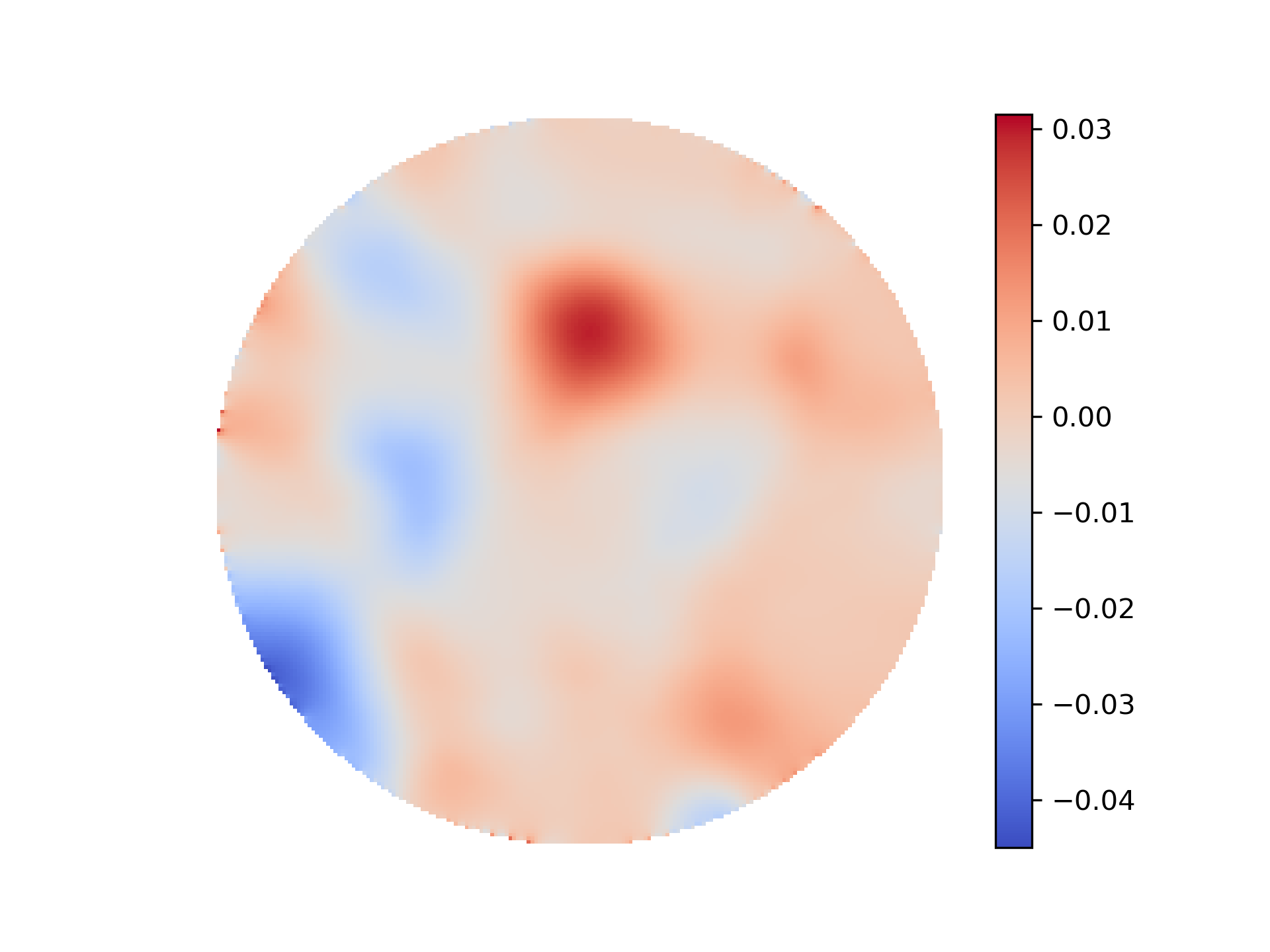}
        \includegraphics[scale=0.21]{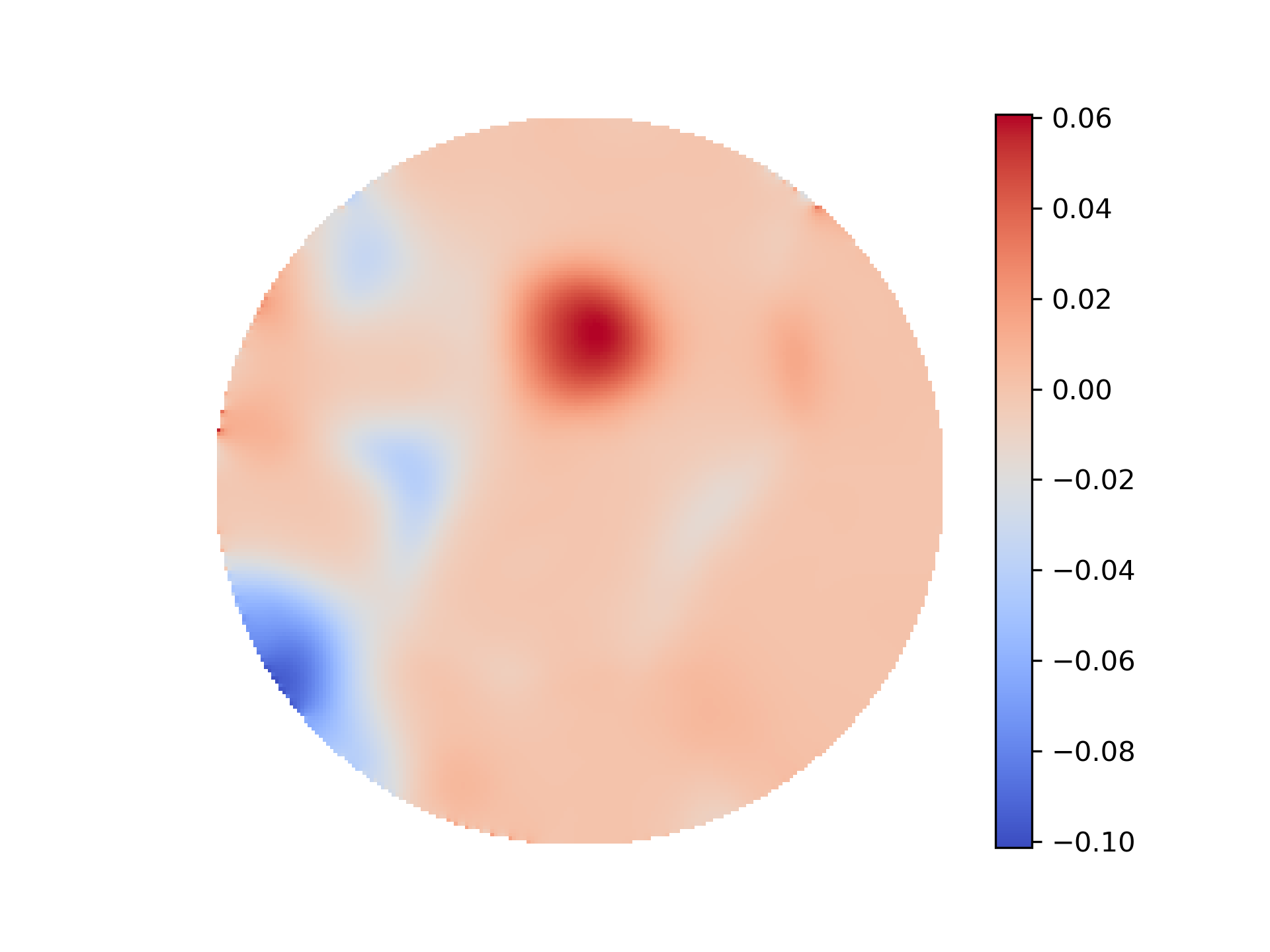}
        \caption*{Error for \( M_b=5 \) at \( t=0, 0.5 \) and \( t=1.0 \).}
    \end{subfigure}

       \begin{subfigure}[t]{1\textwidth}
        \centering
        \includegraphics[scale=0.21]{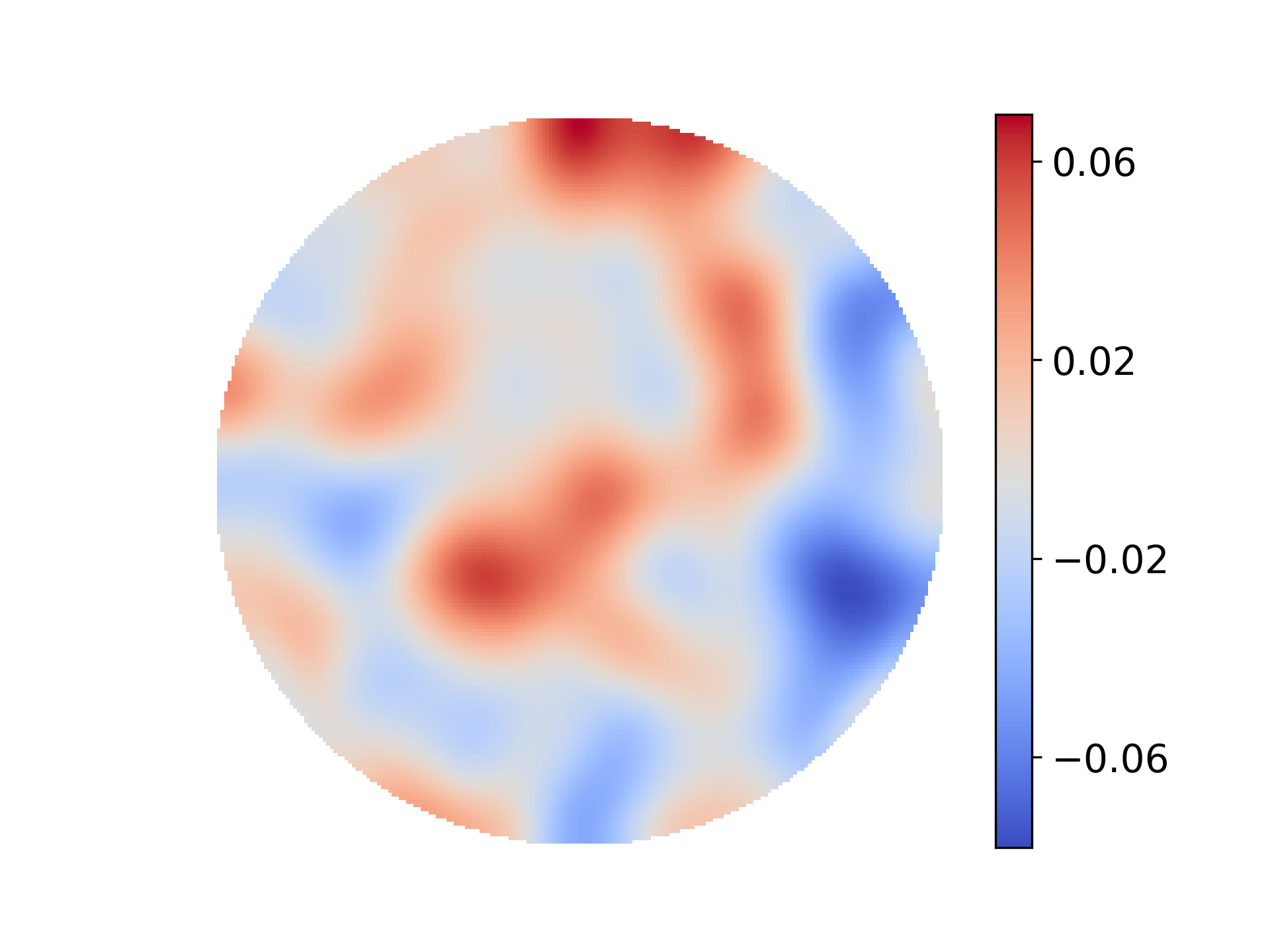}
        \includegraphics[scale=0.21]{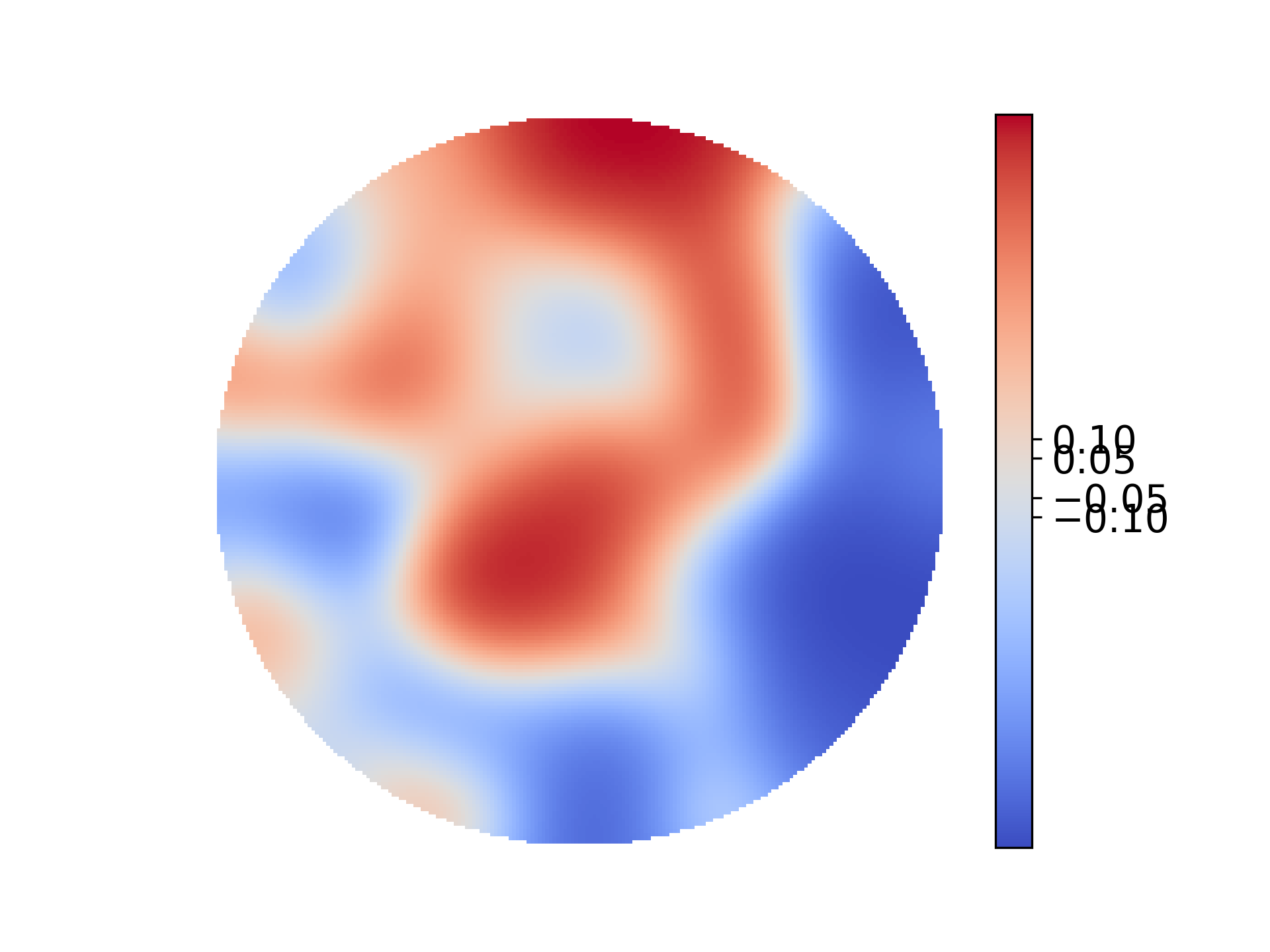}
        \includegraphics[scale=0.21]{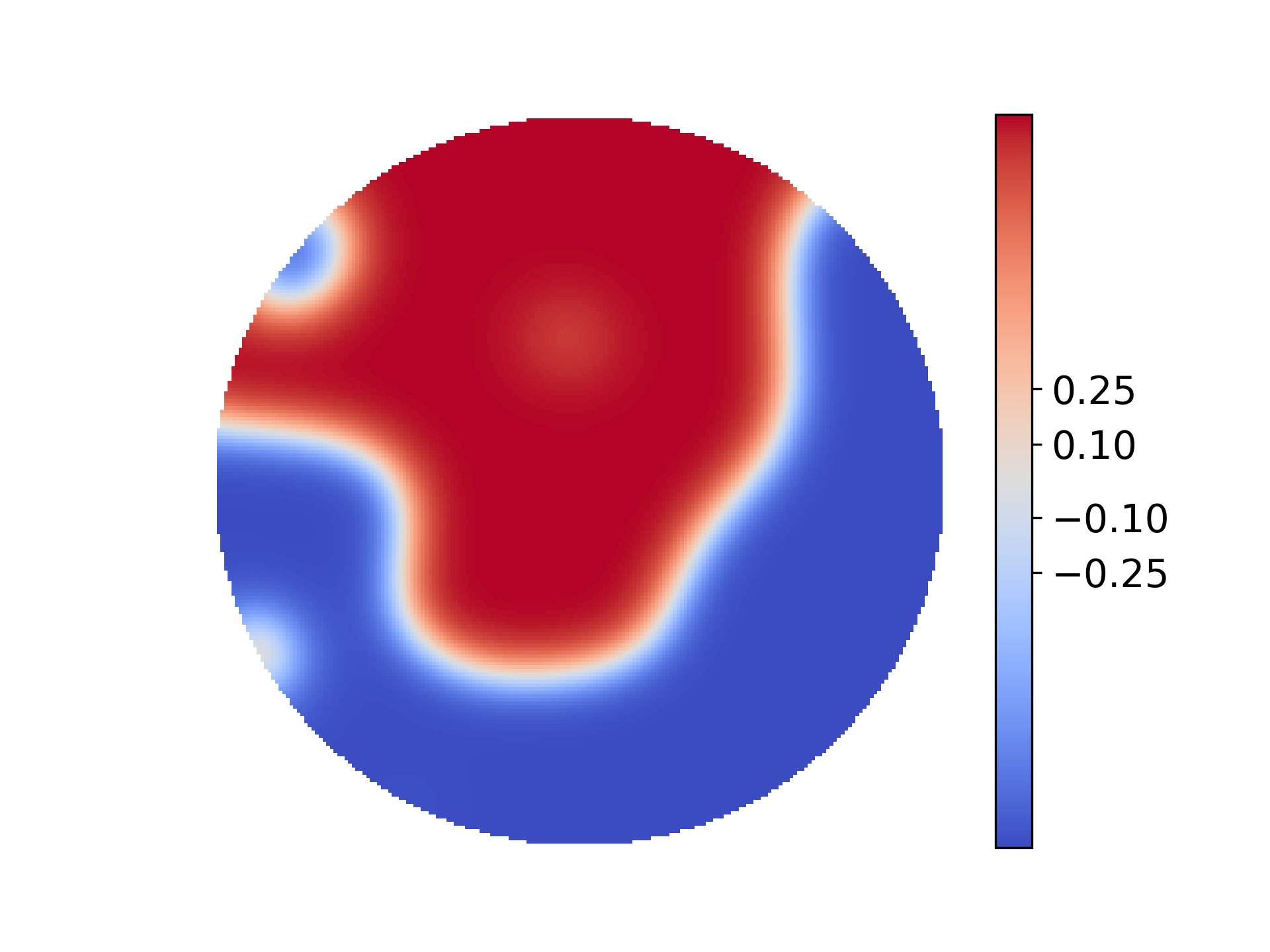}
        \includegraphics[scale=0.21]{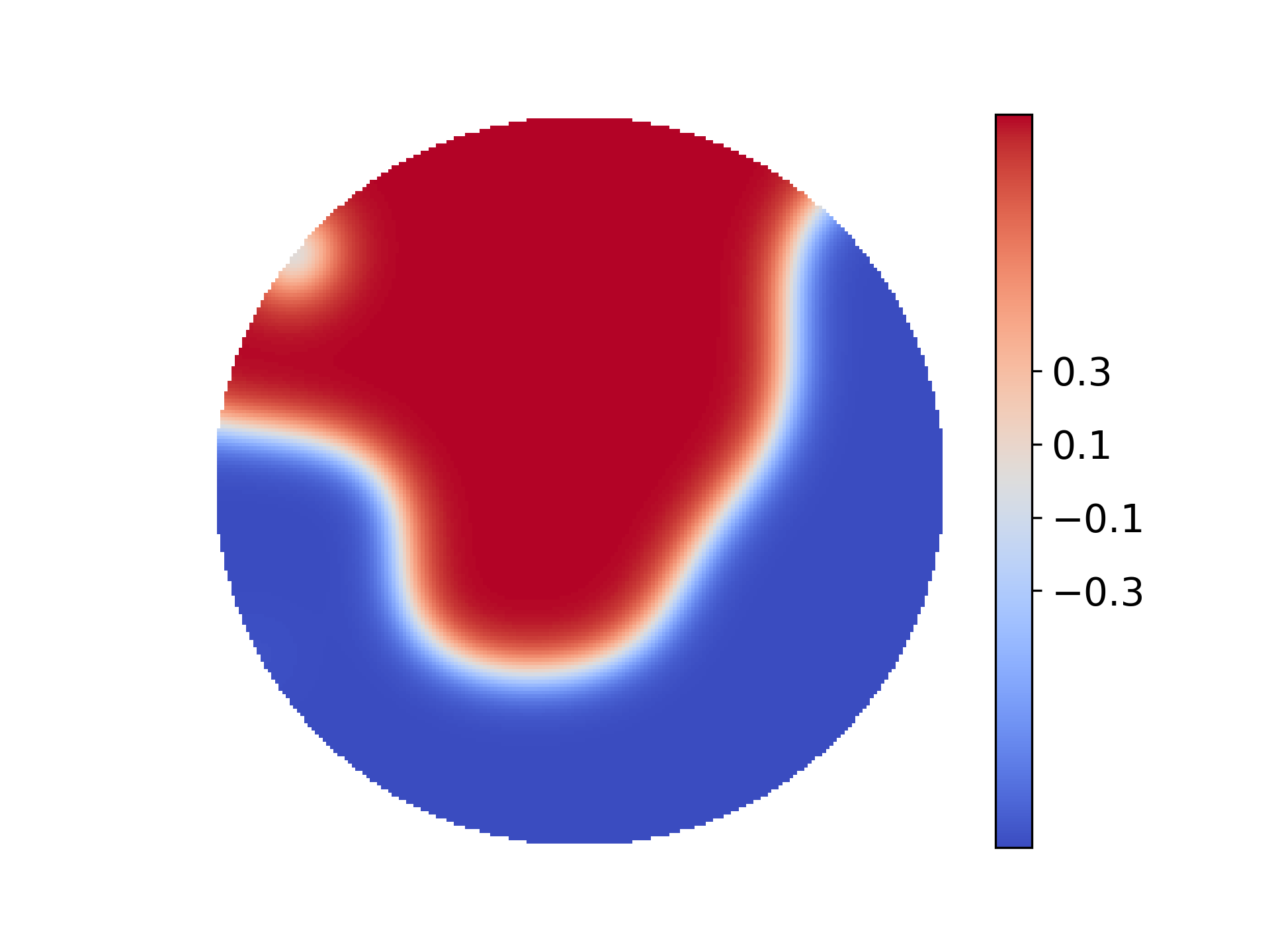}
        \caption*{EDRAS-enhanced PINN solution for \( M_b=10 \) at \( t=0, 0.5 \) and \( t=1.0 \).}
    \end{subfigure}
   \begin{subfigure}[t]{1\textwidth}
        \centering
        \includegraphics[scale=0.21]{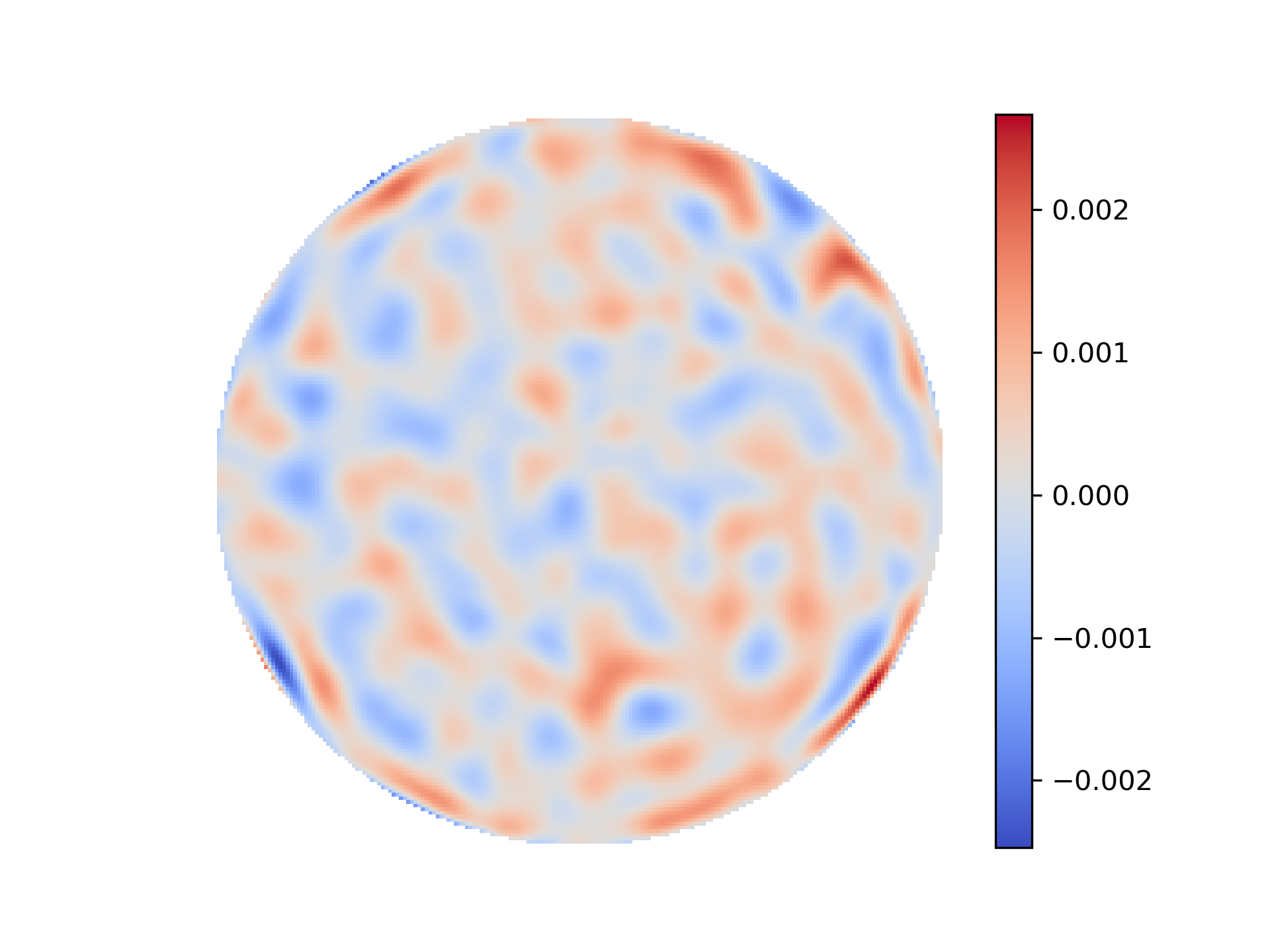}
        \includegraphics[scale=0.21]{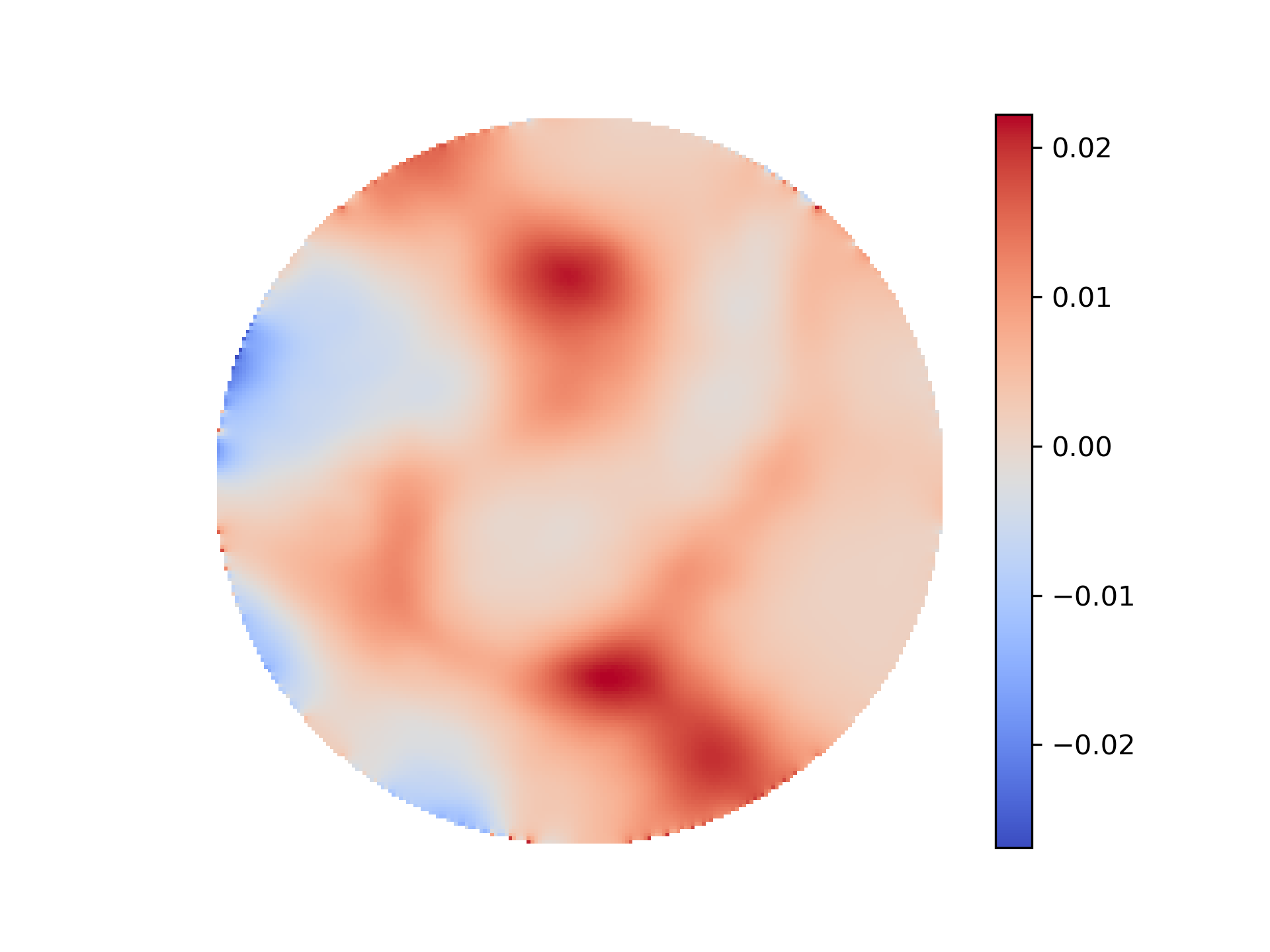}
        \includegraphics[scale=0.21]{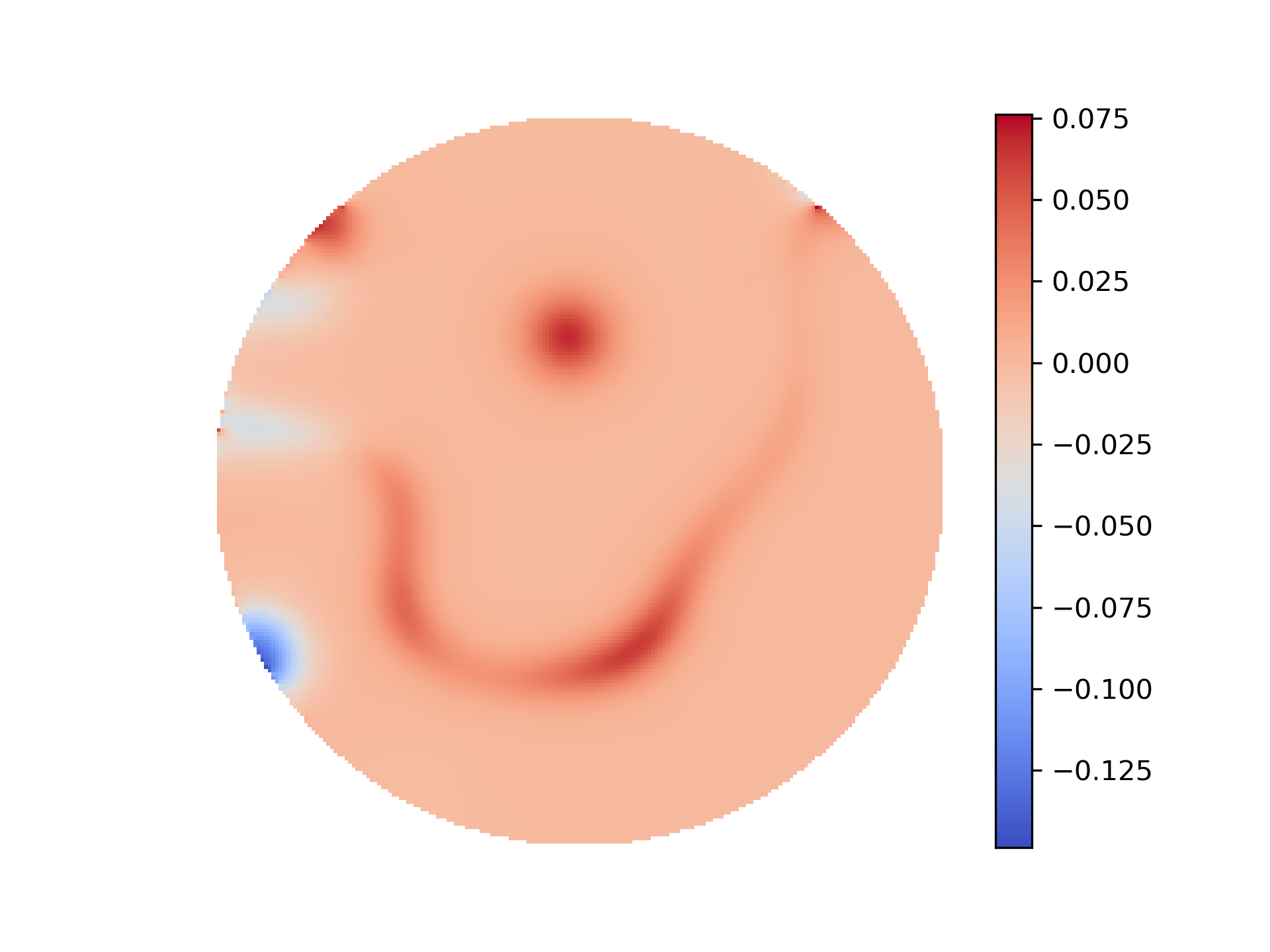}
        \includegraphics[scale=0.21]{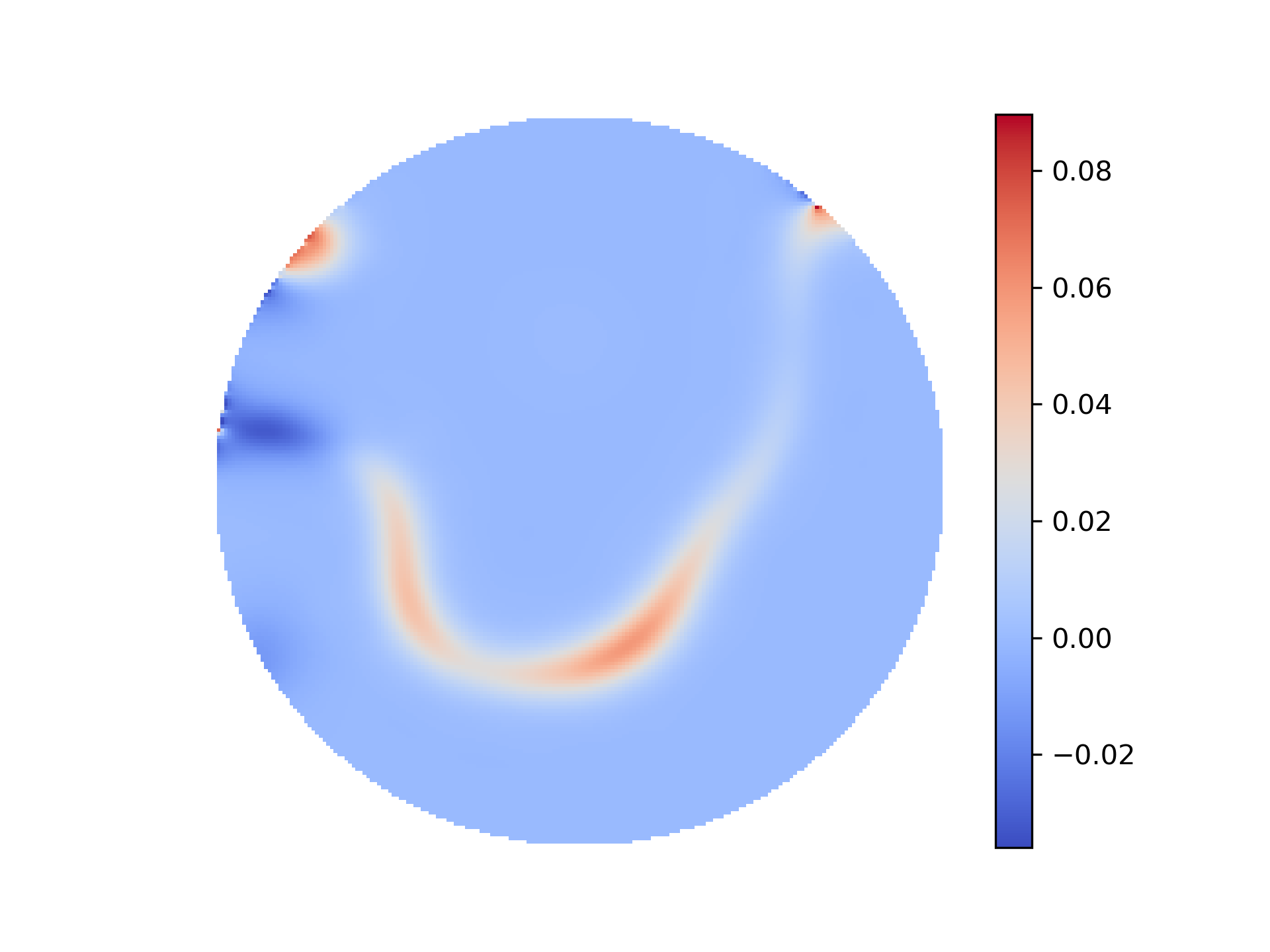}
        \caption*{Error for \( M_b=10 \) at \( t=0, 0.5 \) and \( t=1.0 \).}
    \end{subfigure}

    \caption{Comparison of FDM
     and EDRAS-enhanced PINN solutions of the Allen-Cahn equation under homogeneous Neumann boundary conditions for three mobility values (\( M_b=2, 5, 10 \)) at  (\( t=0, 0.4, 0.8 \) and \( t=1.0 \)), respectively. The results show EDRAS-enhanced PINN solutions and their errors against the FDM solution. The errors are in the order of $O(10^{-2})$.}
    \label{fig:neu-2-5-10}
\end{figure}

To validate the EDRAS-enhanced PINN approach, we benchmark it against a reference solution obtained using finite difference methods \cite{yu2024domain} for the Allen-Cahn model with Neumann boundary conditions. Figure \ref{fig:neu-2-5-10} illustrates a comparison between the reference solution computed by the finite difference method and the  solution generated by the EDRAS-enhanced PINN method for various bulk mobility values under homogeneous Neumann boundary conditions. Comprehensive error metrics across different bulk mobility values are presented in Table \ref{tab:err_2d}. These results clearly demonstrate the effectiveness, accuracy and robust performance of the EDRAS-enhanced PINN approach.
\begin{table}[htbp]
\centering
\caption{Error of the Allen-Cahn equation with the Neumann boundary condition.} 
\label{tab:err_2d}
\begin{tabular}{cccc}
\hline
    $M_b$  & MAE & MSE & Relative MSE  \\
\hline
2   & 0.00168 & $6.75 \times 10^{-5}$ & 0.011200\\
5  & 0.00399 & $1.10 \times 10^{-4}$ & 0.000822 \\
10  & 0.00549 & $1.85 \times 10^{-4}$ & 0.000410\\
\hline
\end{tabular}
\end{table}

We further investigate the impact of hyperparameter $d_f^{(0)}$, the threshold of the  collocation point density, on the performance of the EDRAS enhanced PINNs by carrying out a series of experiments for a 2D Allen-Cahn equation with the Neumann boundary condition. Table \ref{tab:den} shows a non-monotonic relationship between  density threshold ($d^{(0)}_f$) and prediction errors (MAE/MSE/RMSE), highlighting two key trends:
\begin{itemize}
\item First, moderate density improves accuracy. For a fixed $M_b$ (e.g., $M_b = 2$), increasing $d^{(0)}_f$ from 150 to 200 reduces both MAE ($3.32 \times 10^{-3} \rightarrow 2.32 \times 10^{-3}$) and MSE ($7.96 \times 10^{-5} \rightarrow 7.00 \times 10^{-5}$), suggesting that higher sampling density better constrains the PDE solution. Similar improvements occur for $M_b = 5$ (e.g., $d^{(0)}_f = 150 \rightarrow 200$ lowers Relative MSE by $\sim$50\%) and $M_b = 10$ (e.g., $d^{(0)}_f = 100 \rightarrow 150$ lowers Relative MSE as well as MSE by $\sim$40\%).
    \item Second, saturation and rebound occurs at the high density. Beyond a threshold, further increases in $d^{(0)}_f$ degrade the performance adversely. For $M_b = 2$, raising $d^{(0)}_f$ to 300 MAE/MSE  \textit{increases} compared to the case of $d^{(0)}_f = 200$. Hence, one needs to balance the sampling efficiency and computational feasibility when choosing the threshold, $d_f^{(0)}$.
        \end{itemize}

\begin{table}[h]
\centering
\caption{Performance metrics for different configurations}
\label{tab:den}
\begin{tabular}{cccccc}
\hline
$M_b$ & $d^{(0)}_f$ & MAE & MSE & Relative MSE \\
\hline
2 & 150 & $3.32 \times 10^{-3}$ & $7.96 \times 10^{-5}$ & $1.32 \times 10^{-2}$ \\
  & 200 & $2.32 \times 10^{-3}$ & $7.00 \times 10^{-5}$ & $1.16 \times 10^{-2}$ \\
  & 300 & $2.86 \times 10^{-3}$ & $7.36 \times 10^{-5}$ & $1.22 \times 10^{-2}$ \\
\hline
5 & 150 & $6.64 \times 10^{-3}$ & $1.29 \times 10^{-4}$ & $9.62 \times 10^{-4}$ \\
  & 200 & $5.02 \times 10^{-3}$ & $7.93 \times 10^{-5}$ & $5.92 \times 10^{-4}$ \\
  & 300 & $6.36 \times 10^{-3}$ & $1.25 \times 10^{-4}$ & $9.35 \times 10^{-4}$ \\
\hline
10 & 100 & $9.65 \times 10^{-3}$ & $1.08 \times 10^{-3}$ & $2.40 \times 10^{-3}$ \\
   & 150 & $6.74 \times 10^{-3}$ & $2.56 \times 10^{-4}$ & $5.67 \times 10^{-4}$ \\
   & 200 & $8.56 \times 10^{-3}$ & $4.95 \times 10^{-4}$ & $1.10 \times 10^{-3}$ \\
\hline
\end{tabular}
\end{table}

\begin{figure}[H]
    \centering
    \begin{subfigure}[t]{1\textwidth}
        \centering
        \includegraphics[scale=0.22]{fig/ac2d/neumann-Mb2/u_pred_0.0000.png}
        \includegraphics[scale=0.22]{fig/ac2d/neumann-Mb2/u_pred_0.4000.png}
        \includegraphics[scale=0.22]{fig/ac2d/neumann-Mb2/u_pred_0.8000.png}
        \includegraphics[scale=0.22]{fig/ac2d/neumann-Mb2/u_pred_1.0000.png}
        \caption{Snapshots at $t=0, 0.4, 0.8, 1.0$ for Neumann boundary conditions with \( M_b=2\).}
        \label{fig:static-dy-1}
    \end{subfigure}

    \begin{subfigure}[t]{1\textwidth}
        \centering
        \includegraphics[scale=0.22]{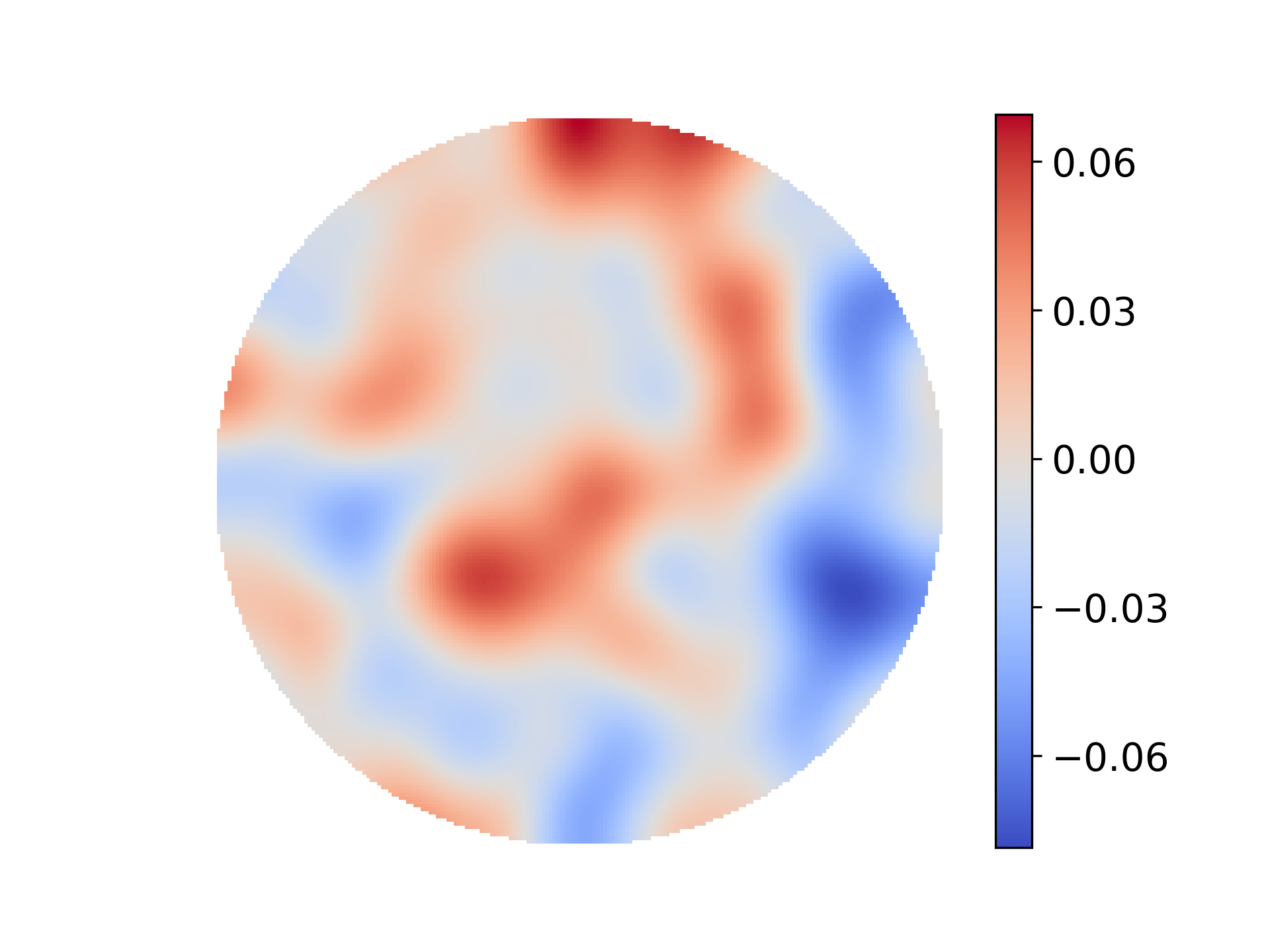}
        \includegraphics[scale=0.22]{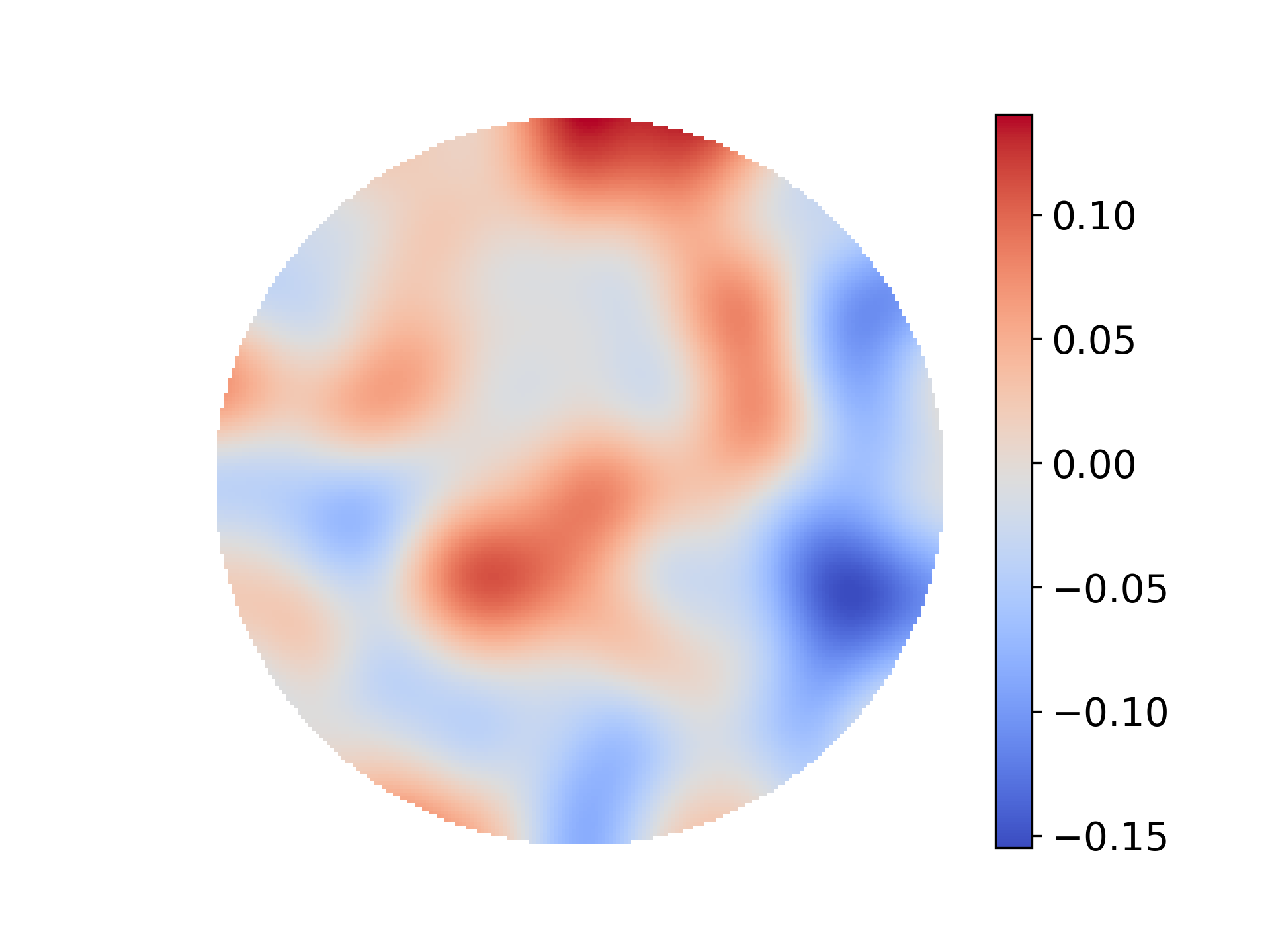}
        \includegraphics[scale=0.22]{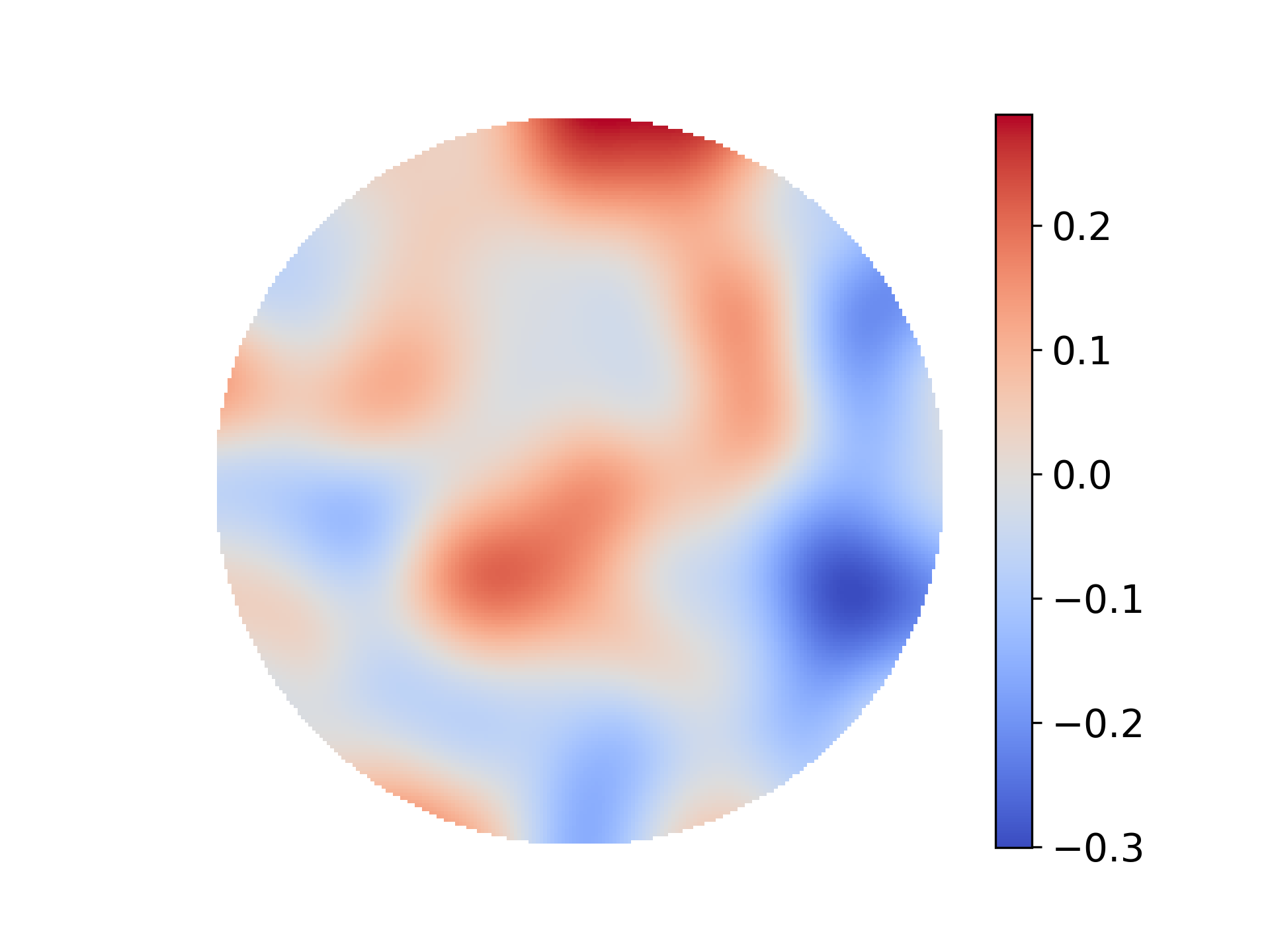}
        \includegraphics[scale=0.22]{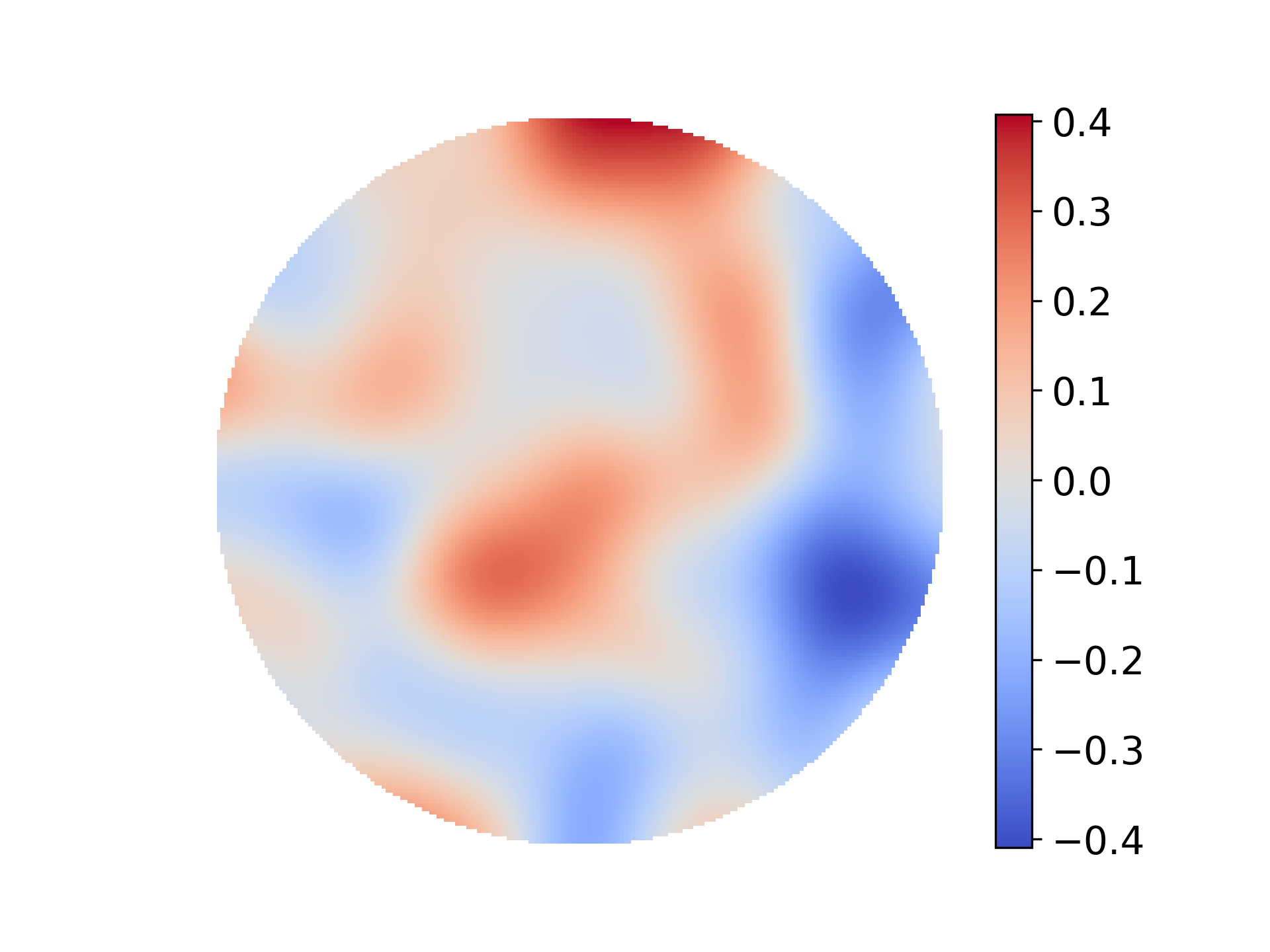}
        \caption{Snapshots at $t=0, 0.4, 0.8, 1.0$ for dynamic boundary conditions with \( M_b=M_s=2 \).}
        \label{fig:static-dy-2}
    \end{subfigure}

    \begin{subfigure}[t]{1\textwidth}
        \centering
        \includegraphics[scale=0.22]{fig/ac2d/neumann-Mb5/u_pred_0.0000.png}
        \includegraphics[scale=0.22]{fig/ac2d/neumann-Mb5/u_pred_0.4000.png}
        \includegraphics[scale=0.22]{fig/ac2d/neumann-Mb5/u_pred_0.8000.png}
        \includegraphics[scale=0.22]{fig/ac2d/neumann-Mb5/u_pred_1.0000.png}
        \caption{Snapshots at $t=0, 0.4, 0.8, 1.0$ for Neumann boundary conditions with \( M_b=5 \).}
        \label{fig:static-dy-3}
    \end{subfigure}

    \begin{subfigure}[t]{1\textwidth}
        \centering
        \includegraphics[scale=0.22]{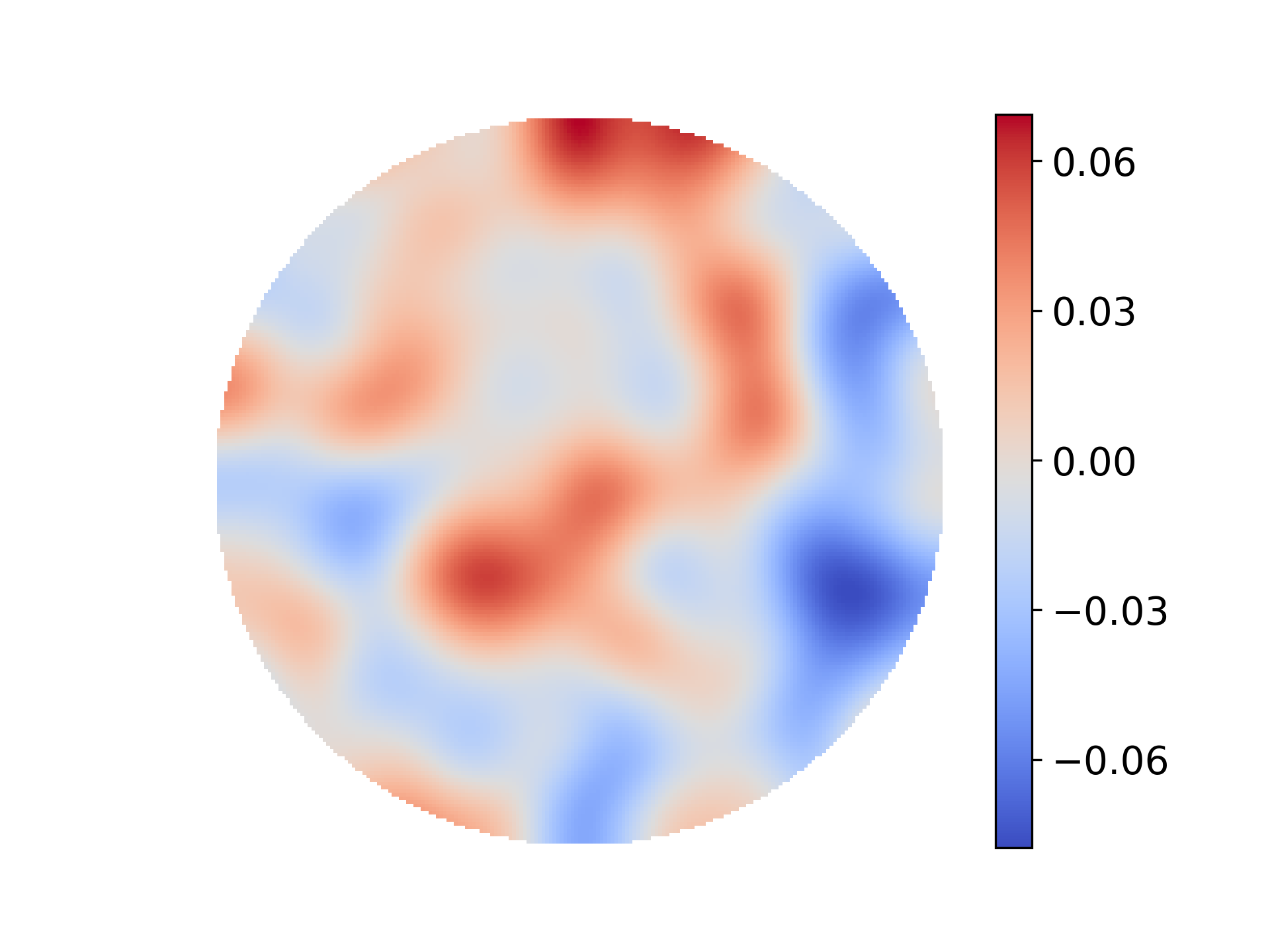}
        \includegraphics[scale=0.22]{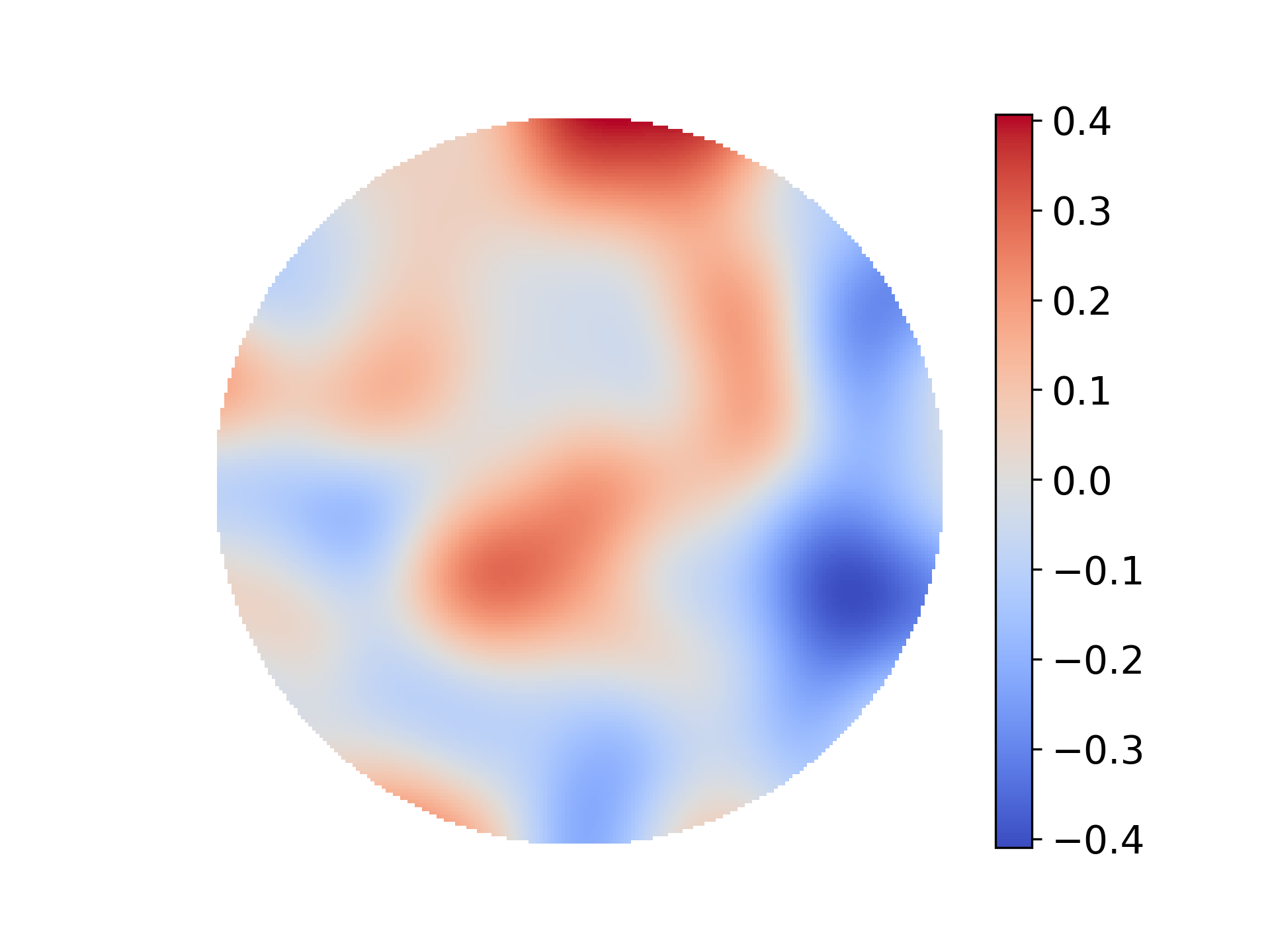}
        \includegraphics[scale=0.22]{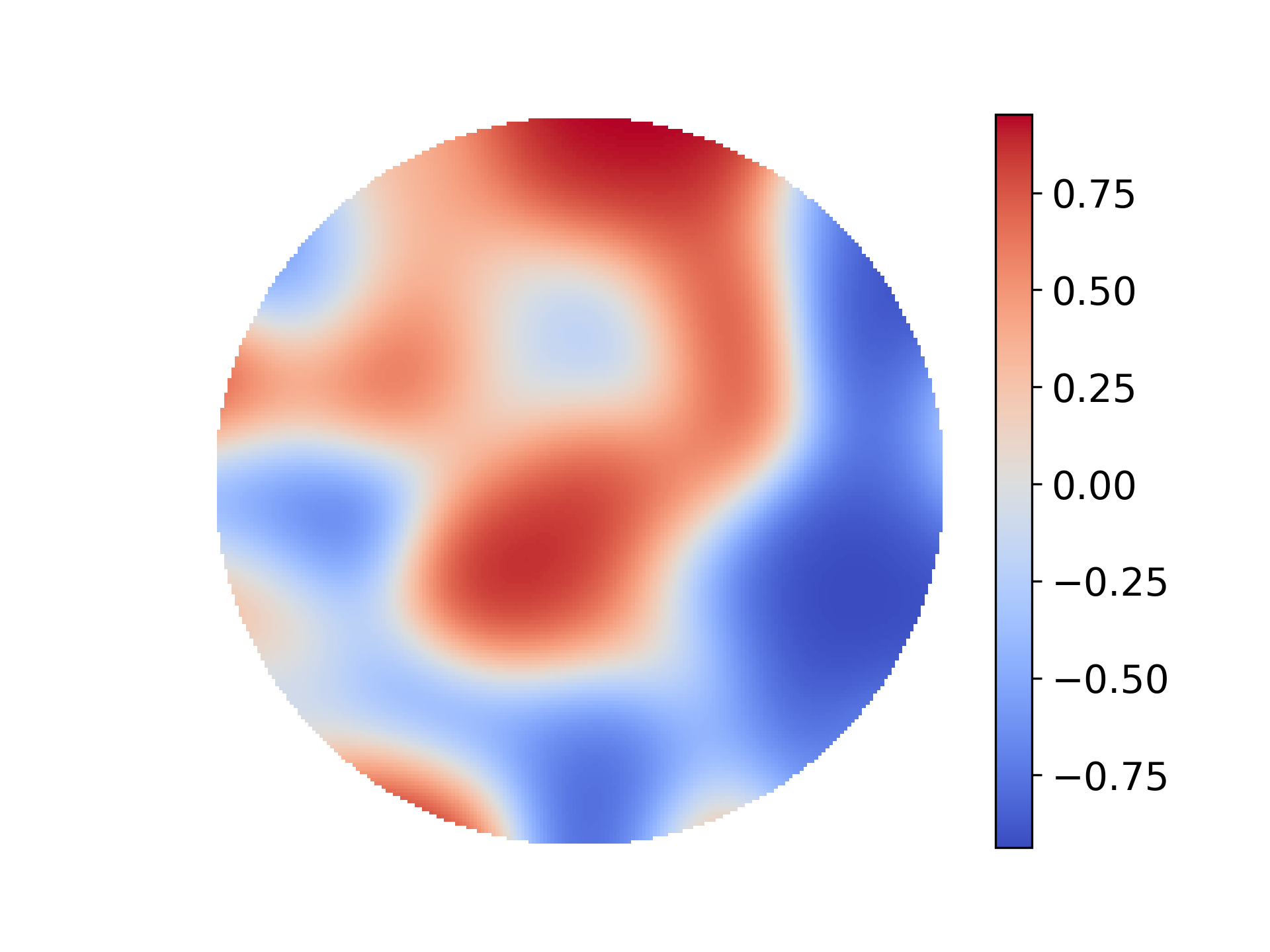}
        \includegraphics[scale=0.22]{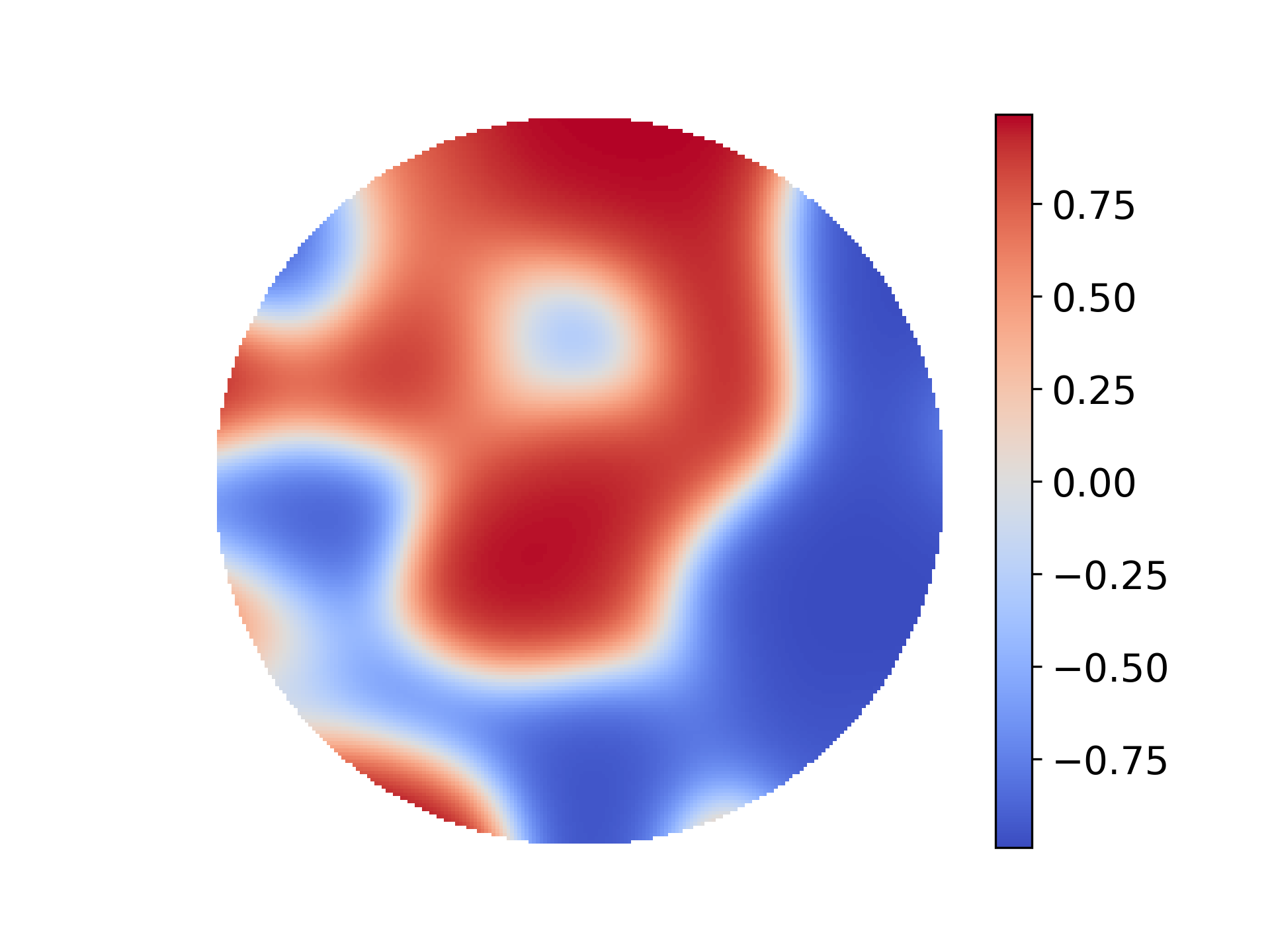}
        \caption{Snapshots at $t=0, 0.4, 0.8, 1.0$ for dynamic boundary conditions with \( M_b=M_s=5 \).}
        \label{fig:static-dy-4}
    \end{subfigure}

    \begin{subfigure}[t]{1\textwidth}
        \centering
        \includegraphics[scale=0.22]{fig/ac2d/neumann-Mb10/u_pred_0.0000.png}
        \includegraphics[scale=0.22]{fig/ac2d/neumann-Mb10/u_pred_0.4000.png}
        \includegraphics[scale=0.22]{fig/ac2d/neumann-Mb10/u_pred_0.8000.png}
        \includegraphics[scale=0.22]{fig/ac2d/neumann-Mb10/u_pred_1.0000.png}
        \caption{Snapshots at $t=0, 0.4, 0.8, 1.0$ for Neumann boundary conditions with \( M_b=10 \).}
        \label{fig:static-dy-5}
    \end{subfigure}

    \begin{subfigure}[t]{1\textwidth}
        \centering
        \includegraphics[scale=0.22]{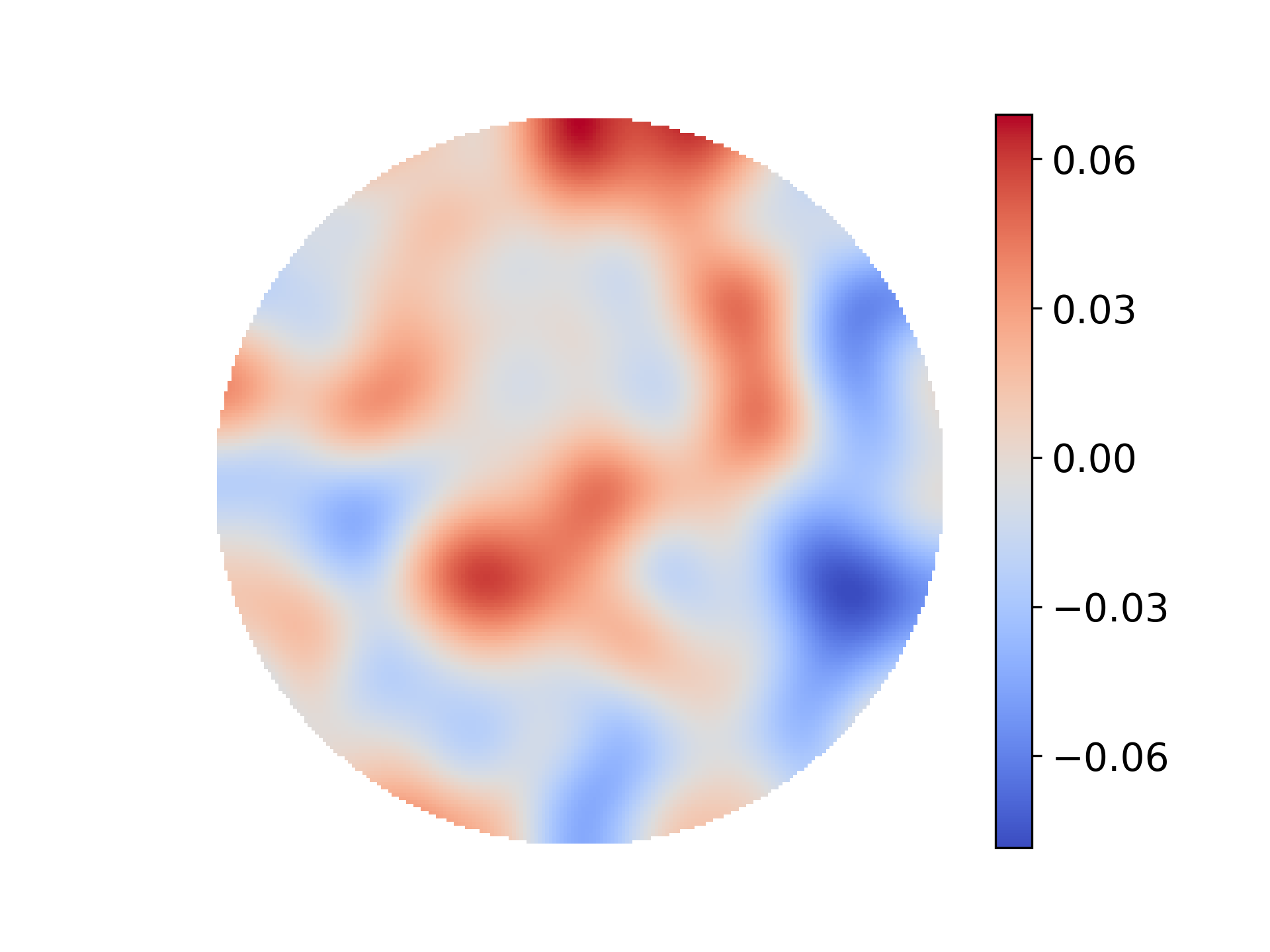}
        \includegraphics[scale=0.22]{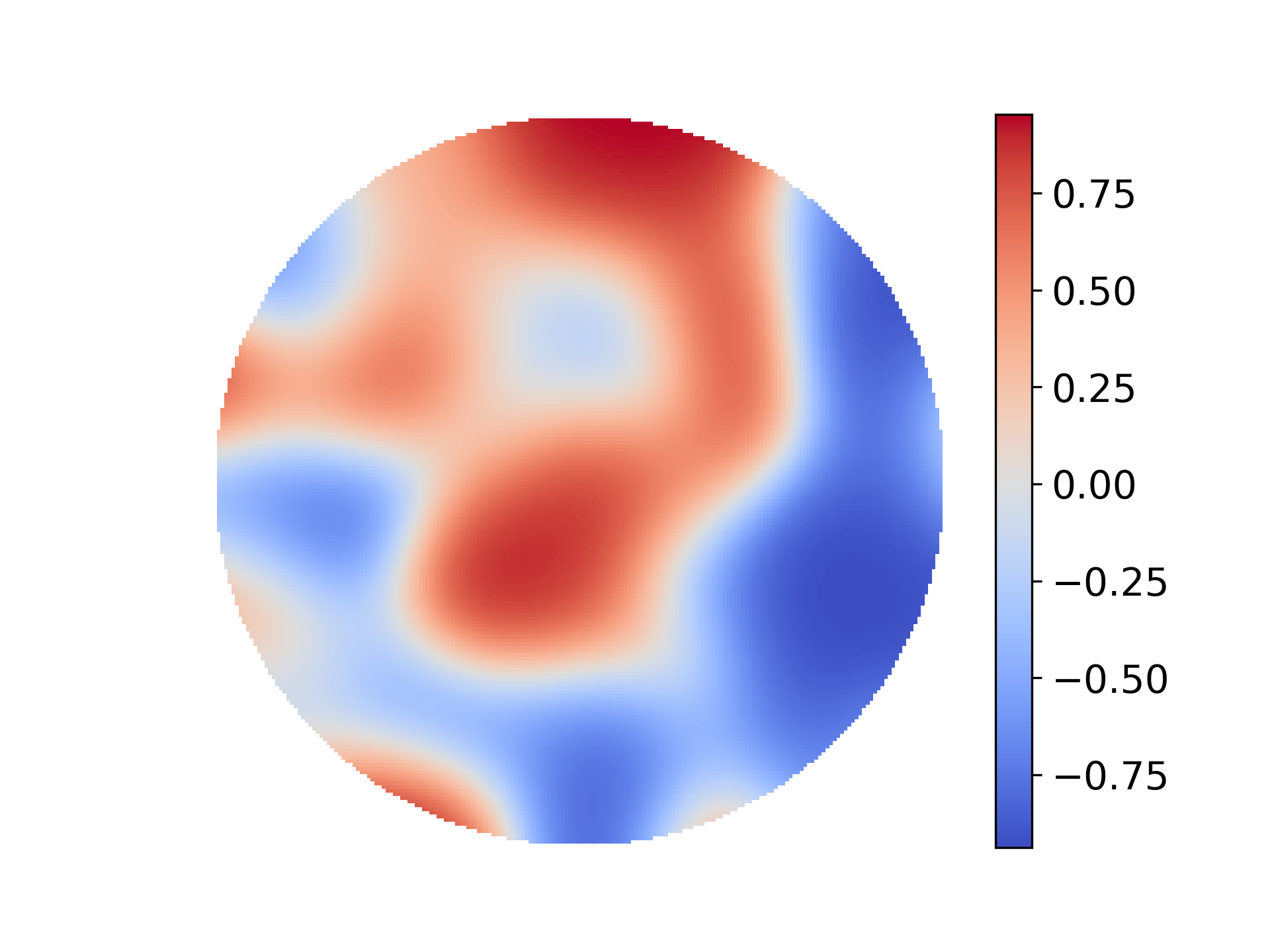}
        \includegraphics[scale=0.22]{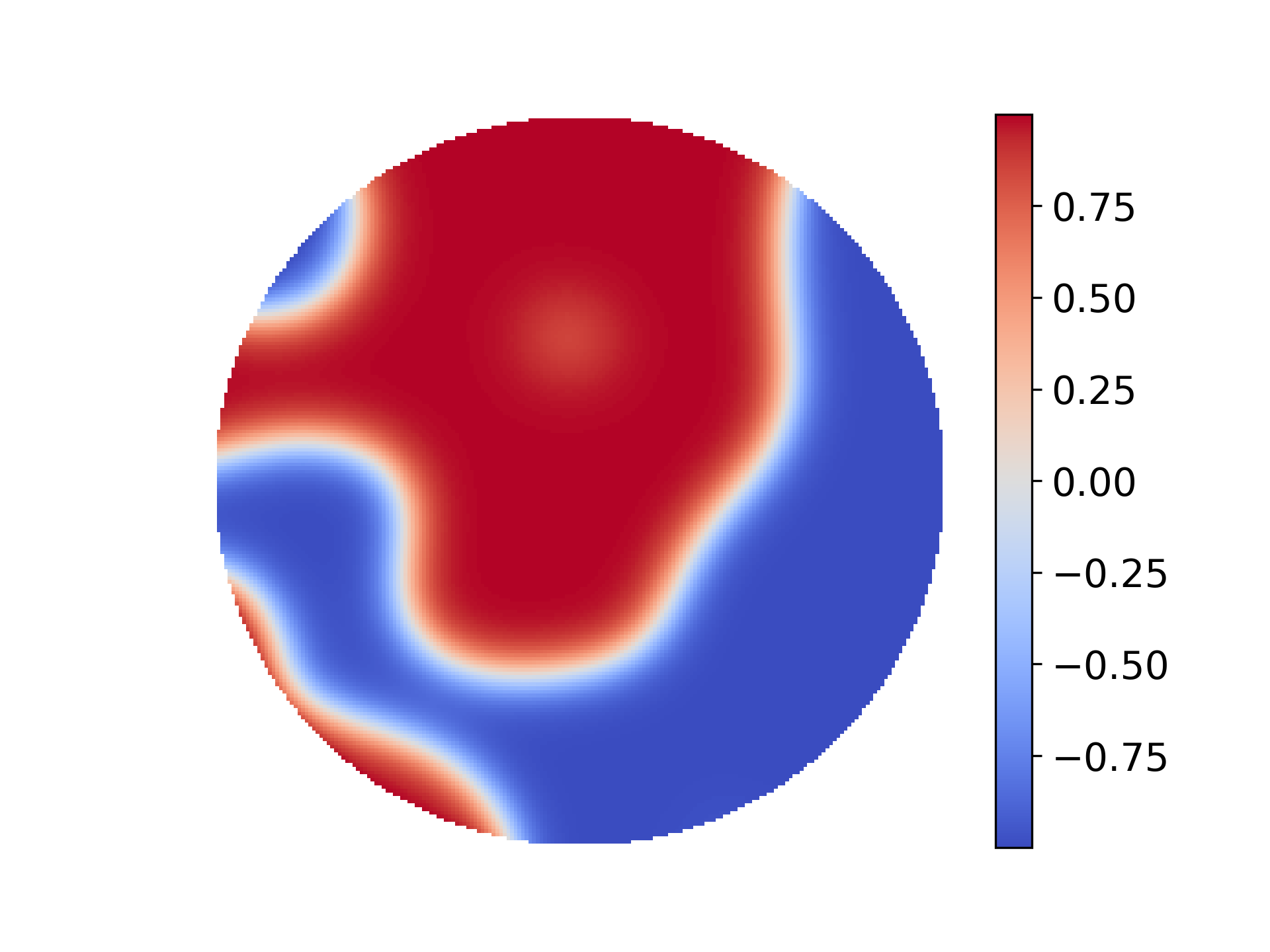}
        \includegraphics[scale=0.22]{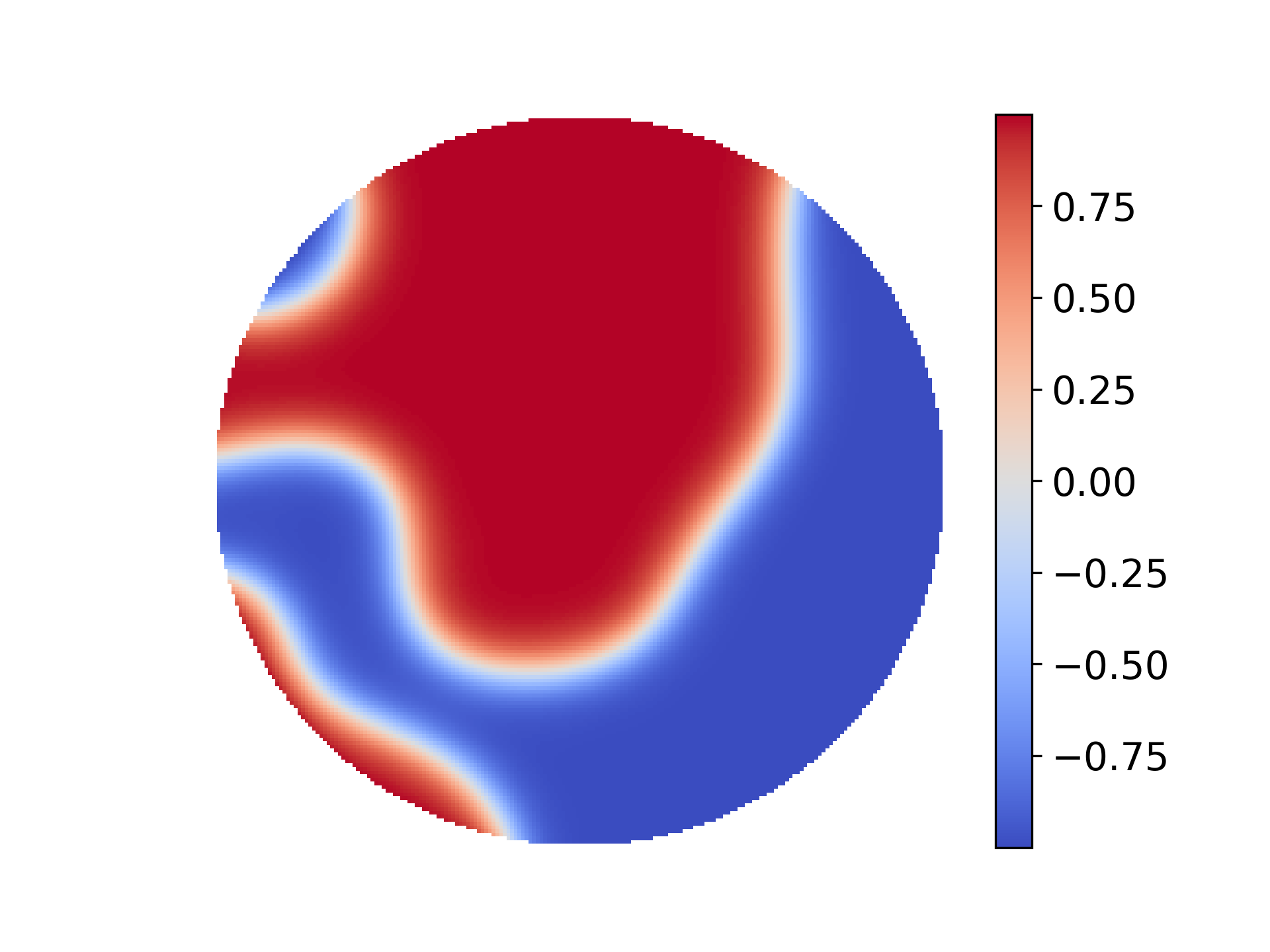}
        \caption{Snapshots at $t=0, 0.4, 0.8, 1.0$ for dynamic boundary conditions with \( M_b=M_s=10 \).}
        \label{fig:static-dy-6}
    \end{subfigure}

    \caption{Snapshots of solutions at selected time with different mobility values to showcase the difference between the solution obtained using  static Neumann boundary conditions and  that using dynamic boundary conditions. The results show that  the impact of surface dynamics enhances as the surface mobility increases.}
    \label{fig:static-dy}
\end{figure}

We next conduct a comparative analysis of static Neumann and dynamic boundary conditions by examining boundary effects on the bulk and boundary dynamics across three mobility regimes. Our investigation focuses on six scenarios:
\begin{enumerate}
\item For mobility $M_b=2$:
\begin{itemize}
\item Static case: Neumann boundary conditions (see Figure \ref{fig:static-dy-1});
\item Dynamic case: Dynamic boundary conditions with $M_b=M_s=2$ (see Figure \ref{fig:static-dy-2}).
\end{itemize}

\item For mobility $M_b=5$:
\begin{itemize}
\item Static case: Neumann boundary conditions (see Figure \ref{fig:static-dy-3});
\item Dynamic case: Dynamic boundary conditions with $M_b=M_s=5$ (see Figure \ref{fig:static-dy-4}).
\end{itemize}

\item For mobility $M_b=10$:
\begin{itemize}
\item Static case: Neumann boundary conditions (see Figure \ref{fig:static-dy-5});
\item Dynamic case: Dynamic boundary conditions with $M_b=M_s=10$ (see Figure \ref{fig:static-dy-6}).
\end{itemize}

\end{enumerate}

\begin{figure}[H]
    \centering
     \includegraphics[scale=0.45]{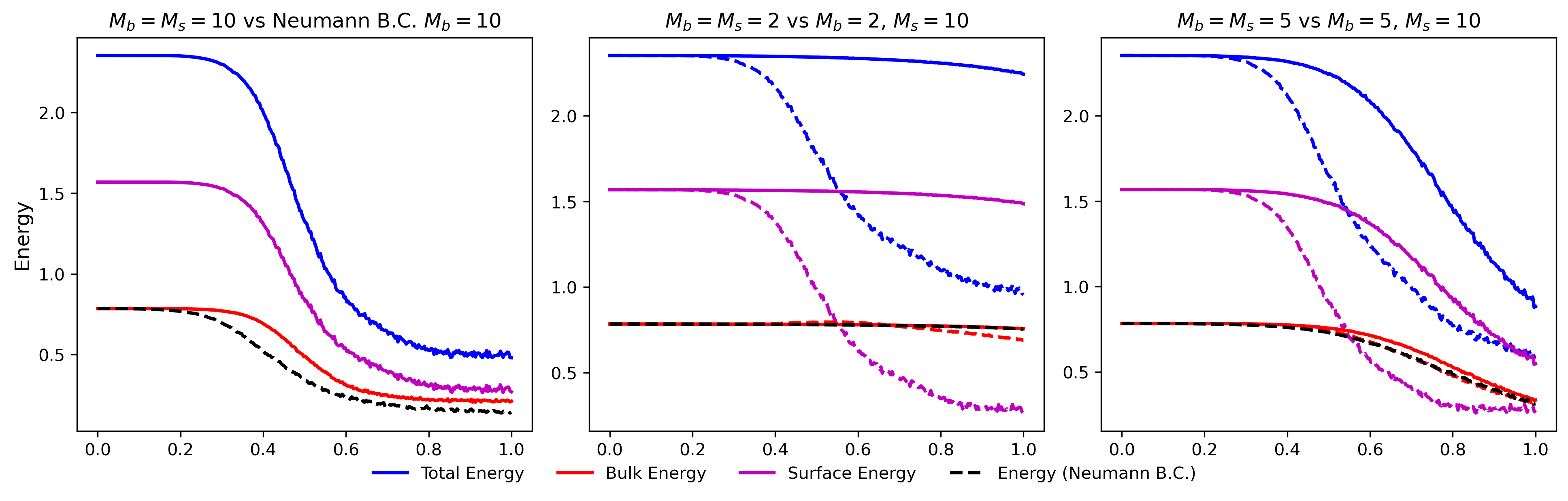}
    \caption{Plots of the total energy (blue), bulk energy (red), and surface energy (purple) with respect to time in a disk for various cases. Left: comparison between dynamic boundary conditions with $M_b=M_s=10$ and Neumann boundary conditions with $M_b=10$ (dark dash line). Middle: Dynamic boundary conditions with $M_b=M_s=2$ (solid lines) and $M_b=2, M_s=10$ (dash lines) v.s. Neumann boundary conditions with $M_b=2$ (dark dash line). Right: Dynamic boundary conditions with $M_b=M_s=5$ (solid lines) and $M_b=5, M_s=10$ (dash lines) v.s. Neumann boundary conditions with $M_b=5$ (dark dash line). The surface energy driven dissipation dominates the dynamics near the boundary. }
    \label{fig:energy-t}
\end{figure}

The comparison reveals some noticeable differences existing between static and dynamic boundary conditions:
\begin{itemize}
\item Boundary behavior: Distinct patterns emerge at the domain boundaries.
\item Near-boundary bulk dynamics: The influence of boundary conditions extends noticeably into the bulk region.
\item Energy dissipation: As shown in Figure \ref{fig:energy-t}, while all scenarios exhibit energy dissipation, the Neumann boundary condition (black dashed line in the left subplot) leads to a faster bulk energy decay compared to the case with a dynamic boundary condition (red solid line in the left subplot).
\end{itemize}


\subsection{Allen-Cahn Equation with Dynamic Boundary Conditions in 2D}\label{subsec:AC_dy}

\begin{figure}[H]
    \centering
    \begin{subfigure}[t]{1\textwidth}
        \centering
        \includegraphics[scale=0.22]{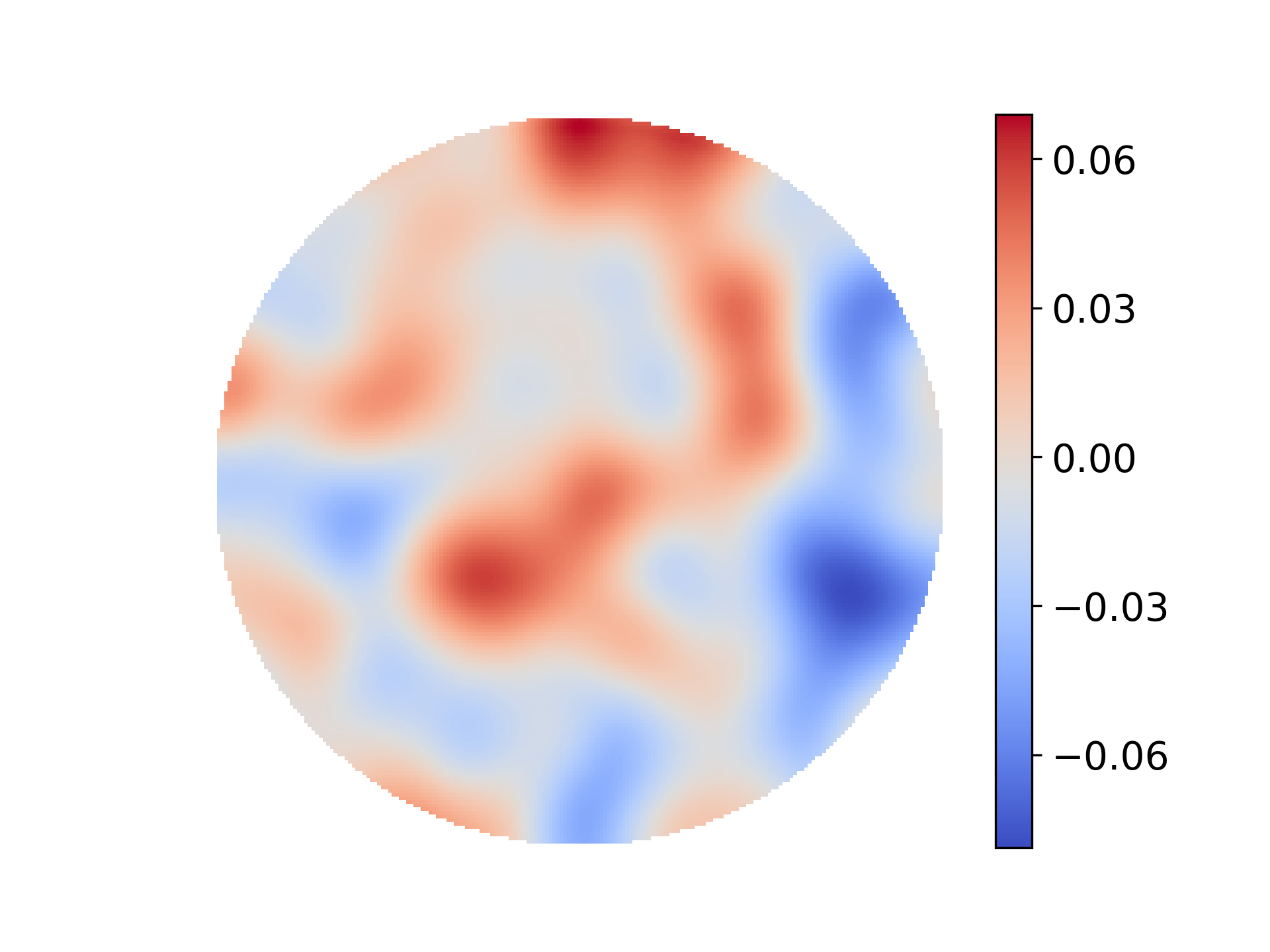}
        \includegraphics[scale=0.22]{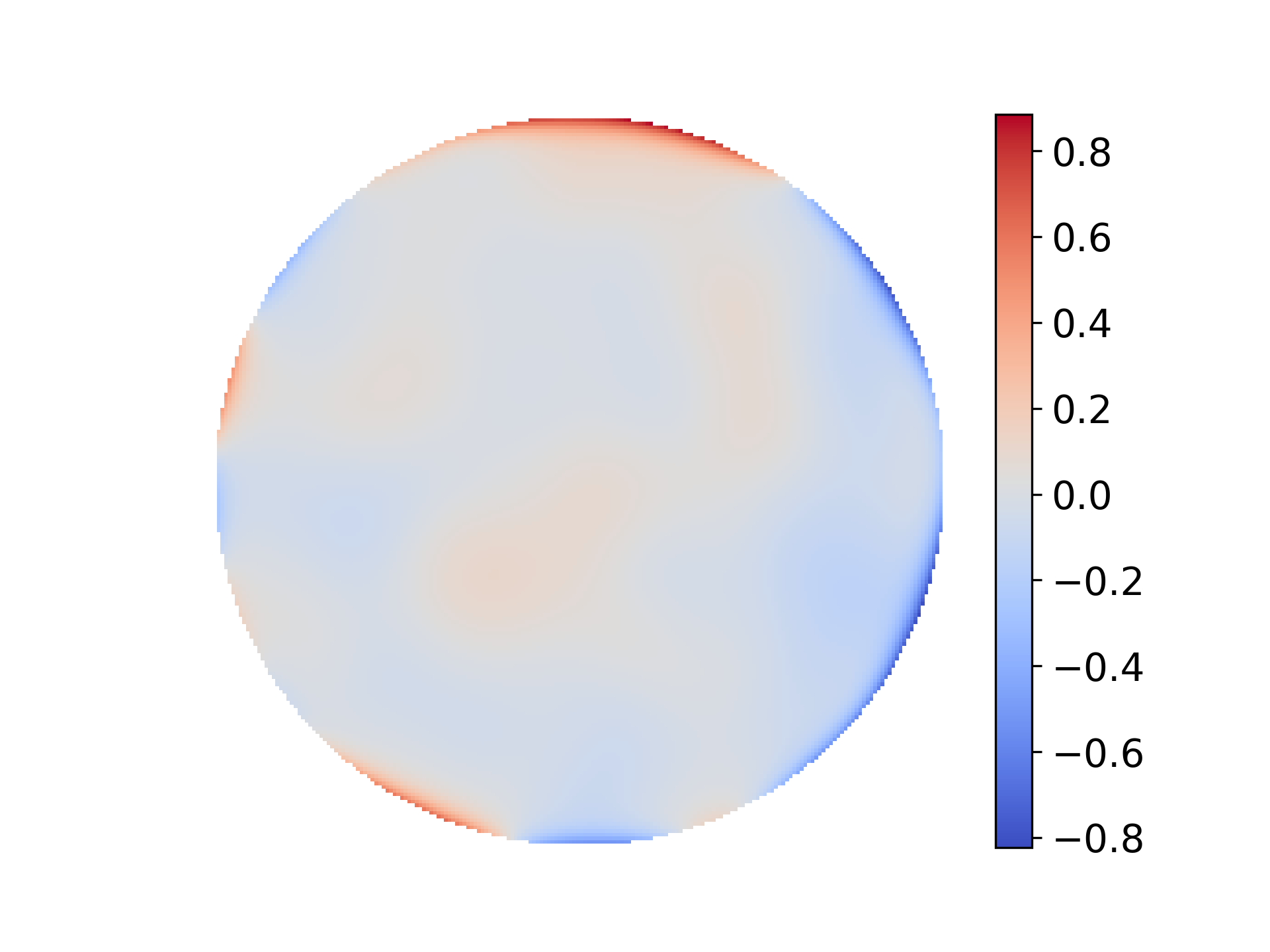}
        \includegraphics[scale=0.22]{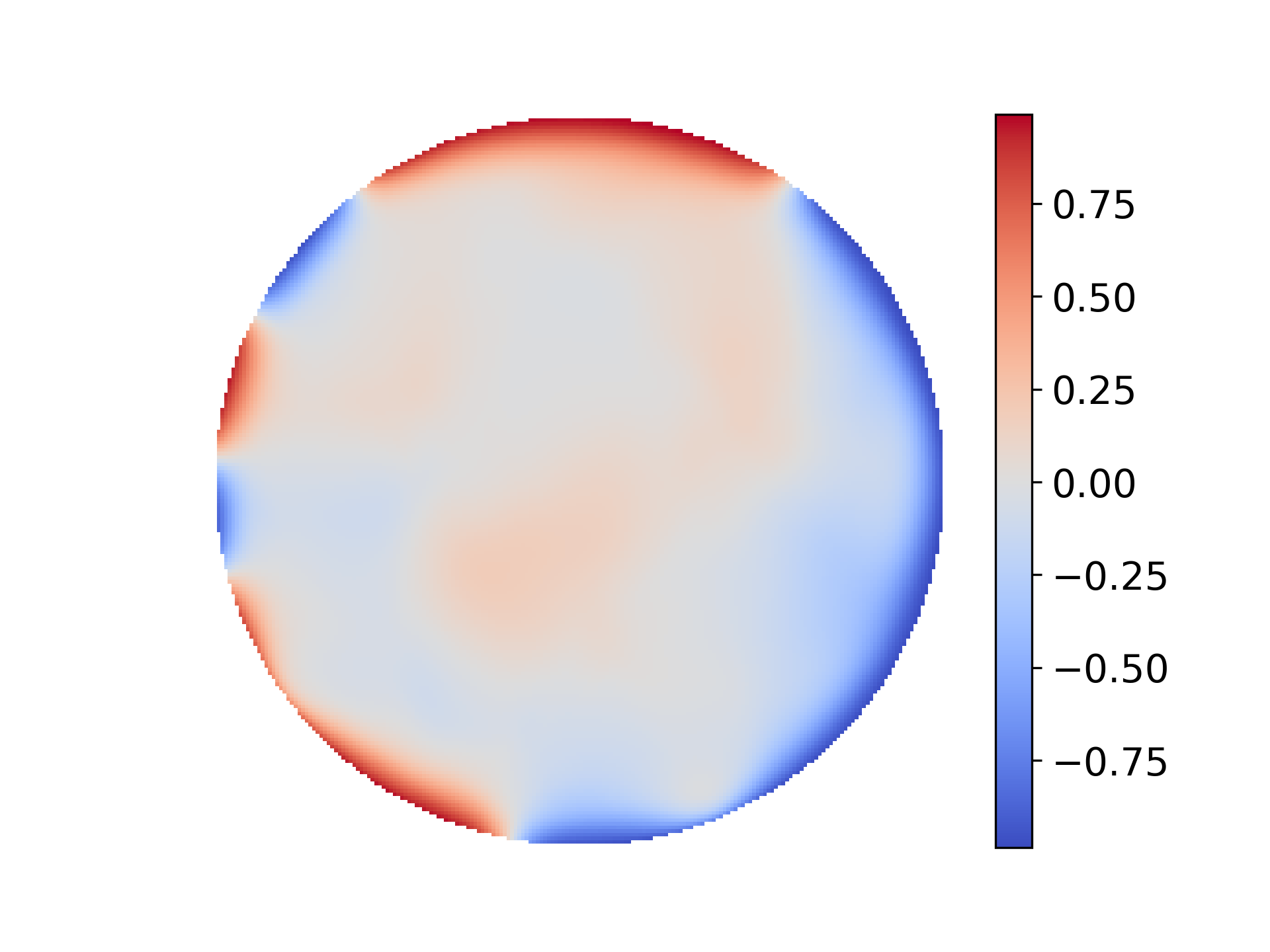}
        \includegraphics[scale=0.22]{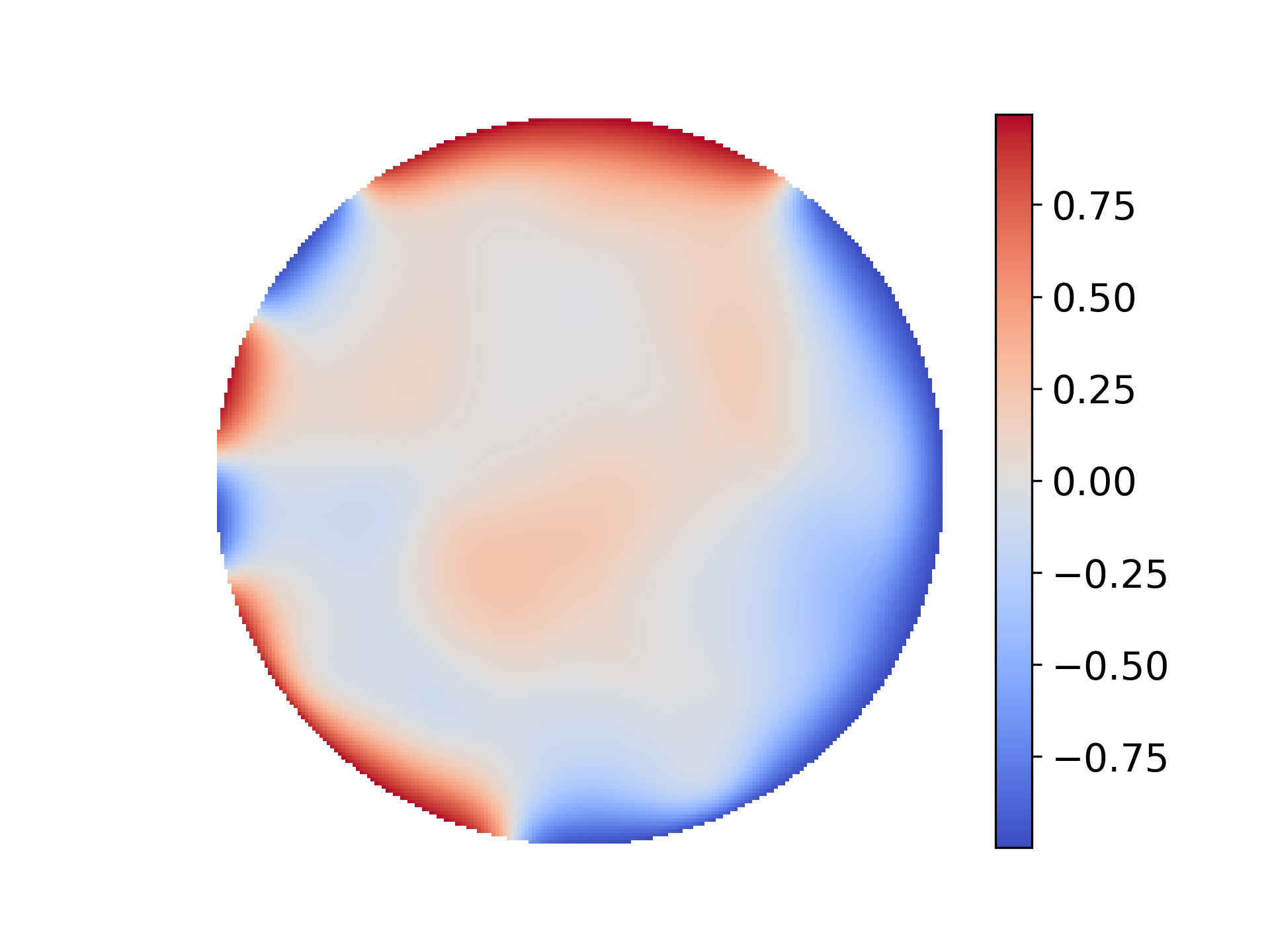}
        \caption{\( M_b=2 \) and \( M_s=10 \).}
    \end{subfigure}

    \begin{subfigure}[t]{1\textwidth}
        \centering
        \includegraphics[scale=0.22]{fig/ac2d/dy-Mb2-Ms2/pred_u_0.0000.png}
        \includegraphics[scale=0.22]{fig/ac2d/dy-Mb2-Ms2/pred_u_0.4000.png}
        \includegraphics[scale=0.22]{fig/ac2d/dy-Mb2-Ms2/pred_u_0.8000.png}
        \includegraphics[scale=0.22]{fig/ac2d/dy-Mb2-Ms2/pred_u_1.0000.png}
        \caption{\( M_b=2 \) and \( M_s=2 \).}
    \end{subfigure}

    \begin{subfigure}[t]{1\textwidth}
        \centering
        \includegraphics[scale=0.22]{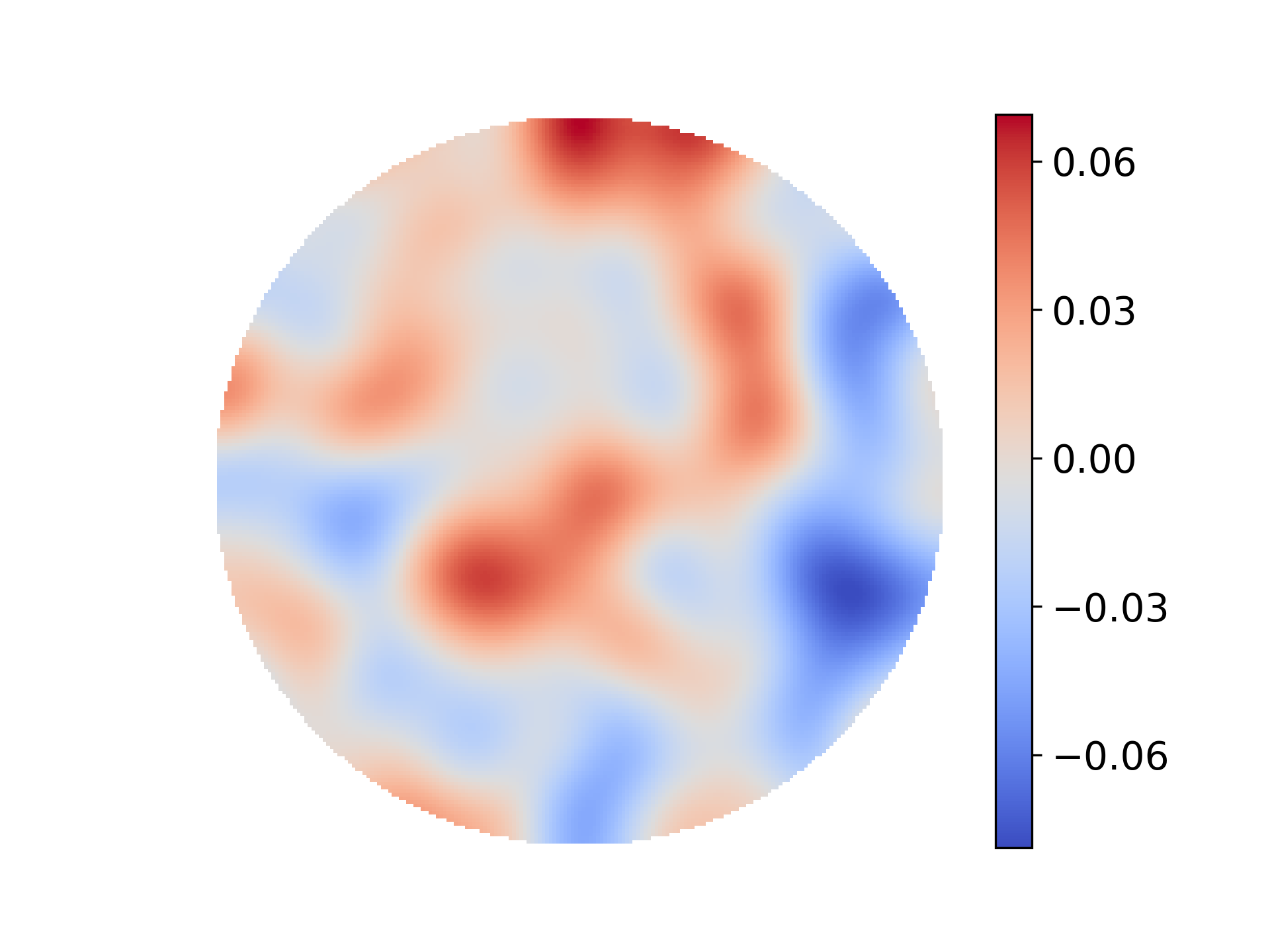}
        \includegraphics[scale=0.22]{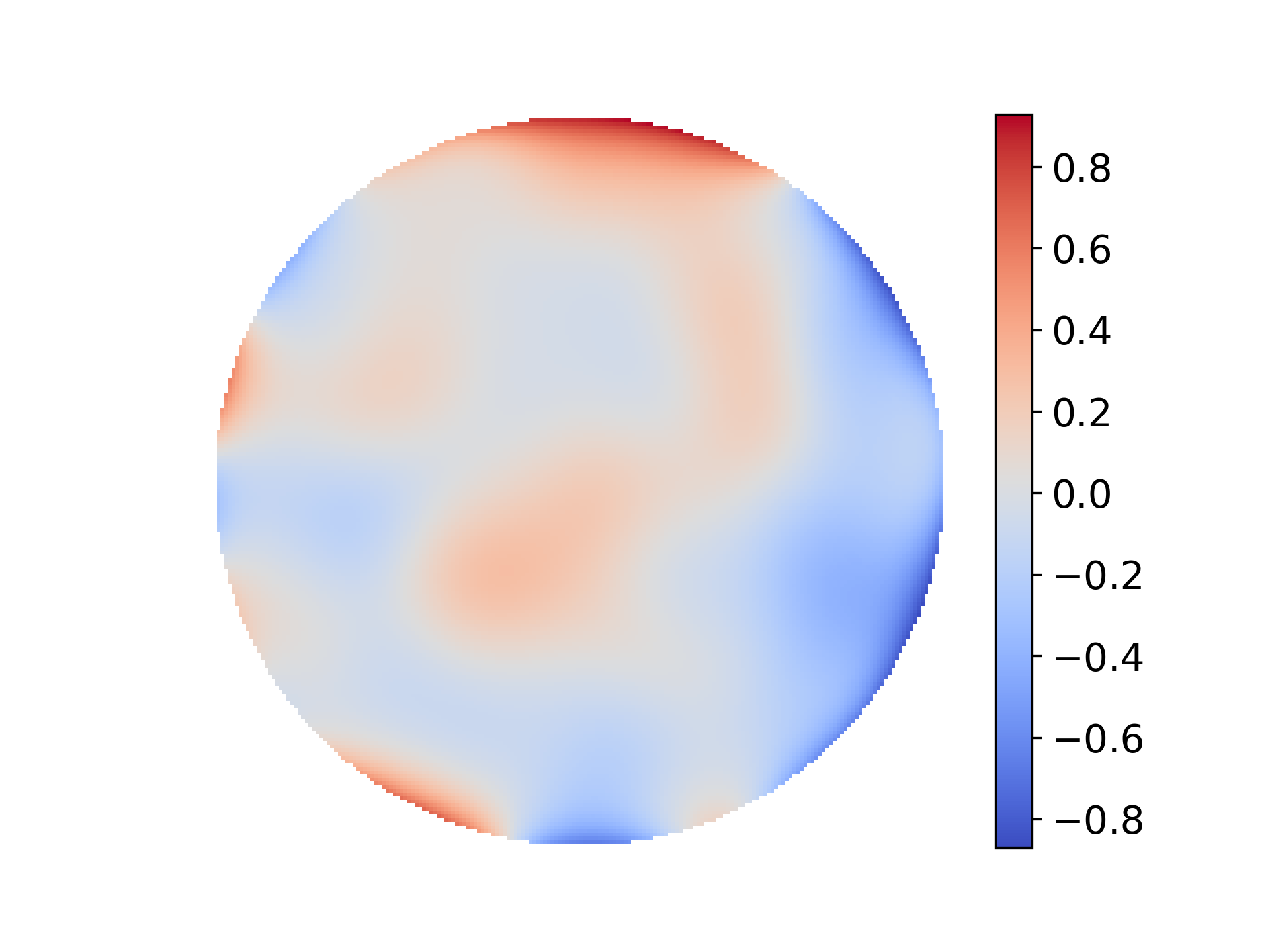}
        \includegraphics[scale=0.22]{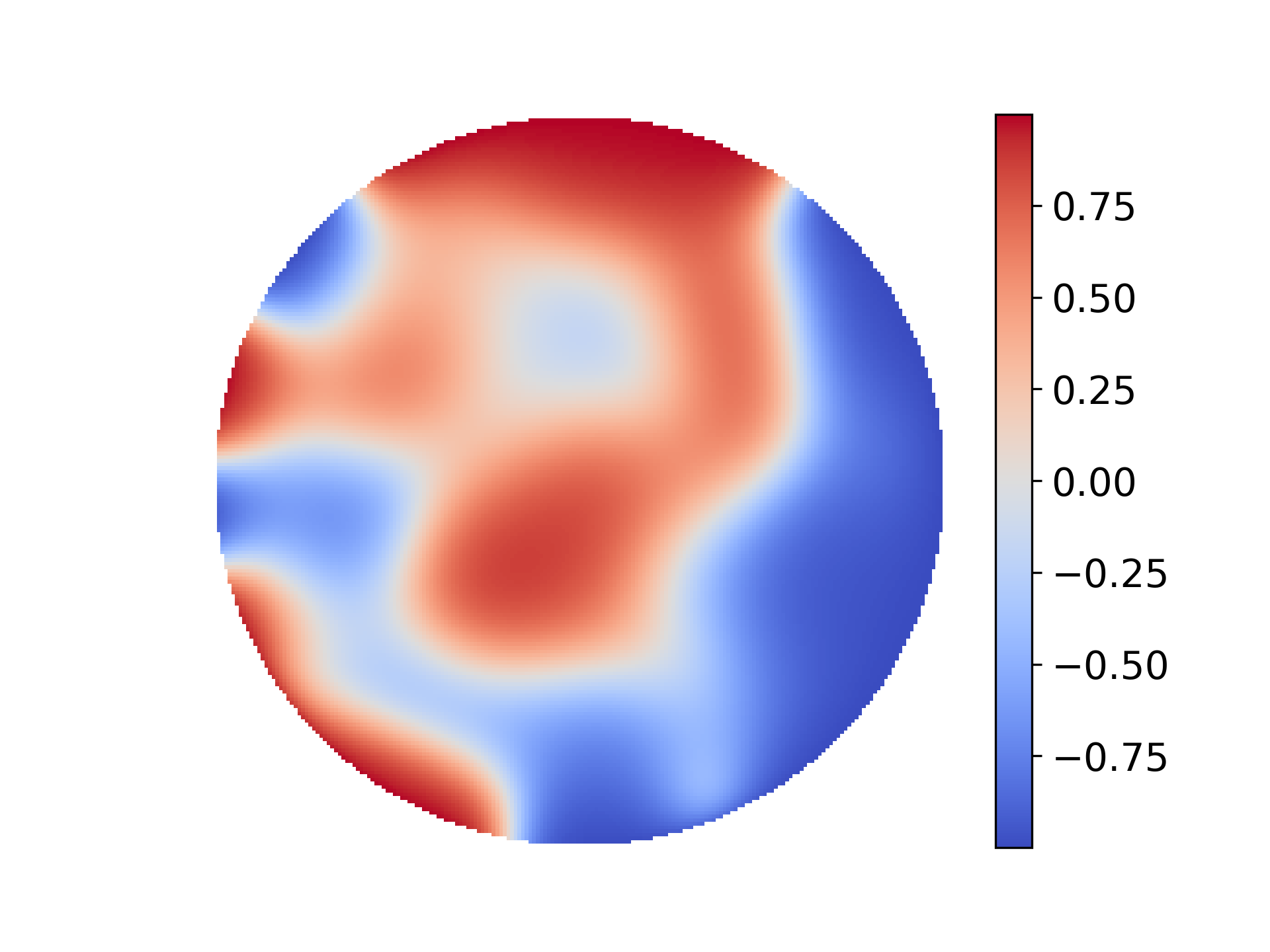}
        \includegraphics[scale=0.22]{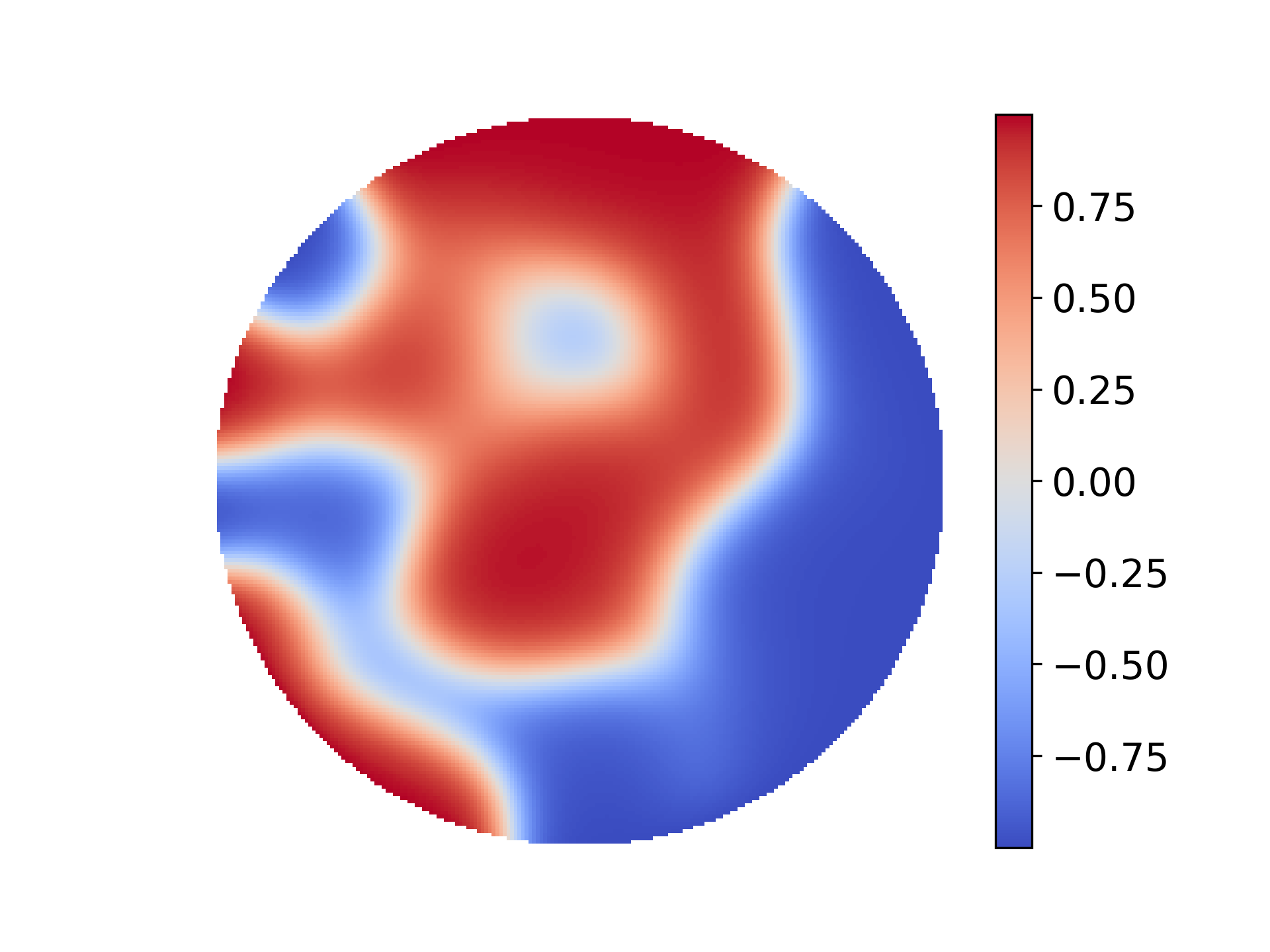}
        \caption{\( M_b=5 \) and \( M_s=10 \).}
    \end{subfigure}

    \begin{subfigure}[t]{1\textwidth}
        \centering
        \includegraphics[scale=0.22]{fig/ac2d/dy-Mb5-Ms5/pred_u_0.0000.png}
        \includegraphics[scale=0.22]{fig/ac2d/dy-Mb5-Ms5/pred_u_0.4000.png}
        \includegraphics[scale=0.22]{fig/ac2d/dy-Mb5-Ms5/pred_u_0.8000.png}
        \includegraphics[scale=0.22]{fig/ac2d/dy-Mb5-Ms5/pred_u_1.0000.png}
        \caption{\( M_b=5 \) and \( M_s=5 \).}
    \end{subfigure}

    \caption{Snapshots of solutions at selected time at $t=0, 0.4, 0.8, 1.0$ with various mobility values in the case of dynamic boundary conditions. The results show that the impact of dynamic boundary conditions is noticeable near the boundary while the region of influence expands as the surface mobility enhances. }
    \label{fig:dy-dy}
\end{figure}

\begin{figure}[H]
    \centering
    \begin{subfigure}[t]{1\textwidth}
        \centering
        \includegraphics[scale=0.22]{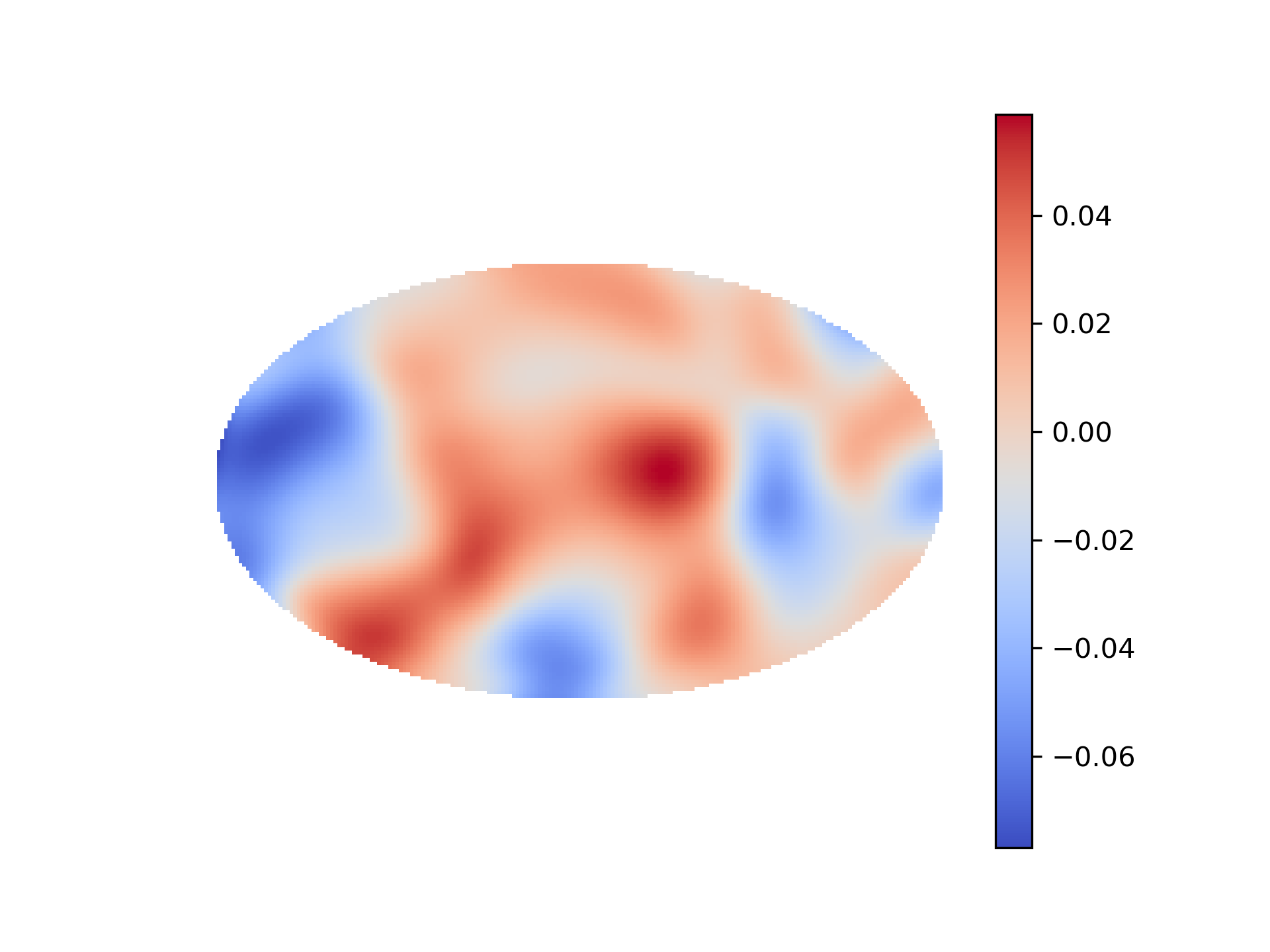}
        \includegraphics[scale=0.22]{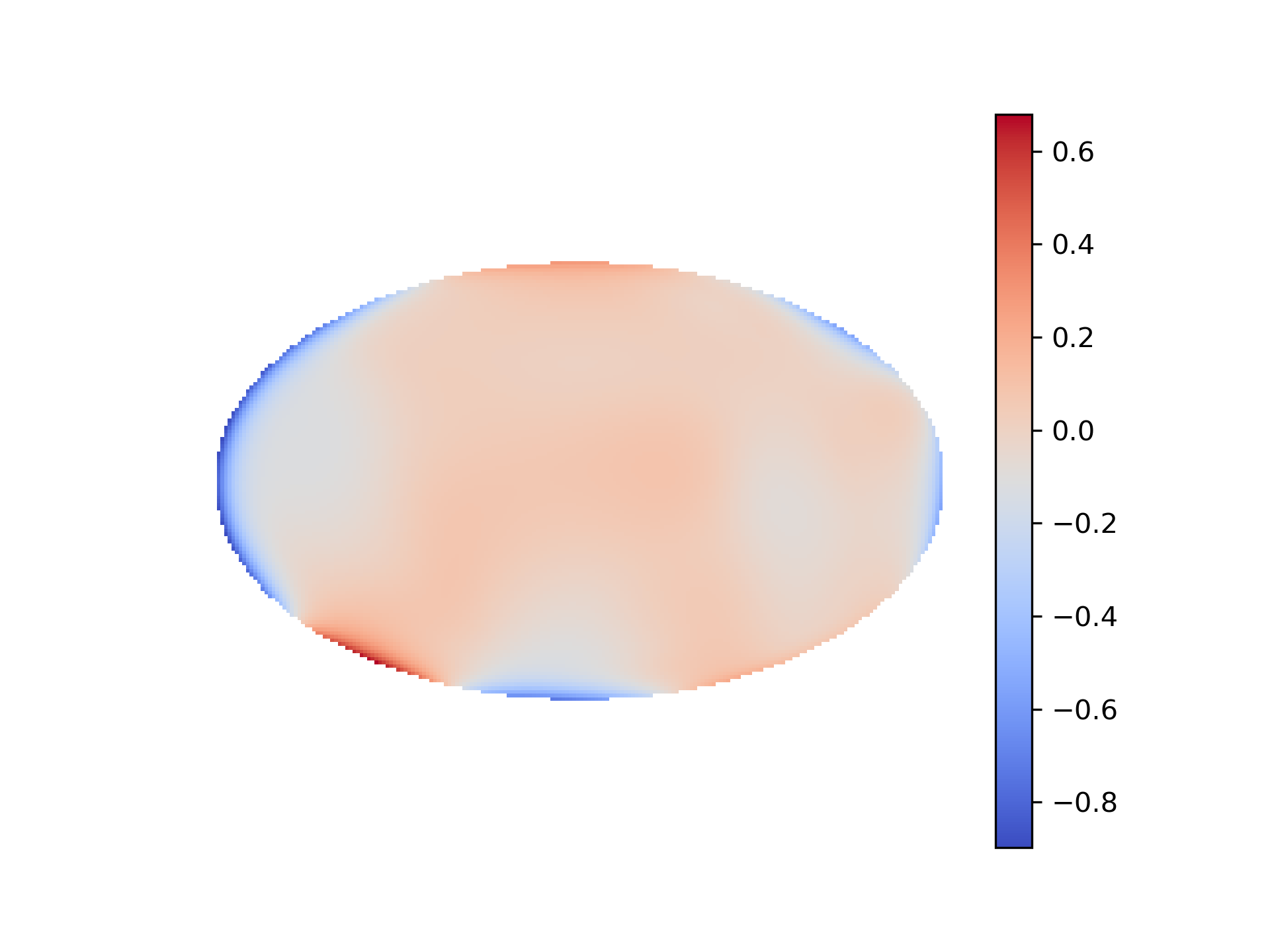}
        \includegraphics[scale=0.22]{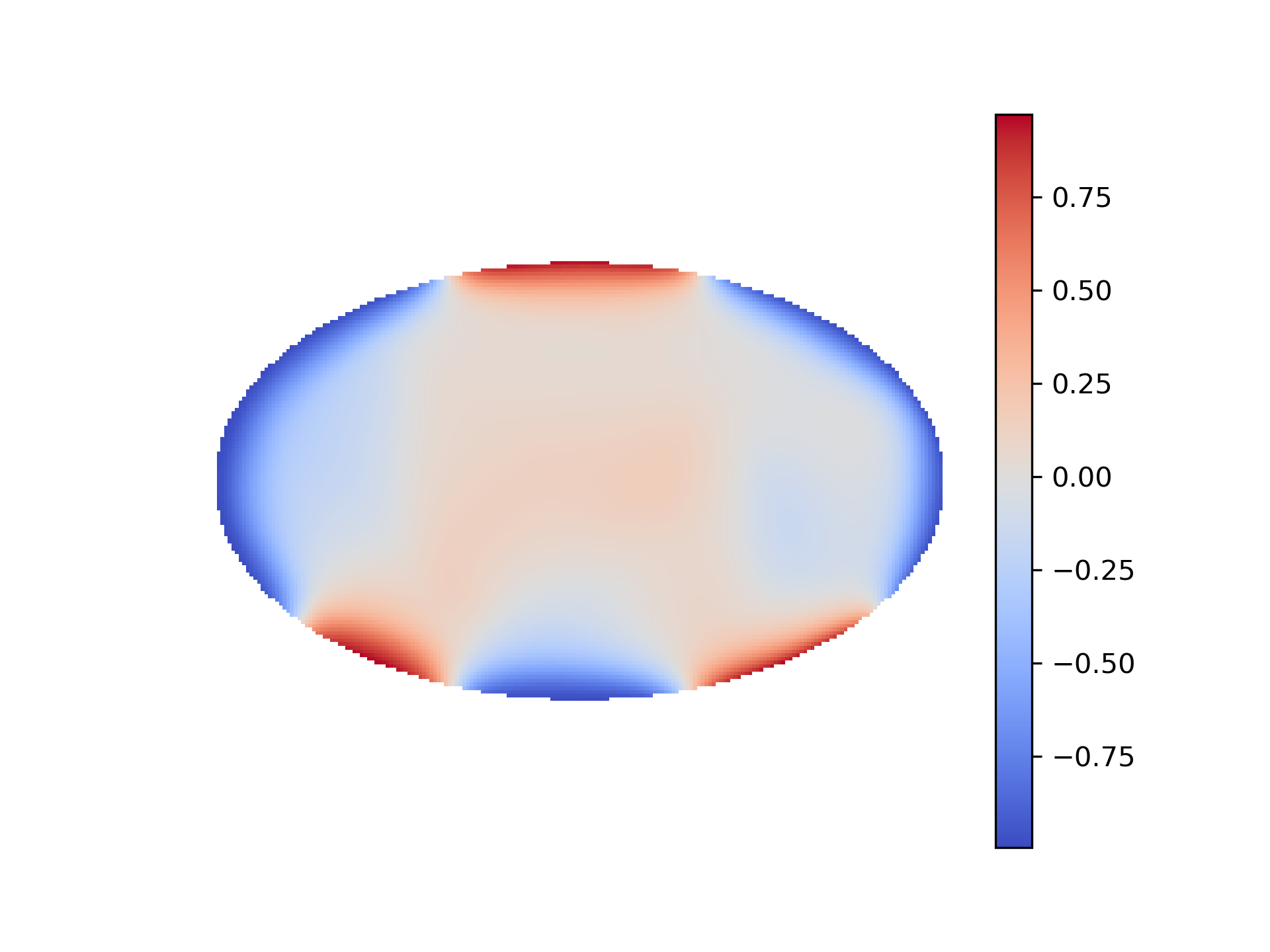}
        \includegraphics[scale=0.22]{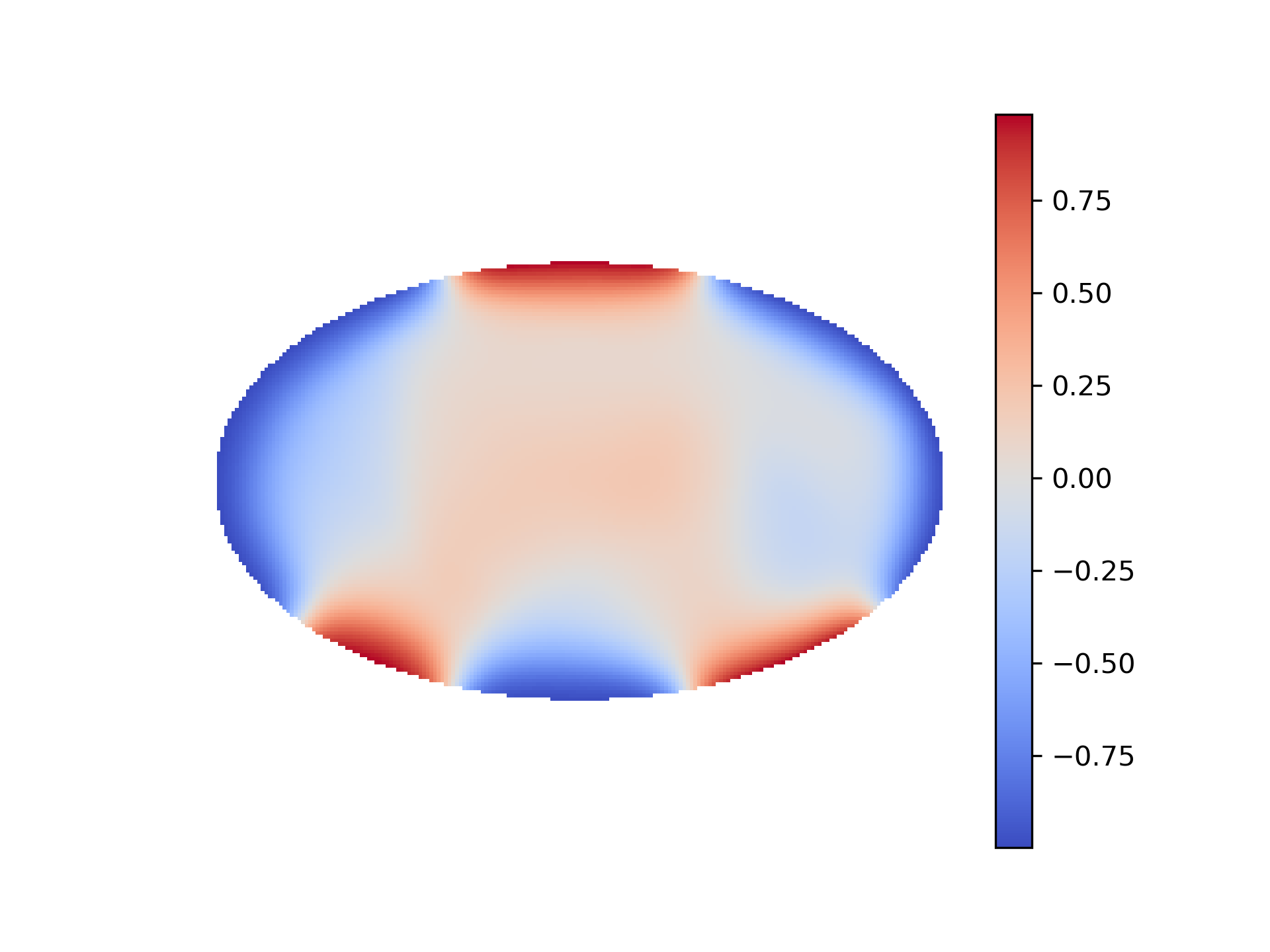}
        \caption{\( M_b=2 \) and \( M_s=10 \).}
    \end{subfigure}

    \begin{subfigure}[t]{1\textwidth}
        \centering
        \includegraphics[scale=0.22]{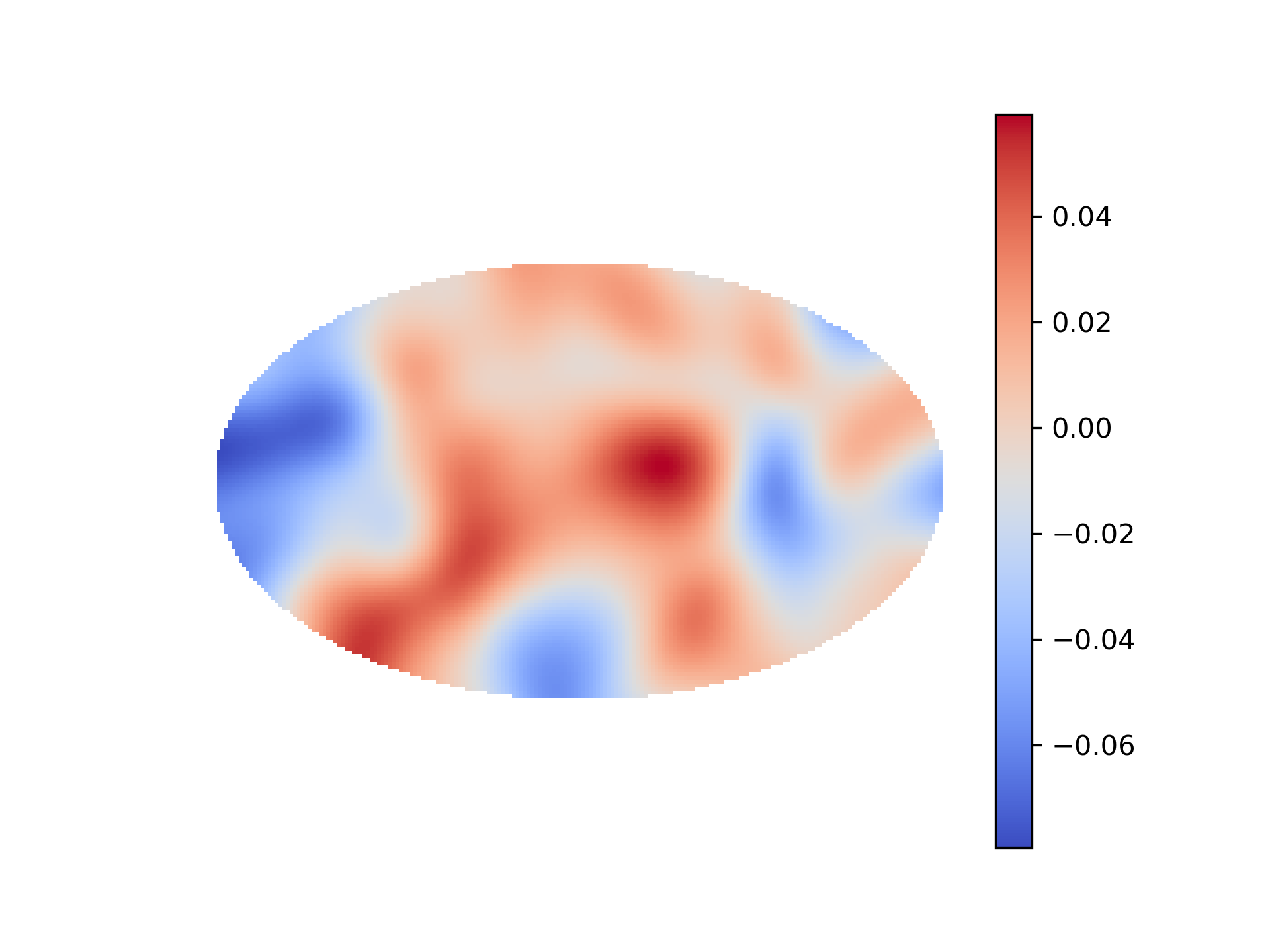}
        \includegraphics[scale=0.22]{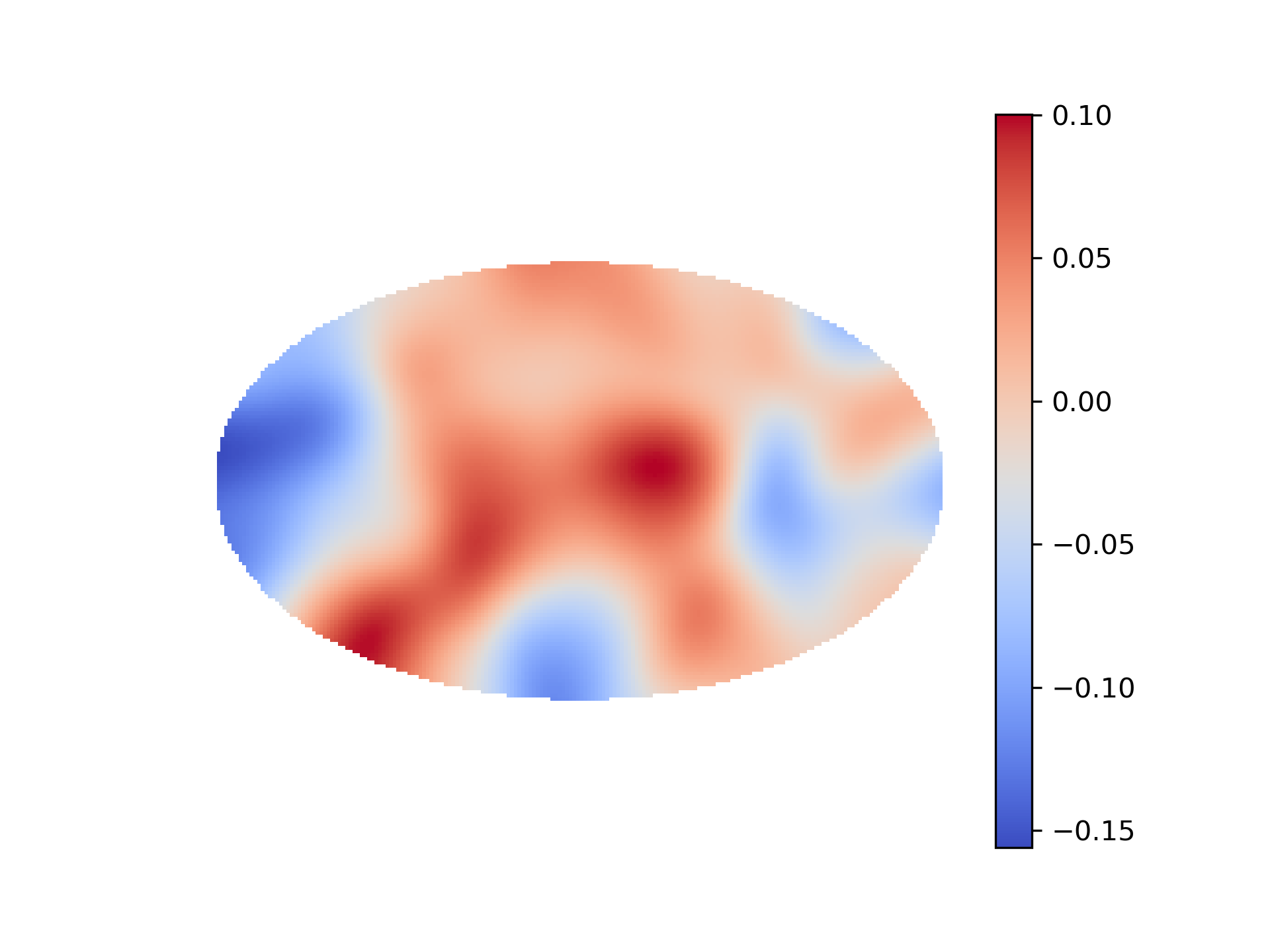}
        \includegraphics[scale=0.22]{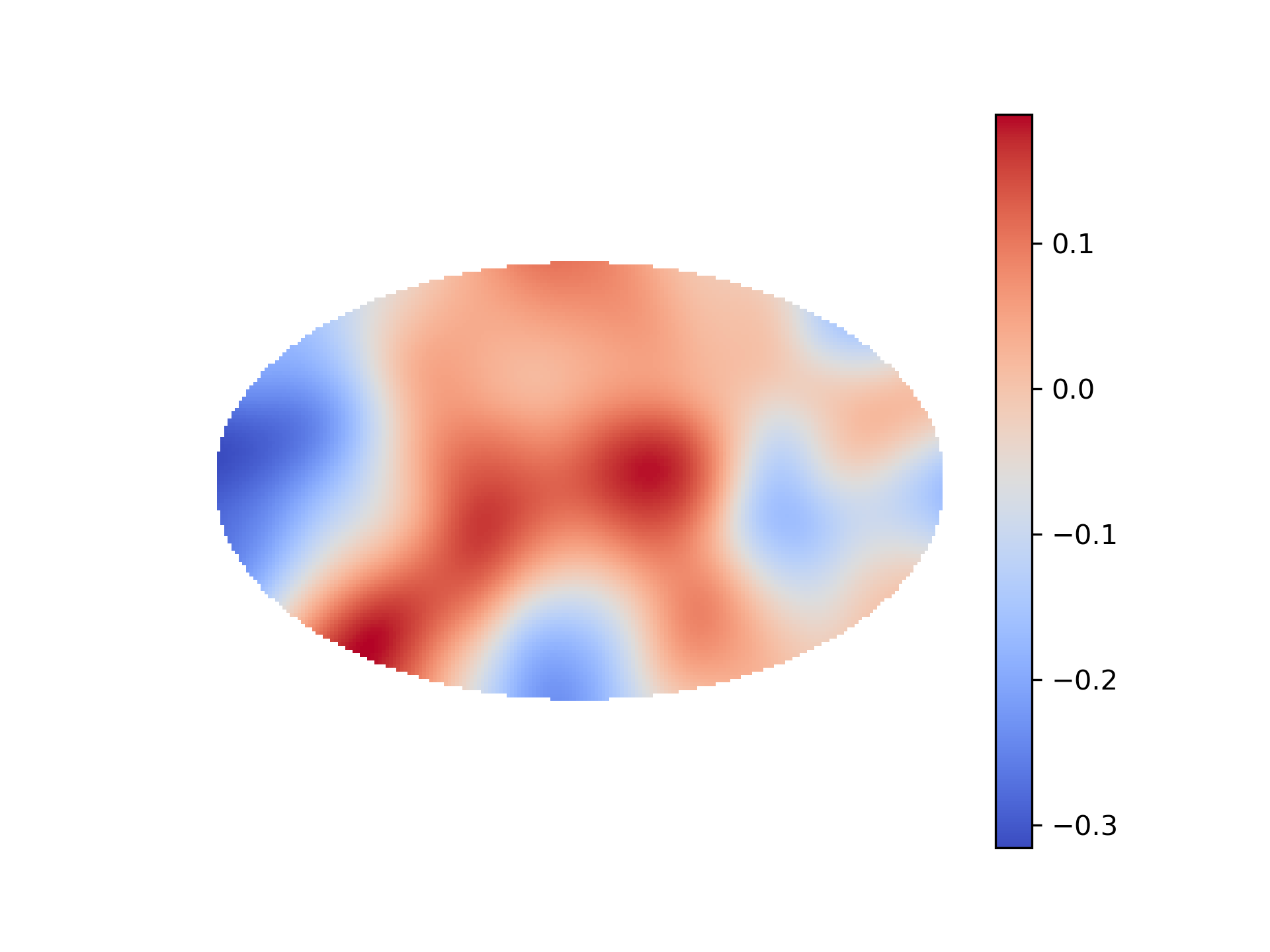}
        \includegraphics[scale=0.22]{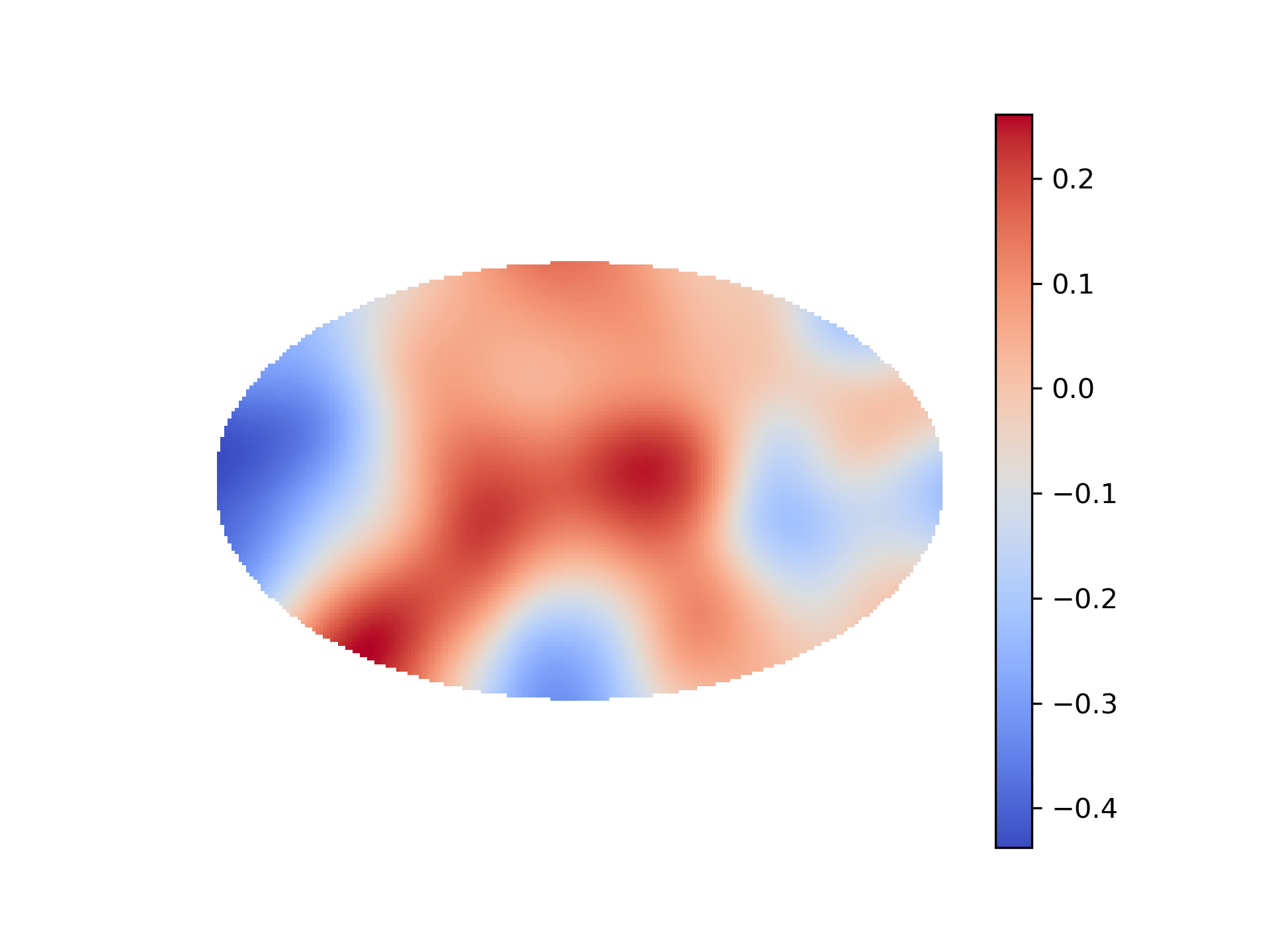}
        \caption{\( M_b=2 \) and \( M_s=2 \).}
    \end{subfigure}

    \begin{subfigure}[t]{1\textwidth}
        \centering
        \includegraphics[scale=0.22]{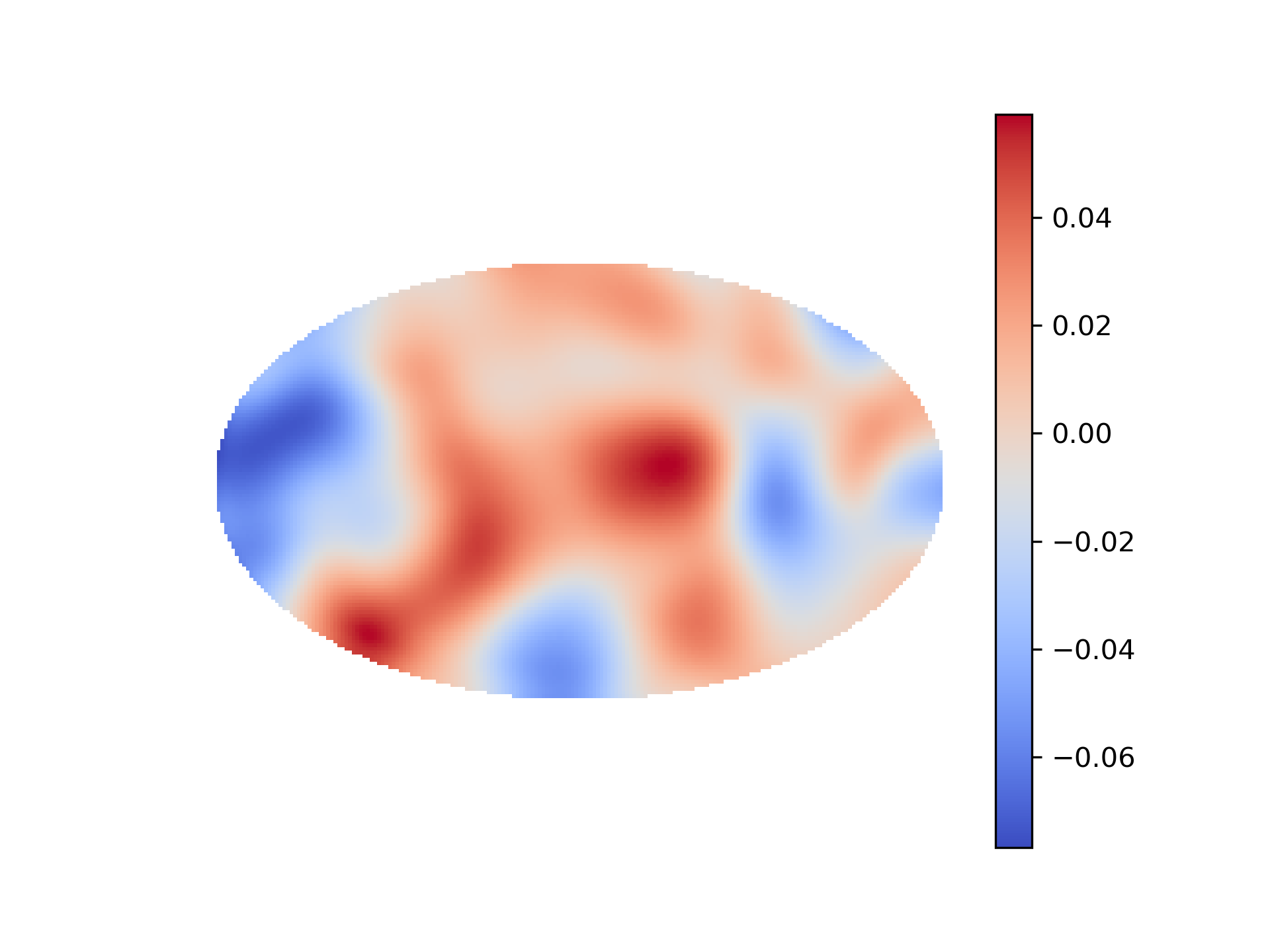}
        \includegraphics[scale=0.22]{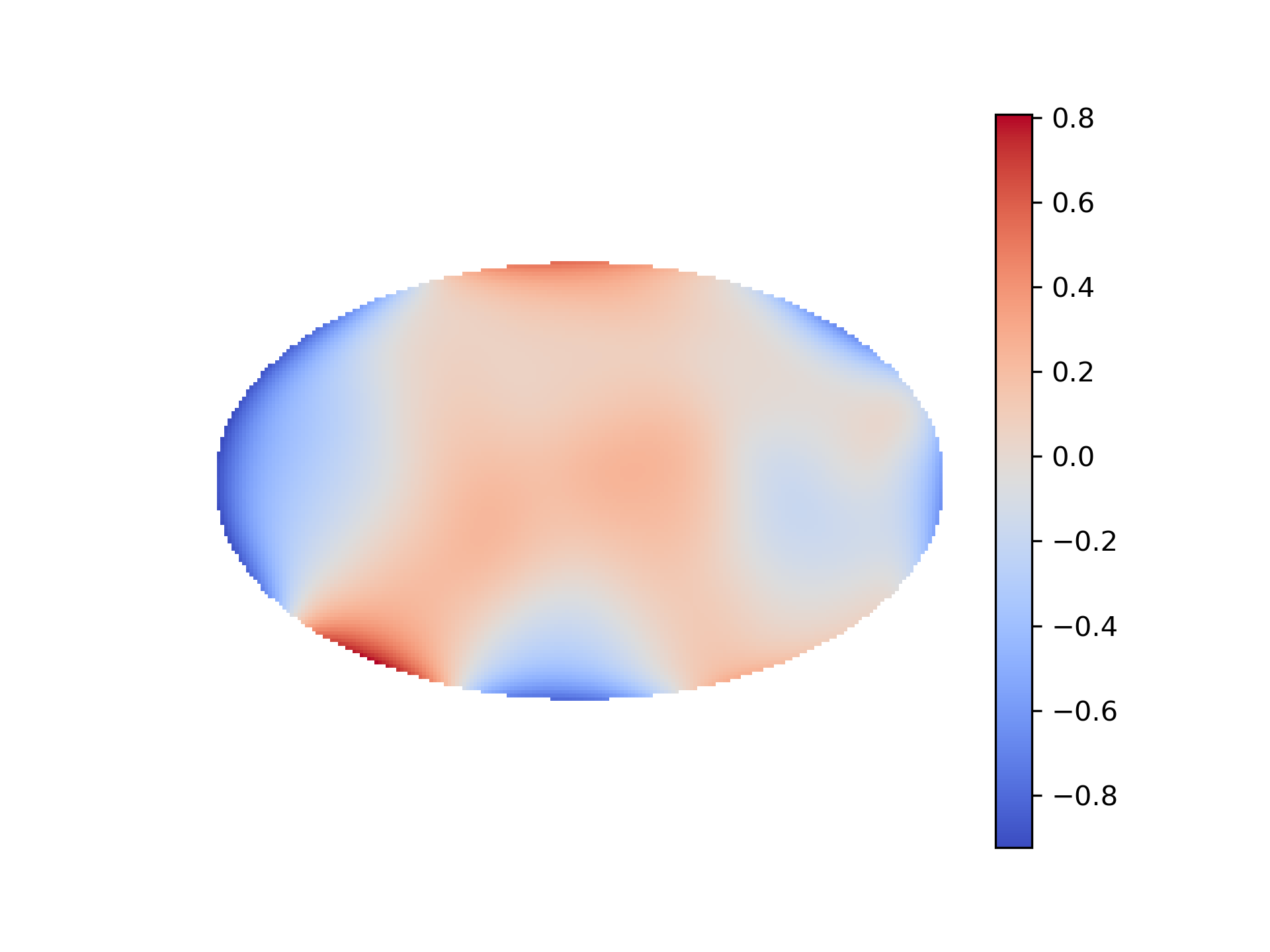}
        \includegraphics[scale=0.22]{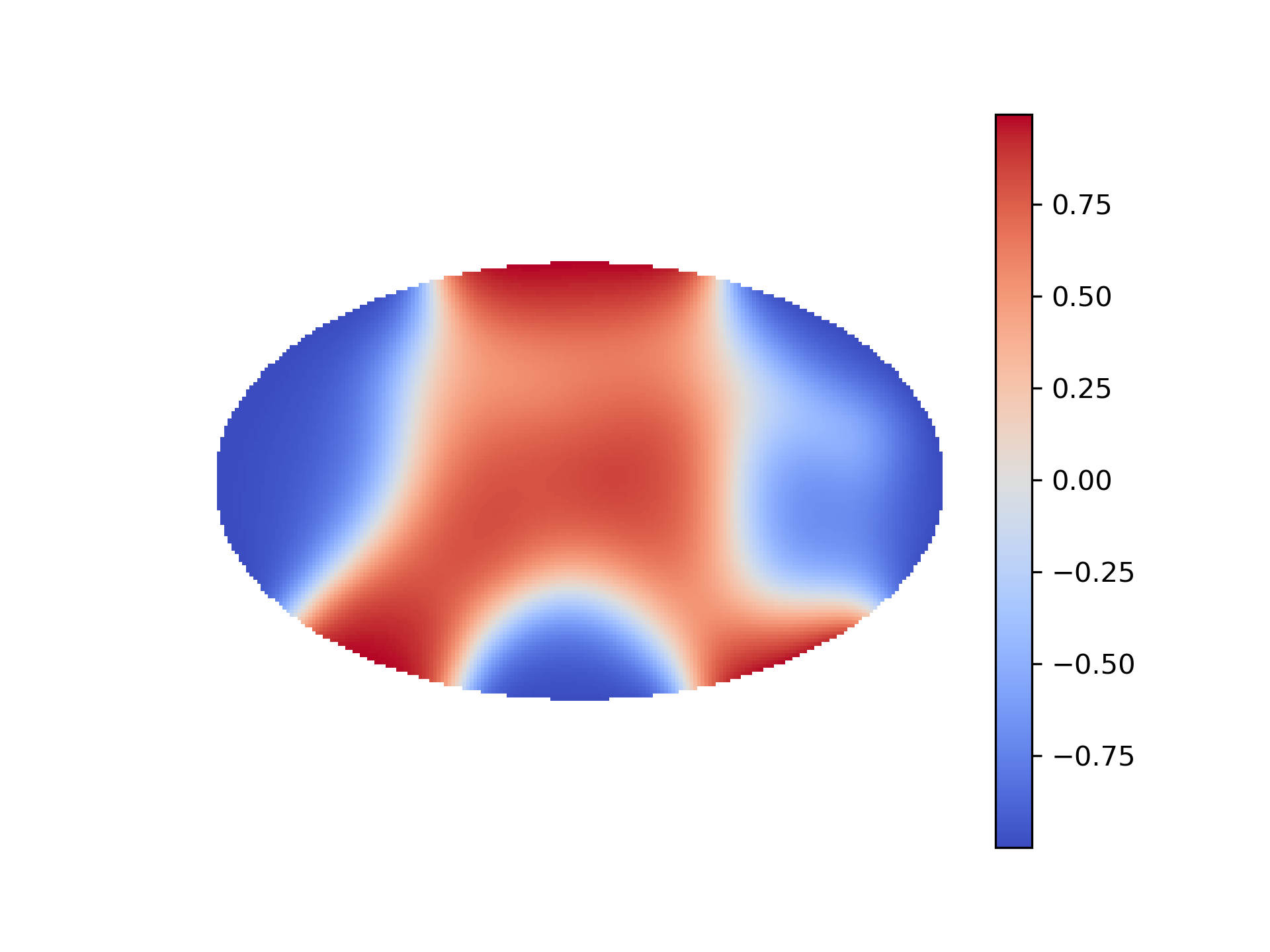}
        \includegraphics[scale=0.22]{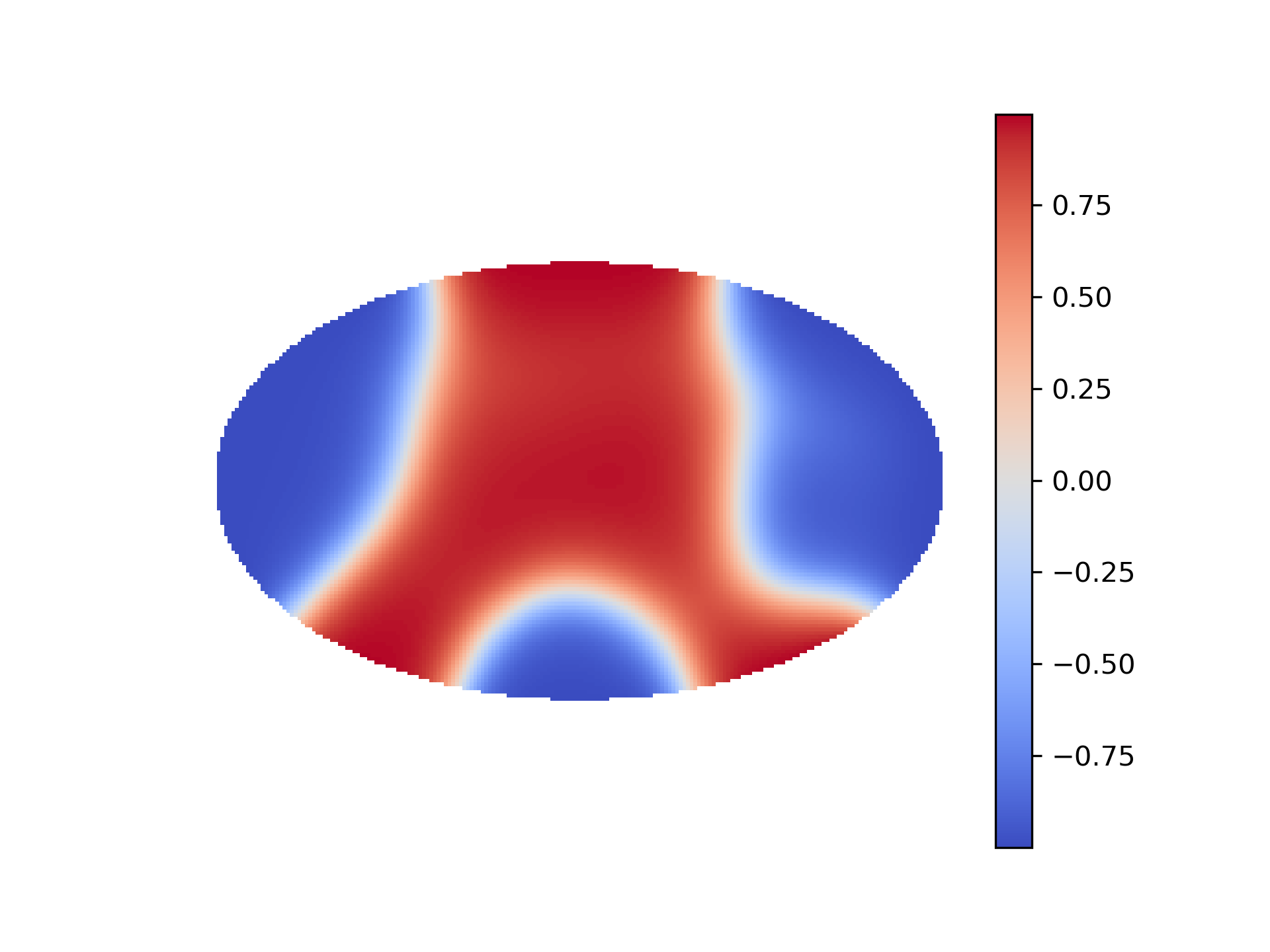}
        \caption{\( M_b=5 \) and \( M_s=10 \).}
    \end{subfigure}

    \begin{subfigure}[t]{1\textwidth}
        \centering
        \includegraphics[scale=0.22]{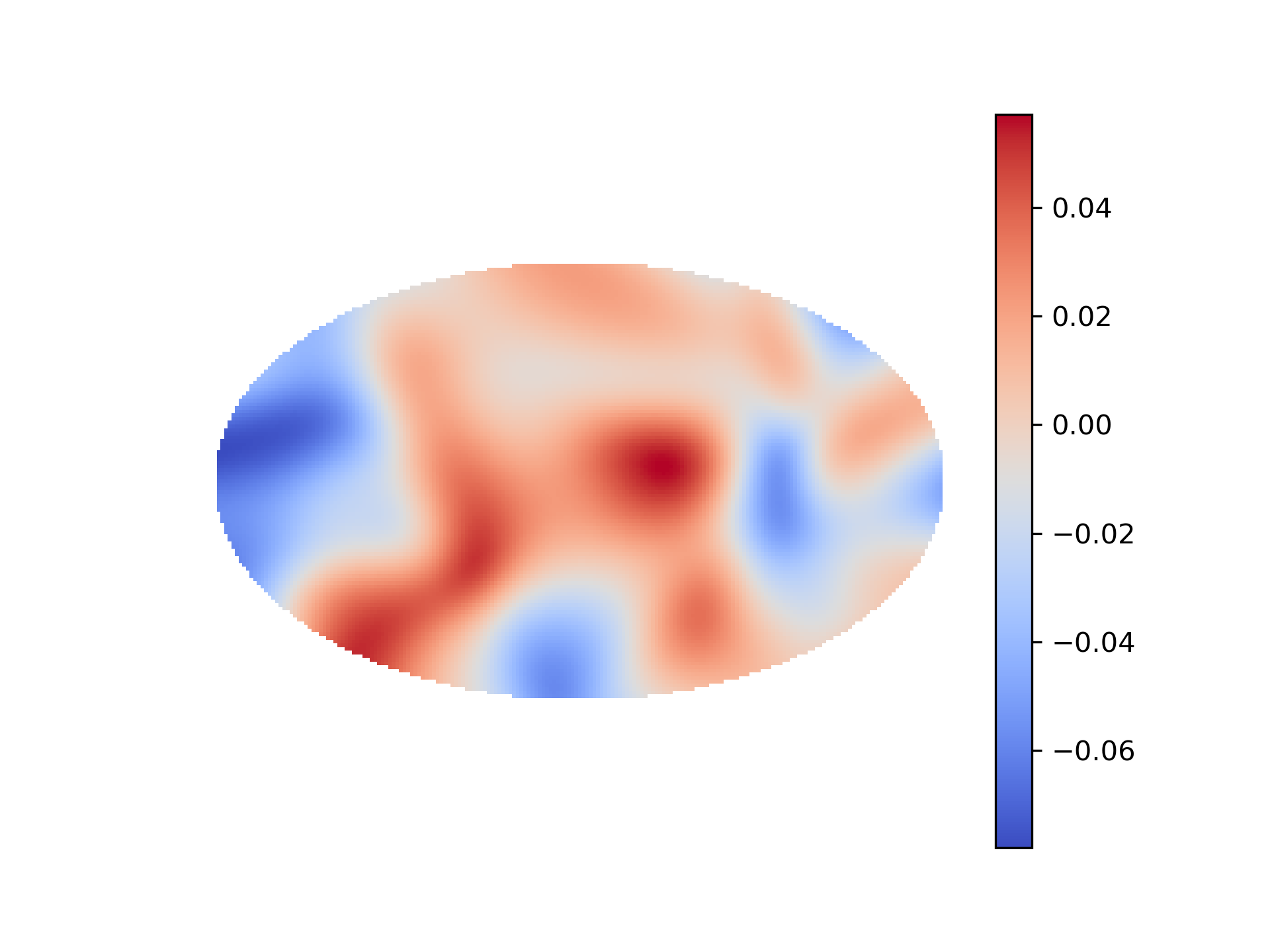}
        \includegraphics[scale=0.22]{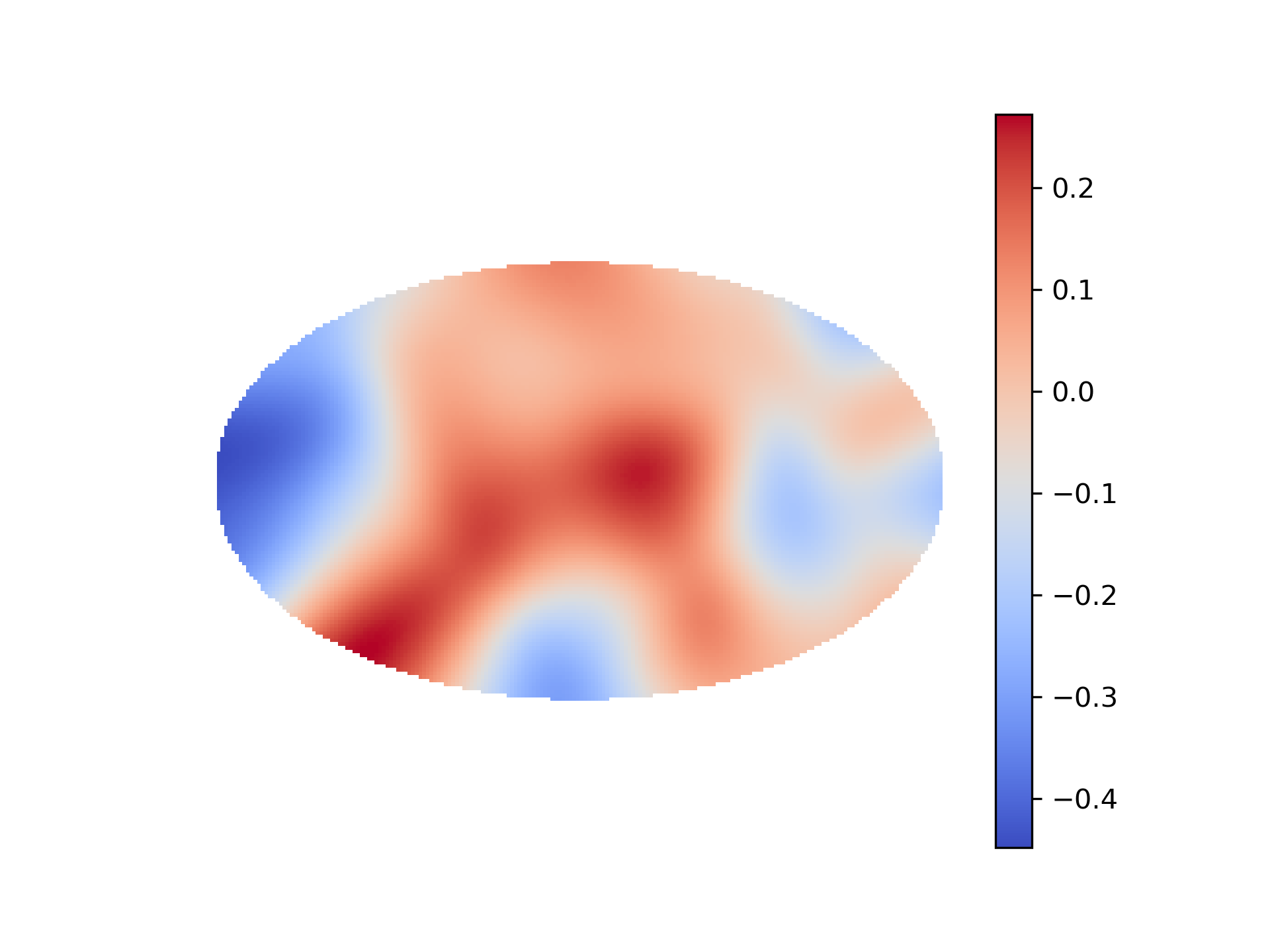}
        \includegraphics[scale=0.22]{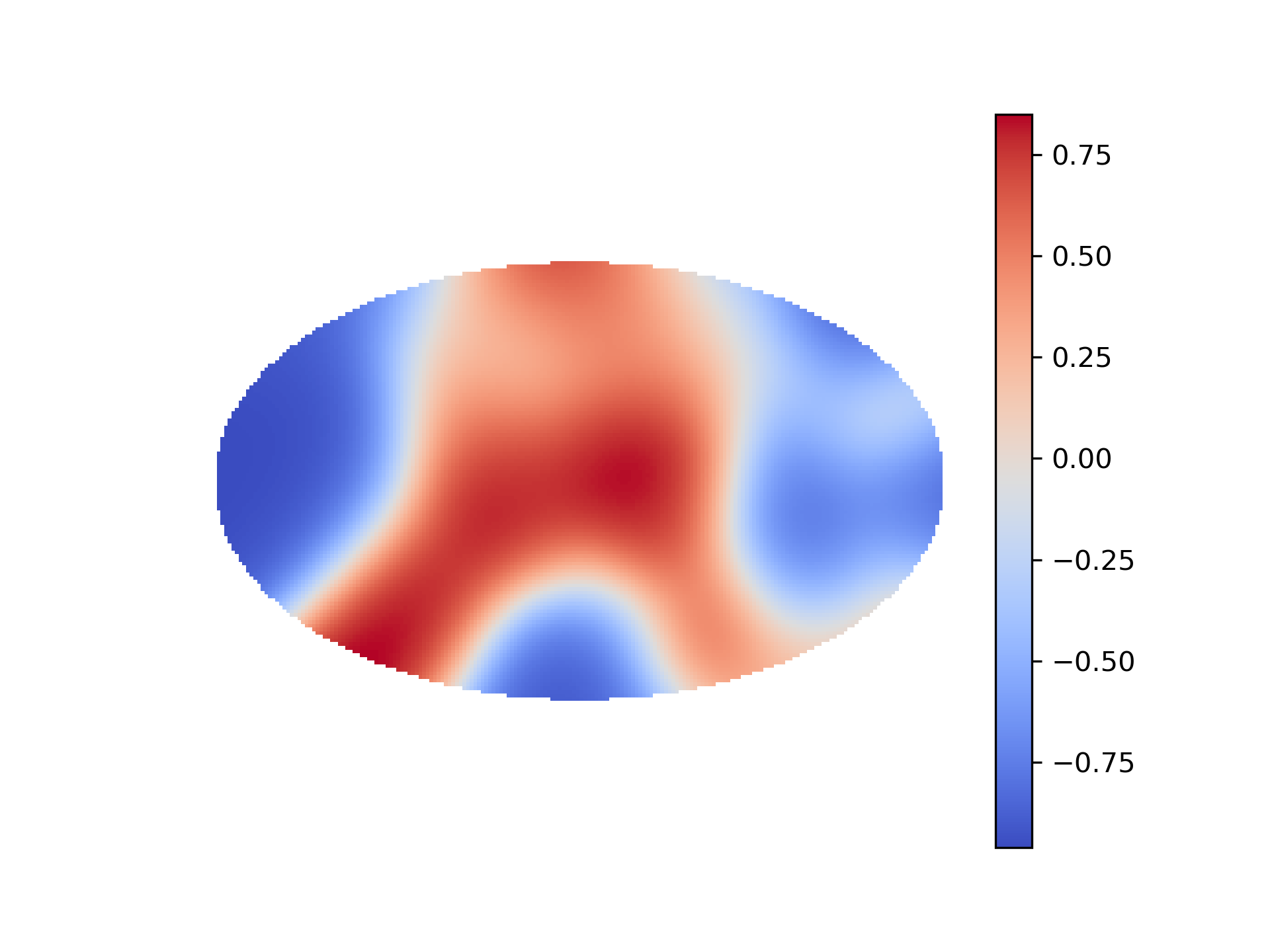}
        \includegraphics[scale=0.22]{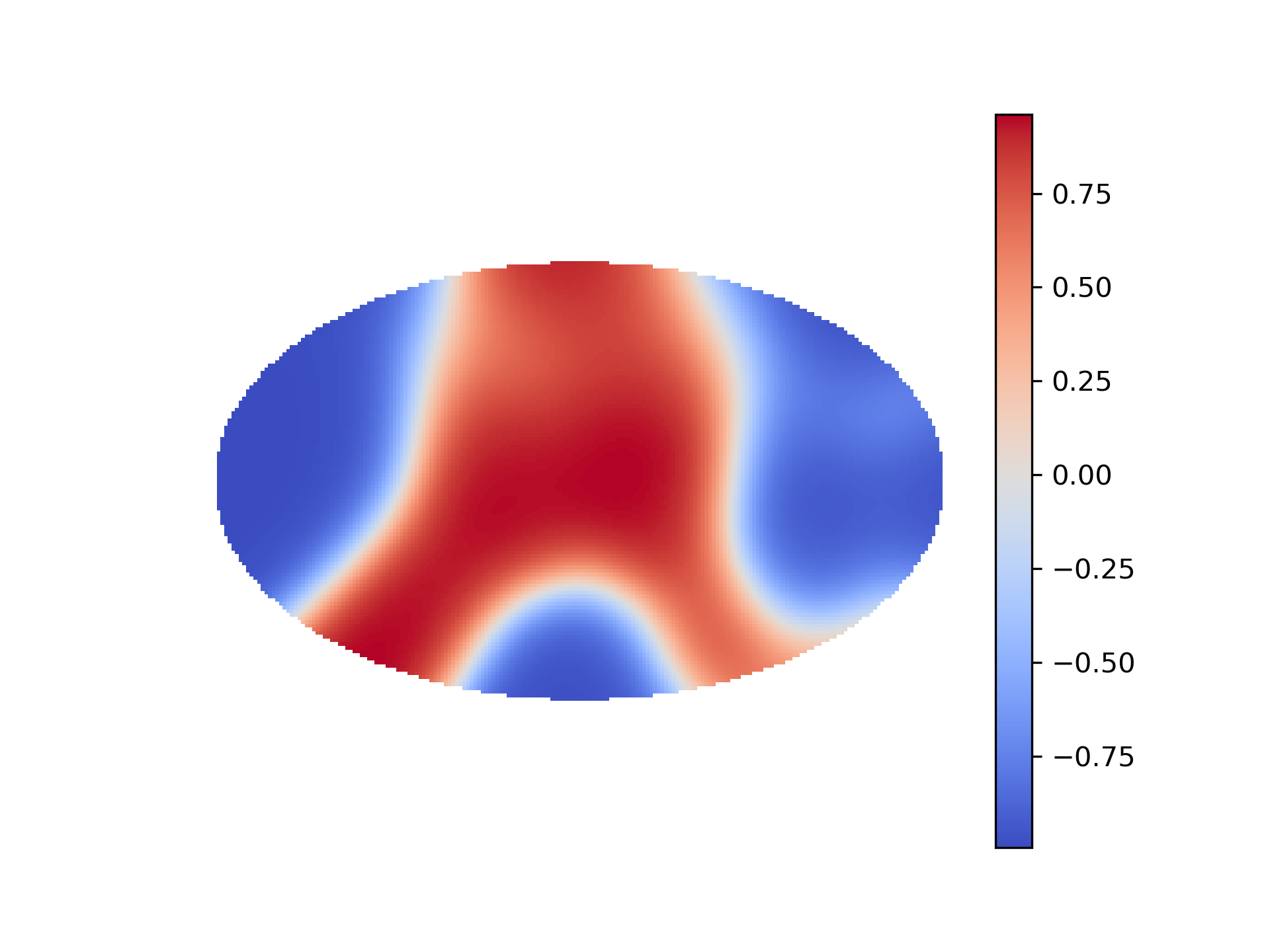}
        \caption{\( M_b=5 \) and \( M_s=5 \).}
    \end{subfigure}
   \caption{Snapshots of solutions at selected time at $t=0, 0.4, 0.8, 1.0$ with various mobility values in the case of dynamic boundary conditions in an ellipse domain. The results show that the impact of dynamic boundary conditions is noticeable near the boundary while the region of influence expands as the surface mobility enhances. }
    \label{fig:dy-dy-ellipse}
\end{figure}

We further investigate the impact of dynamic boundary conditions on bulk dynamics by examining the effect of varying surface mobility $M_s$ in model \eqref{eq:ac_2d} across two bulk mobility regimes. Our analysis encompasses both circular and elliptical domains. The results are presented in Figure \ref{fig:dy-dy} and Figure \ref{fig:dy-dy-ellipse}, respectively.

For the lower bulk mobility regime ($M_b=2$), we contrast scenarios with a high surface mobility ($M_s=10$) and matched surface-bulk mobility ($M_s=2$), respectively. Similarly, in the higher bulk mobility regime ($M_b=5$), we compare cases with an elevated surface mobility ($M_s=10$) and a matched mobility ($M_s=5$), respectively. In both cases, systems with a higher surface mobility ($M_s=10$) demonstrate accelerated boundary coarsening, which in turn drives faster the bulk dynamics near the boundary compared to their lower surface mobility counterparts, leading to seemingly weakened phase separation in the bulk.

\begin{figure}[H]
    \centering
    \includegraphics[scale=0.3]{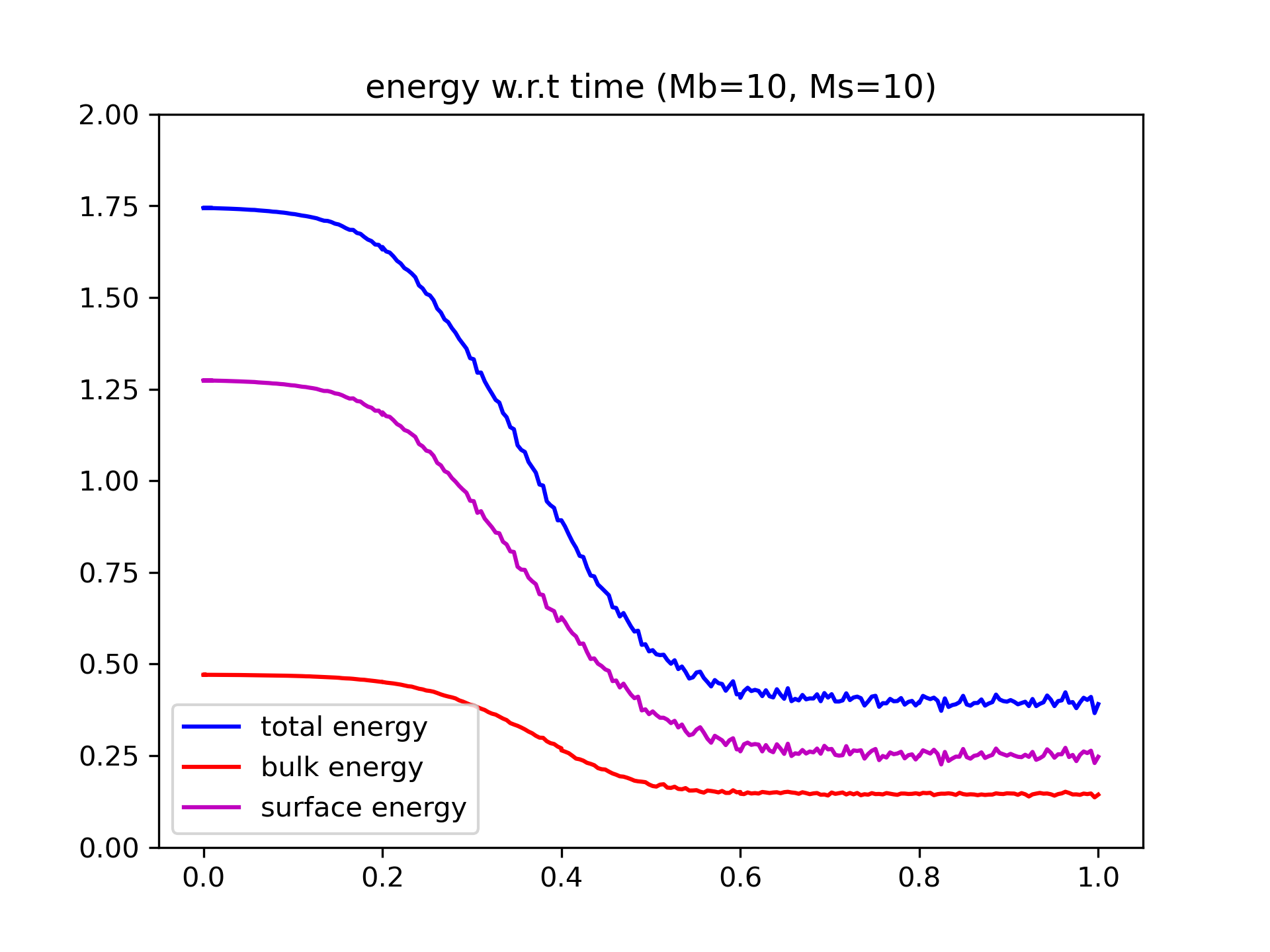}
    \includegraphics[scale=0.3]{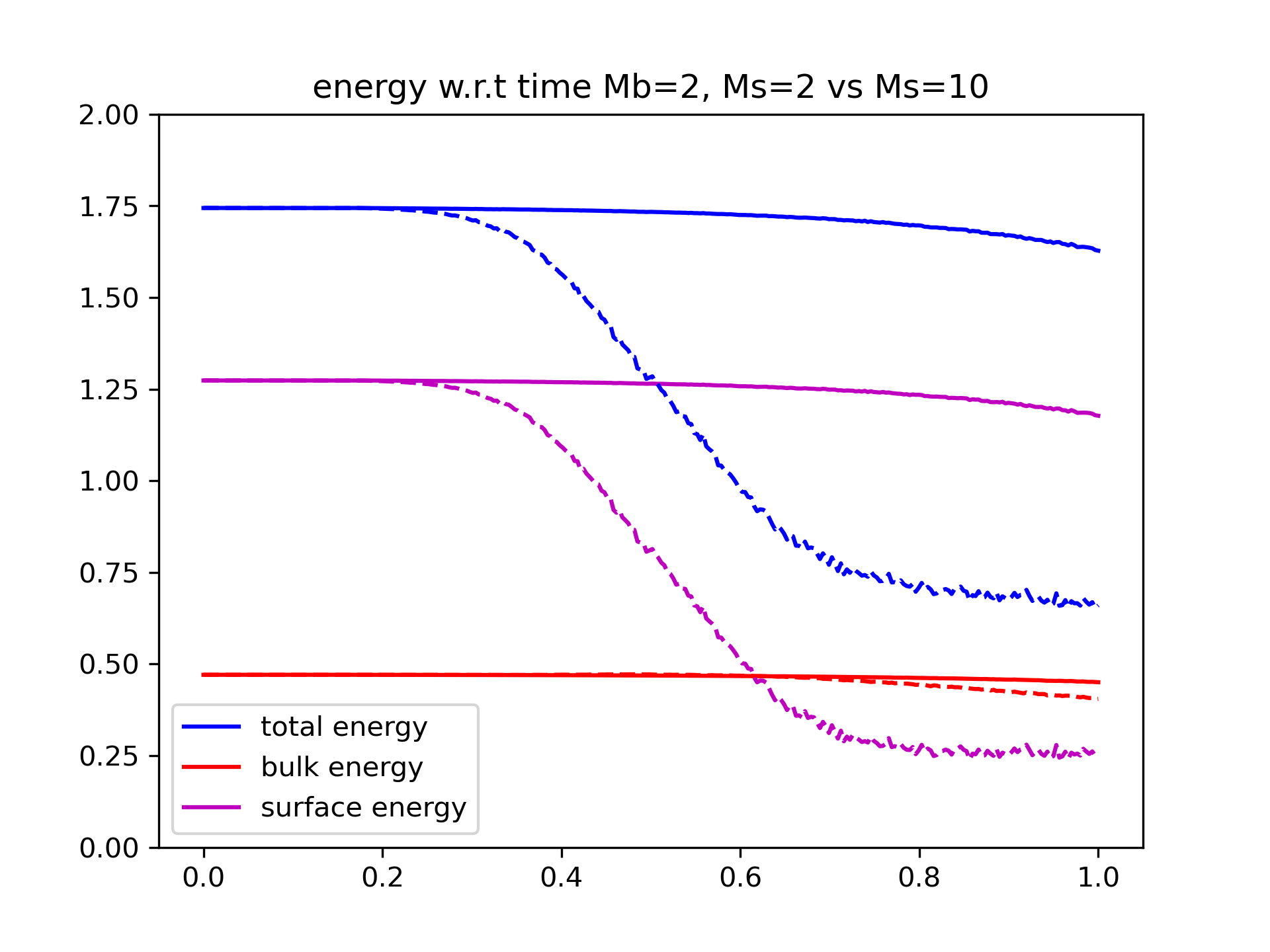}
    \includegraphics[scale=0.3]{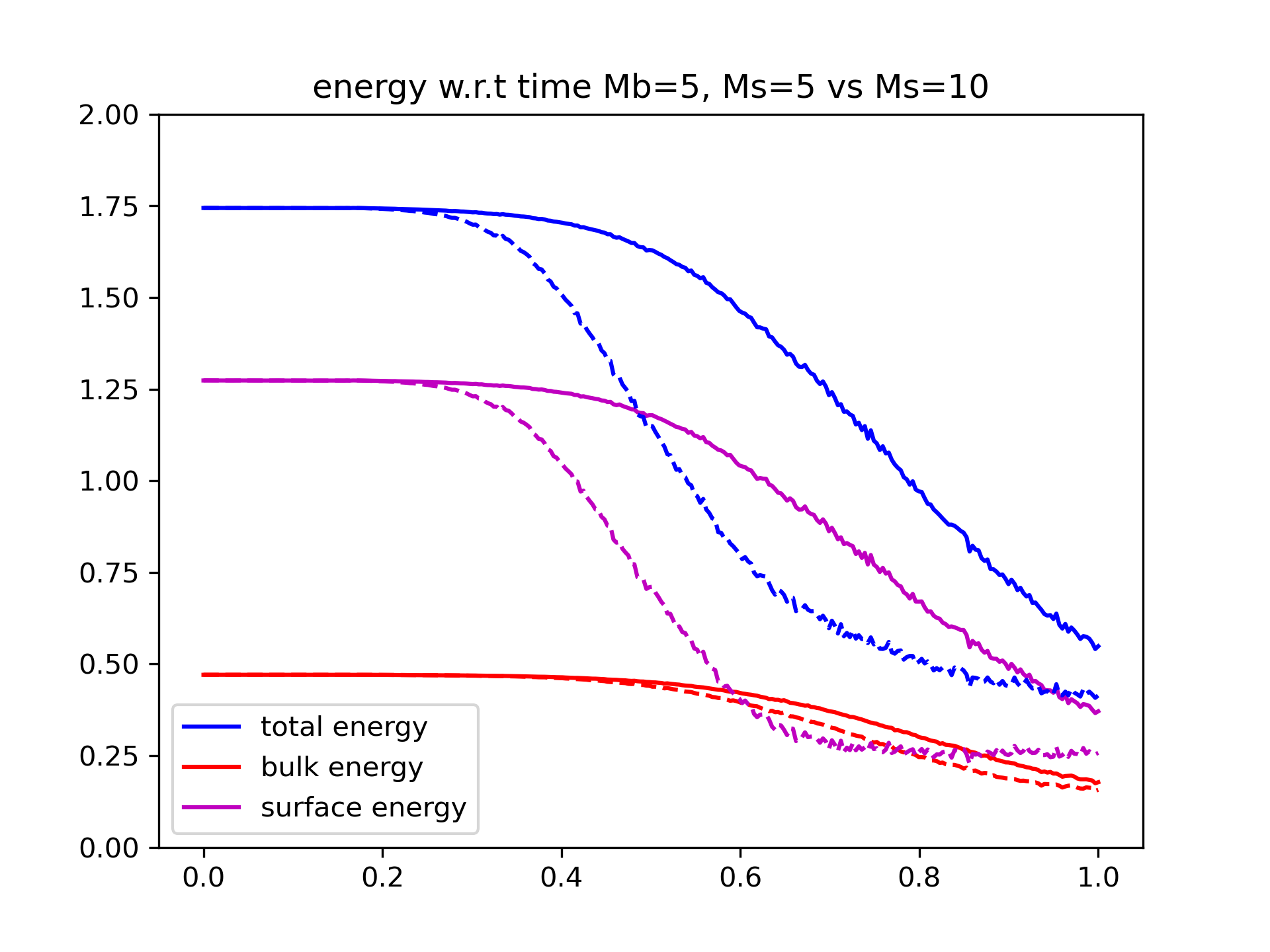}
    \caption{Time evolves of total energy, bulk energy, and surface energy for Allen-Cahn equation with dynamic boundary condition (solid line) and time evolves of total energy (bulk energy) for Allen-Cahn equation with Neumann boundary condition (dash line) in ellipse-shaped 2D computational domain. The plot of the total energy (blue), bulk energy (red), and surface energy (purple) with respect to time in an ellipse for a few selected mobility values. Left:  $M_b=M_s=10$. Middle: $M_b=M_s=2$ (solid lines) v.s. $M_b=2, M_s=10$ (dash lines). Right: $M_b=M_s=5$ (solid lines) v.s. $M_b=5, M_s=10$ (dash lines). The surface energy decay in fact  accelerates the bulk energy decay slightly.}
    \label{fig:energy-t-ellipse}
\end{figure}

The energy evolution plots (Figure \ref{fig:energy-t}, middle and right) quantitatively confirm these observations. While all systems exhibit energy dissipation, those with a higher surface mobility ($M_s=10$) show slightly faster bulk energy decays compared to systems with a lower surface mobility ($M_s=2$ or $M_s=5$). These patterns persist across both circular and elliptical domains, as evidenced by results in Figure \ref{fig:energy-t-ellipse} as well.

These findings carry significant implications for modeling systems where surface and bulk scales are comparable or where surface effects dominate. The marked influence of dynamic boundary conditions on bulk dynamics underscores the importance of selecting appropriate boundary conditions for accurate system characterization. Furthermore, our results reveal how the surface mobility can serve as a crucial control parameter for bulk dynamics, offering new insights for the system design and optimization.

The comparison between the impact of dynamic boundary conditions on bulk dynamics in elliptical domains is presented in Figure \ref{fig:dy-dy-ellipse}. The different impacts can be investigated by comparing the system's dynamics with different surface mobility values under fixed bulk mobility. The results are analogous to the case that  we have discussed in a circular domain.

\section{Conclusion}\label{sec:conclusion}

In this paper, we introduce the Energy Dissipation Rate-Guided Adaptive Sampling (EDRAS) method, a novel approach that noticeably enhances the capabilities of Physics-Informed Neural Networks (PINNs) in solving thermodynamically consistent partial differential equations. EDRAS utilizes energy dissipation rate density as a strategic guide for dynamically sampling collocation points, consistently outperforming traditional Residual-based Adaptive Refinement (RAR) methods across multiple performance metrics. As a comparison, we identify the drawbacks of the RAR method and the advantage of the EDRAS method from the probability perspective and validated the point numerically. The EDRAS method's effectiveness stems from its ability to precisely identify critical sampling points aimed at reducing the error while maintaining computational efficiency through the evaluation of lower-order derivatives  and avoiding estimate the global distribution—a particular advantage for high-dimensional PDEs where the computational cost of evaluating higher-order derivatives can be prohibitive.

We demonstrate EDRAS's practical value through a detailed investigation of the Allen-Cahn model with dynamic boundary conditions in arbitrary domains. The mesh-free nature of EDRAS-enhanced PINNs proves particularly advantageous for problems involving complex geometries, enabling us to reveal fundamental differences between static and dynamic boundary conditions' influence on bulk dynamics governed by the Allen-Cahn equation. Our results illuminate the crucial role of surface phenomena in determining system behavior, with broad implications for computational science and engineering applications where surface-bulk interactions dominate.

This work makes two significant contributions: 1) it establishes EDRAS as a powerful tool for enhancing PINNs in thermodynamically consistent systems, and 2) it provides new insights into the complex interplay between boundary conditions and bulk dynamics. These advances open new possibilities for studying thermodynamic systems in complex geometries, particularly in applications where traditional numerical methods face significant challenges.

\section*{CRediT authorship contribution statement}

CL designed and implemented EDRAS-enhanced PINN method, conducted the related simulations and wrote the draft.

WY conducted the simulations of the Allen-Cahn equation with boundary conditions on arbitrary domains using finite difference methods that he developed recently and drafted part of the paper.

QW designed the project, carried out the analysis, drafted and edited the paper.

\section*{Declaration of interest}
None

\section*{Acknowledgements}
Chunyan Li's research is partially supported by a SPARC Graduate Research Grant from the Office of the Vice President for Research at the University of South Carolina. Qi Wang's research is partially supported by
NSF awards DMS-2038080 and OIA-2242812. The funders played no role in study design, data collection, analysis and interpretation of data, or the writing of this manuscript.

\section*{Data availability}
No data was used for the research described in the article.

\bibliographystyle{plain}

\begin{thebibliography}{10}

\bibitem{chen2025learn}
Chuqi Chen, Yahong Yang, Yang Xiang, and Wenrui Hao.
\newblock Learn sharp interface solution by homotopy dynamics.
\newblock {\em arXiv preprint arXiv:2502.00488}, 2025.

\bibitem{chill2004convergence}
Ralph Chill, Eva Fa{\v{s}}angov{\'a}, and Jan Pr{\"u}ss.
\newblock {\em Convergence to Steady States of Solutions of the Cahn-Hilliard
  Equation with Dynamic Boundary Conditions}.
\newblock Univ., Fachbereich Mathematik und Informatik, 2004.

\bibitem{chunyan2023thesis}
L.~Chunyan.
\newblock {\em Deep Learning for Studying Materials Stability and Solving
  Thermodynamically Consistent PDEs with Dynamic Boundary Conditions in
  Arbitrary Domains}.
\newblock Doctoral dissertation, 2023.

\bibitem{daw2022rethinking}
Arka Daw, Jie Bu, Sifan Wang, Paris Perdikaris, and Anuj Karpatne.
\newblock Rethinking the importance of sampling in physics-informed neural
  networks.
\newblock {\em arXiv preprint arXiv:2207.02338}, 2022.

\bibitem{espath2023continuum}
Luis Espath.
\newblock A continuum framework for phase field with bulk-surface dynamics.
\newblock {\em Partial Differential Equations and Applications}, 4(1):1, 2023.

\bibitem{faure2009generalized}
Henri Faure and Christiane Lemieux.
\newblock Generalized halton sequences in 2008: A comparative study.
\newblock {\em ACM Transactions on Modeling and Computer Simulation (TOMACS)},
  19(4):1--31, 2009.

\bibitem{gal2006cahn}
Ciprian~G Gal.
\newblock A cahn--hilliard model in bounded domains with permeable walls.
\newblock {\em Mathematical methods in the applied sciences},
  29(17):2009--2036, 2006.

\bibitem{gal2008non}
Ciprian~G Gal, Maurizio Grasselli, et~al.
\newblock The non-isothermal allen-cahn equation with dynamic boundary
  conditions.
\newblock {\em Discrete Contin. Dyn. Syst}, 22(4):1009--1040, 2008.

\bibitem{gilardi2009cahn}
Gianni Gilardi, Alain Miranville, Giulio Schimperna, et~al.
\newblock On the cahn--hilliard equation with irregular potentials and dynamic
  boundary conditions.
\newblock {\em Commun. Pure Appl. Anal}, 8(3):881--912, 2009.

\bibitem{guo2024tcas}
Jia Guo, Haifeng Wang, Shilin Gu, and Chenping Hou.
\newblock Tcas-pinn: Physics-informed neural networks with a novel temporal
  causality-based adaptive sampling method.
\newblock {\em Chinese Physics B}, 33(5):050701, 2024.

\bibitem{halton1960efficiency}
John~H Halton.
\newblock On the efficiency of certain quasi-random sequences of points in
  evaluating multi-dimensional integrals.
\newblock {\em Numerische Mathematik}, 2:84--90, 1960.

\bibitem{hammersley2013monte}
John Hammersley.
\newblock {\em Monte carlo methods}.
\newblock Springer Science \& Business Media, 2013.

\bibitem{hammersley1960monte}
John~M Hammersley.
\newblock Monte carlo methods for solving multivariable problems.
\newblock {\em Annals of the New York Academy of Sciences}, 86(3):844--874,
  1960.

\bibitem{han2022residual}
Jiayue Han, Zhiqiang Cai, Zhiyou Wu, and Xiang Zhou.
\newblock Residual-quantile adjustment for adaptive training of
  physics-informed neural network.
\newblock In {\em 2022 IEEE International Conference on Big Data (Big Data)},
  pages 921--930. IEEE, 2022.

\bibitem{hao2024multiscale}
Wenrui Hao, Rui~Peng Li, Yuanzhe Xi, Tianshi Xu, and Yahong Yang.
\newblock Multiscale neural networks for approximating green's functions.
\newblock {\em arXiv preprint arXiv:2410.18439}, 2024.

\bibitem{xiaoboentropy}
Xiaobo Jing and Qi~Wang.
\newblock Thermodynamically consistent models for coupled bulk and surface
  dynamics.
\newblock {\em Entropy}, 24(11):1683, 2022.

\bibitem{xiaobocms}
Xiaobo Jing and Qi~Wang.
\newblock Thermodynamically consistent dynamic boundary conditions of phase
  field models.
\newblock {\em Communications in Mathematical Sciences}, 21(3), 2023.

\bibitem{joe2008constructing}
Stephen Joe and Frances~Y Kuo.
\newblock Constructing sobol sequences with better two-dimensional projections.
\newblock {\em SIAM Journal on Scientific Computing}, 30(5):2635--2654, 2008.

\bibitem{kingma2014adam}
Diederik~P Kingma.
\newblock Adam: A method for stochastic optimization.
\newblock {\em arXiv preprint arXiv:1412.6980}, 2014.

\bibitem{krishnapriyan2021characterizing}
Aditi Krishnapriyan, Amir Gholami, Shandian Zhe, Robert Kirby, and Michael~W
  Mahoney.
\newblock Characterizing possible failure modes in physics-informed neural
  networks.
\newblock {\em Advances in neural information processing systems},
  34:26548--26560, 2021.

\bibitem{Jia2020adapinn}
Colby L.~Wight and Jia Zhao.
\newblock Solving allen-cahn and cahn-hilliard equations using the adaptive
  physics informed neural networks.
\newblock {\em Communications in Computational Physics}, 29(3):930--954, 2021.

\bibitem{liu2021frontiers}
Xiang~Yang Liu.
\newblock {\em Frontiers and progress of current soft matter research}.
\newblock Springer, 2021.

\bibitem{lu2021deepxde}
Lu~Lu, Xuhui Meng, Zhiping Mao, and George~Em Karniadakis.
\newblock Deepxde: A deep learning library for solving differential equations.
\newblock {\em SIAM review}, 63(1):208--228, 2021.

\bibitem{Yangyahong}
Jiaqi Luo, Yahong Yang, Yuan Yuan, Shixin Xu, and Wenrui Hao.
\newblock An imbalanced learning-based sampling method for physics-informed
  neural networks.
\newblock {\em arXiv preprint arXiv:2501.11222}, 2025.

\bibitem{mao2023}
Zhiping Mao and Xuhui Meng.
\newblock Physics-informed neural networks with residual/gradient-based
  adaptive sampling methods for solving partial differential equations with
  sharp solutions.
\newblock {\em Applied Mathematics and Mechanics}, 44(7):1069--1084, 2023.

\bibitem{bcpinn-mattey2022novel}
Revanth Mattey and Susanta Ghosh.
\newblock A novel sequential method to train physics informed neural networks
  for allen cahn and cahn hilliard equations.
\newblock {\em Computer Methods in Applied Mechanics and Engineering},
  390:114474, 2022.

\bibitem{nabian2021efficient}
Mohammad~Amin Nabian, Rini~Jasmine Gladstone, and Hadi Meidani.
\newblock Efficient training of physics-informed neural networks via importance
  sampling.
\newblock {\em Computer-Aided Civil and Infrastructure Engineering},
  36(8):962--977, 2021.

\bibitem{wright1999numerical}
Jorge Nocedal and Stephen~J Wright.
\newblock {\em Numerical optimization}.
\newblock Springer, 1999.

\bibitem{pruss2006maximal}
Jan Pr{\"u}ss, Reinhard Racke, and Songmu Zheng.
\newblock Maximal regularity and asymptotic behavior of solutions for the
  cahn--hilliard equation with dynamic boundary conditions.
\newblock {\em Annali di Matematica Pura ed Applicata}, 185:627--648, 2006.

\bibitem{racke2003cahn}
Reinhard Racke and Songmu Zheng.
\newblock {The Cahn-Hilliard equation with dynamic boundary conditions}.
\newblock {\em Advances in Differential Equations}, 8(1):83 -- 110, 2003.

\bibitem{raissi2019physics}
Maziar Raissi, Paris Perdikaris, and George~E Karniadakis.
\newblock Physics-informed neural networks: A deep learning framework for
  solving forward and inverse problems involving nonlinear partial differential
  equations.
\newblock {\em Journal of Computational physics}, 378:686--707, 2019.

\bibitem{sobol1967distribution}
Il'ya~Meerovich Sobol'.
\newblock On the distribution of points in a cube and the approximate
  evaluation of integrals.
\newblock {\em Zhurnal Vychislitel'noi Matematiki i Matematicheskoi Fiziki},
  7(4):784--802, 1967.

\bibitem{stein1987large}
Michael Stein.
\newblock Large sample properties of simulations using latin hypercube
  sampling.
\newblock {\em Technometrics}, 29(2):143--151, 1987.

\bibitem{verfurth1994posteriori}
R{\"u}diger Verf{\"u}rth.
\newblock A posteriori error estimation and adaptive mesh-refinement
  techniques.
\newblock {\em Journal of Computational and Applied Mathematics},
  50(1-3):67--83, 1994.

\bibitem{wang2024respecting}
Sifan Wang, Shyam Sankaran, and Paris Perdikaris.
\newblock Respecting causality for training physics-informed neural networks.
\newblock {\em Computer Methods in Applied Mechanics and Engineering},
  421:116813, 2024.

\bibitem{wu2023comprehensive}
Chenxi Wu, Min Zhu, Qinyang Tan, Yadhu Kartha, and Lu~Lu.
\newblock A comprehensive study of non-adaptive and residual-based adaptive
  sampling for physics-informed neural networks.
\newblock {\em Computer Methods in Applied Mechanics and Engineering},
  403:115671, 2023.

\bibitem{wu2004convergence}
Hao Wu and Songmu Zheng.
\newblock Convergence to equilibrium for the cahn--hilliard equation with
  dynamic boundary conditions.
\newblock {\em Journal of differential equations}, 204(2):511--531, 2004.

\bibitem{xu2024}
Zhi-Qin~John Xu, Yaoyu Zhang, and Tao Luo.
\newblock Overview frequency principle/spectral bias in deep learning.
\newblock {\em Communications on Applied Mathematics and Computation}, pages
  1--38, 2024.

\bibitem{yu2024domain}
Wenkai Yu, Qi~Wang, and Tiezheng Qian.
\newblock Variational embedding domain method for thermodynamically consistent
  allen-cahn model in arbitrary domains.
\newblock {\em submitted}.

\bibitem{zienkiewicz2005finite}
Olek~C Zienkiewicz and Robert~L Taylor.
\newblock {\em The finite element method set}.
\newblock Elsevier, 2005.

\bibitem{CiCP-28-5}
Liu Ziqi, Cai Wei, and John Zhi-Qin, Xu.
\newblock Multi-scale deep neural network (mscalednn) for solving
  poisson-boltzmann equation in complex domains.
\newblock {\em Communications in Computational Physics}, 28(5):1970--2001,
  2020.

\end{thebibliography}

\section{ Appendix}
\begin{figure}[ht]
   \centering
     \begin{minipage}[b]{0.9\textwidth}
         \includegraphics[scale=0.4]{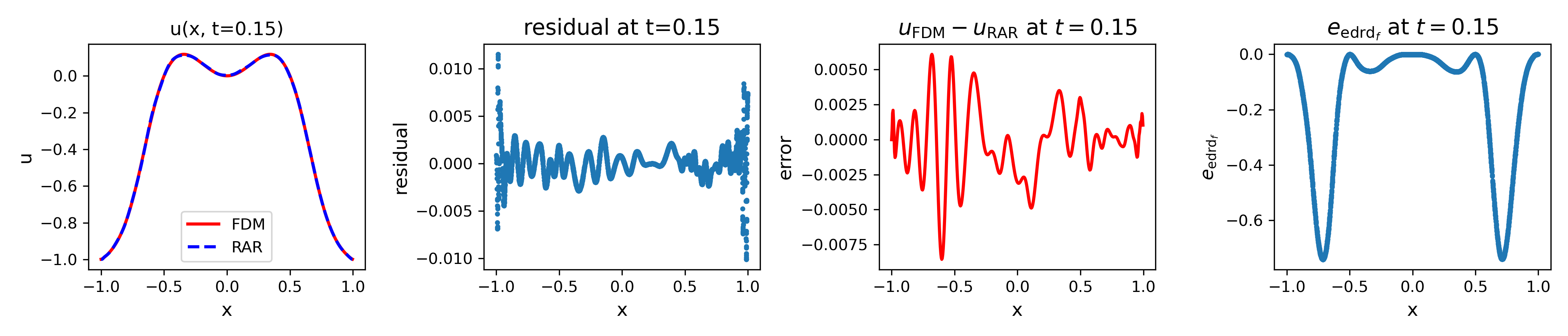}
     \end{minipage}
     \begin{minipage}[b]{0.9\textwidth}
         \includegraphics[scale=0.4]{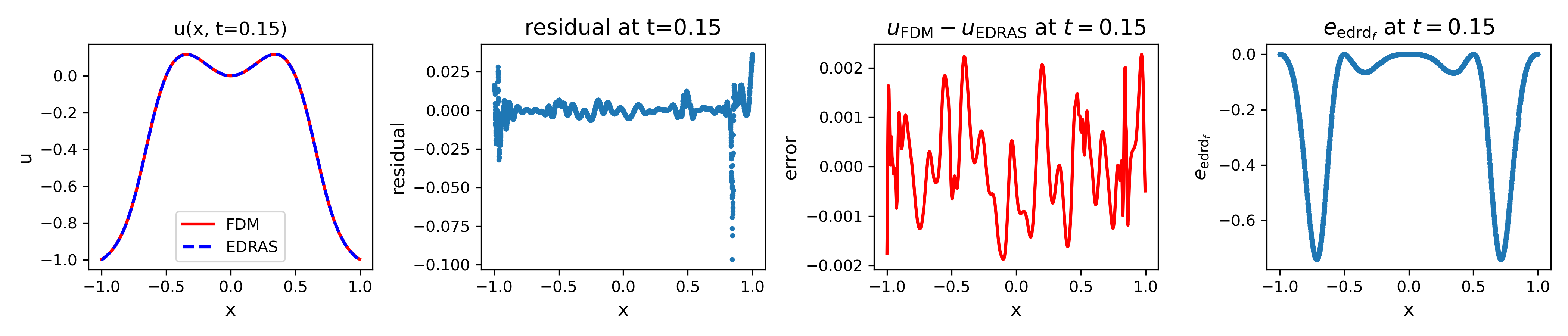}
     \end{minipage}
     	\begin{minipage}[b]{0.9\textwidth}
         \includegraphics[scale=0.4]{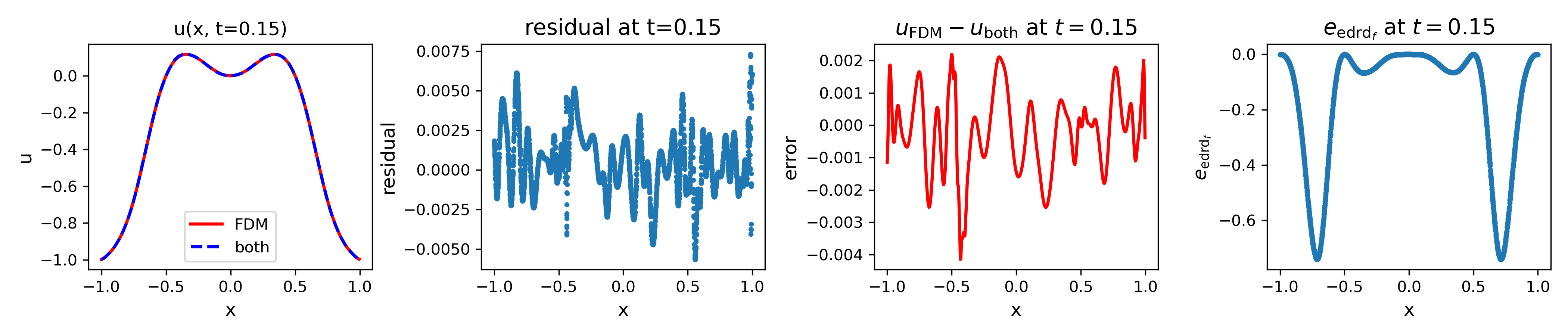}
     \end{minipage}
     \caption{Snapshots of solutions obtained using PINN with  RAR (top), EDRAS (middle) and RAR+EDRAS (bottom) at $t=0.15$.}
     \label{fig:ac1d-015}
     \end{figure}

\begin{figure}[ht]
   \centering
     \begin{minipage}[b]{0.9\textwidth}
         \includegraphics[scale=0.4]{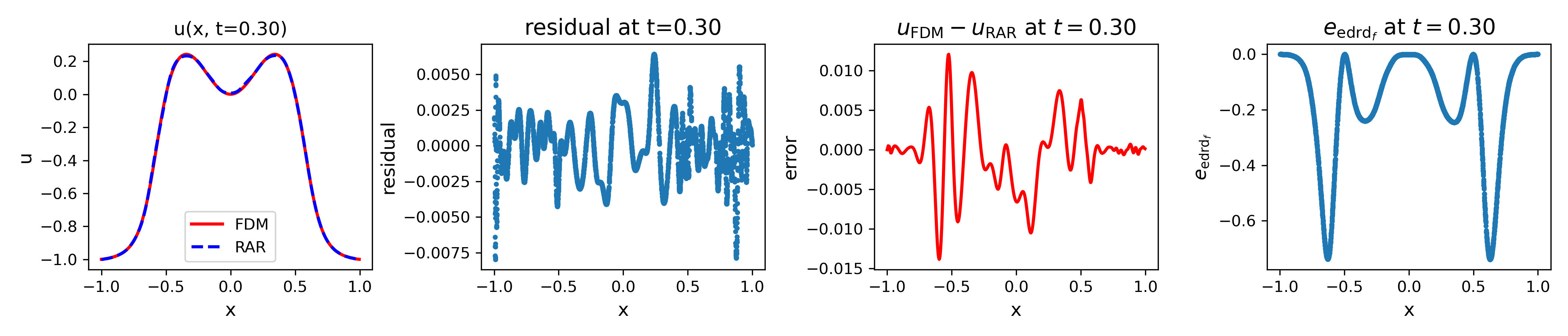}
     \end{minipage}
     \begin{minipage}[b]{0.9\textwidth}
         \includegraphics[scale=0.4]{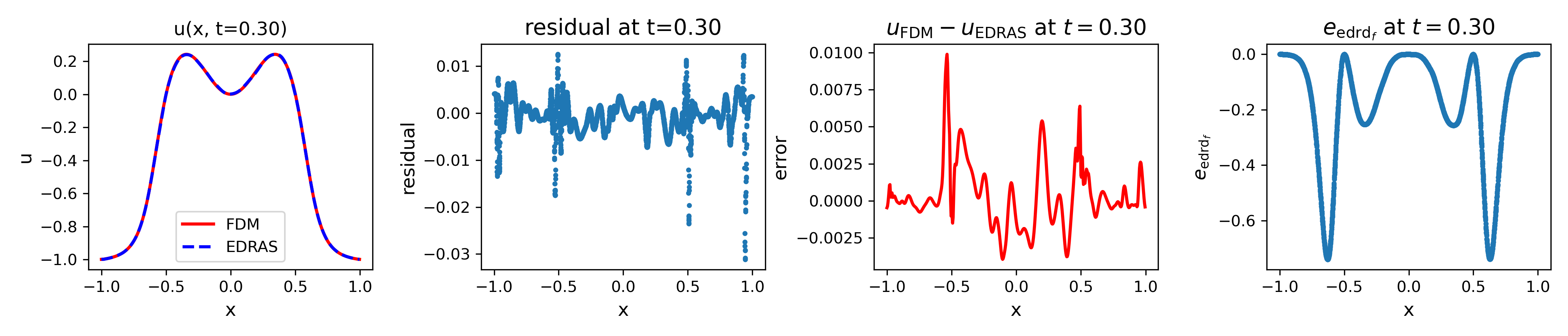}
     \end{minipage}
     \begin{minipage}[b]{0.9\textwidth}
         \includegraphics[scale=0.4]{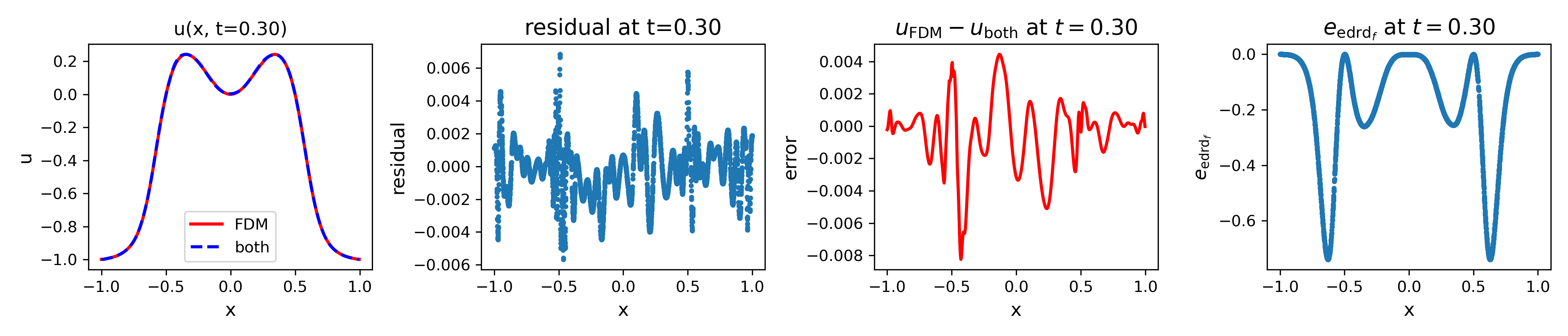}
     \end{minipage}
     \caption{Snapshots of  solutions obtained using PINN with  RAR (top), EDRAS (middle) and RAR+EDRAS (bottom) at $t=0.3$.}
     \label{fig:ac1d-03}
 \end{figure}

\begin{figure}[ht]
     \centering
       \begin{minipage}[b]{0.9\textwidth}
         \includegraphics[scale=0.4]{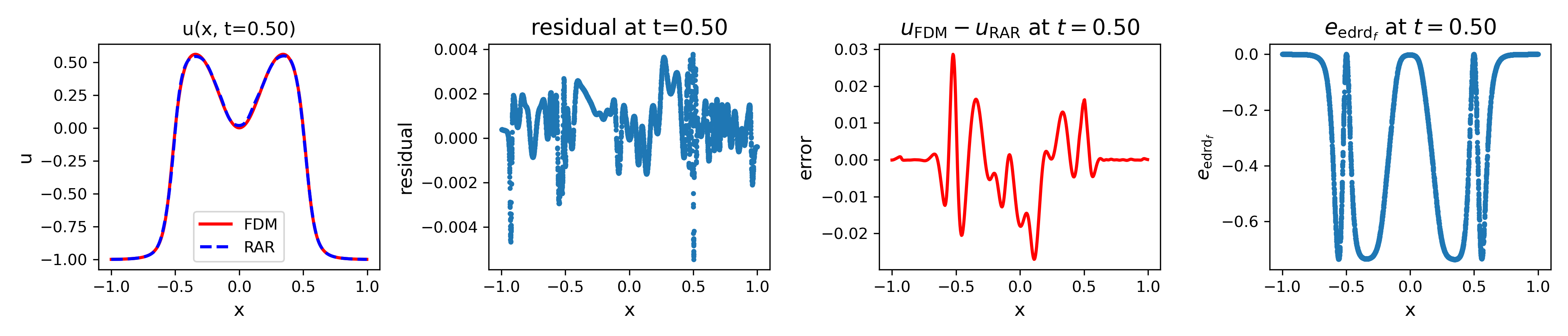}
     \end{minipage}
     \begin{minipage}[b]{0.9\textwidth}
         \includegraphics[scale=0.4]{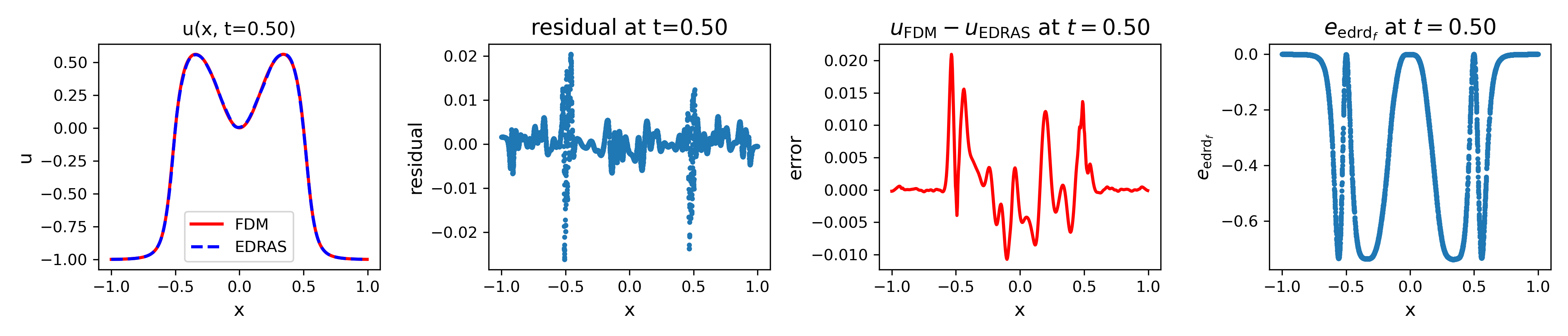}
     \end{minipage}
     \begin{minipage}[b]{0.9\textwidth}
         \includegraphics[scale=0.4]{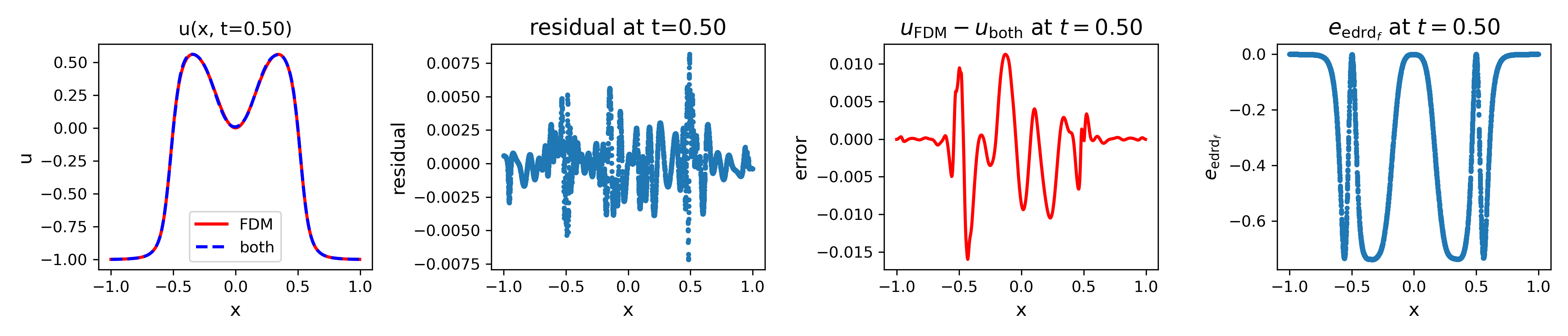}
     \end{minipage}
    \caption{Snapshots of  solutions obtained using PINN with RAR (top), EDRAS (middle) and RAR+EDRAS (bottom) at $t=0.5$.}
     \label{fig:ac1d-5}
 \end{figure}

\begin{figure}[ht]
    \centering
       \begin{minipage}[b]{0.9\textwidth}
         \includegraphics[scale=0.4]{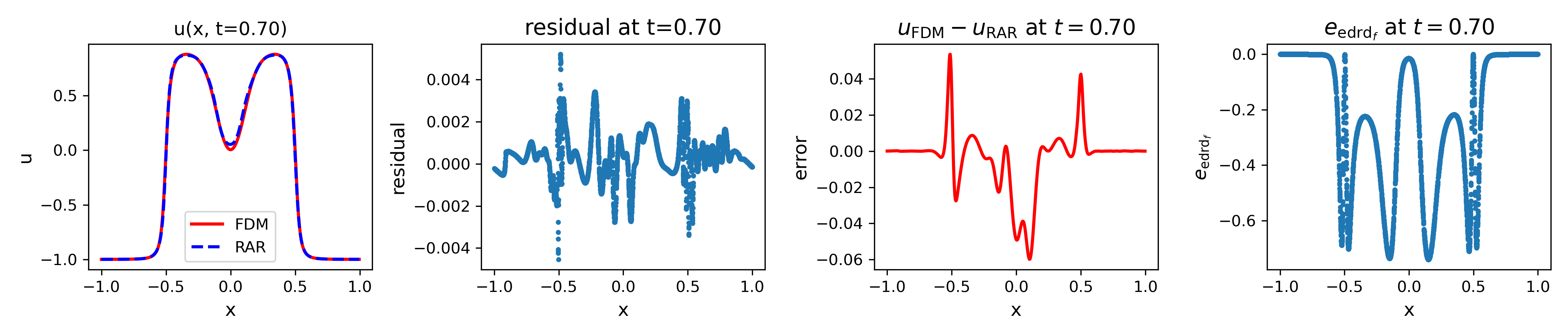}
     \end{minipage}
     \begin{minipage}[b]{0.9\textwidth}
         \includegraphics[scale=0.4]{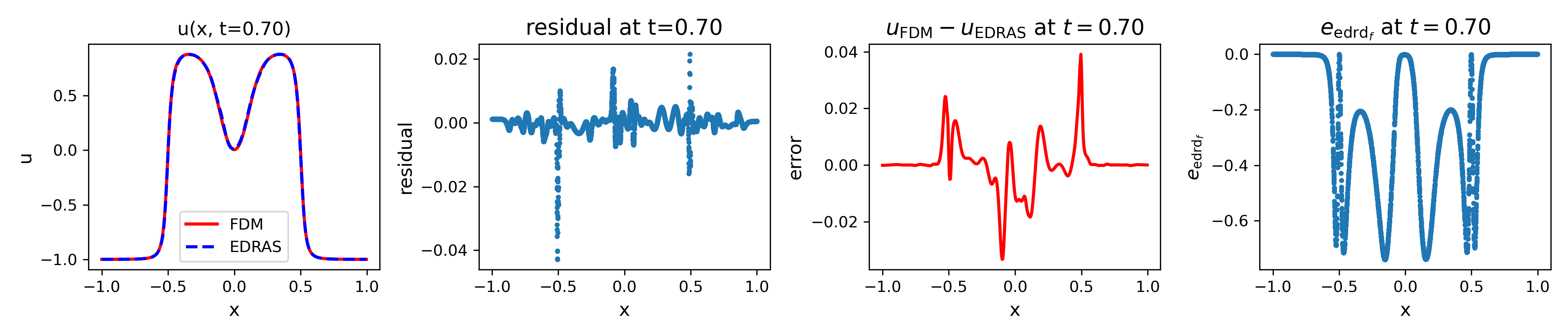}
     \end{minipage}
     \begin{minipage}[b]{0.9\textwidth}
         \includegraphics[scale=0.4]{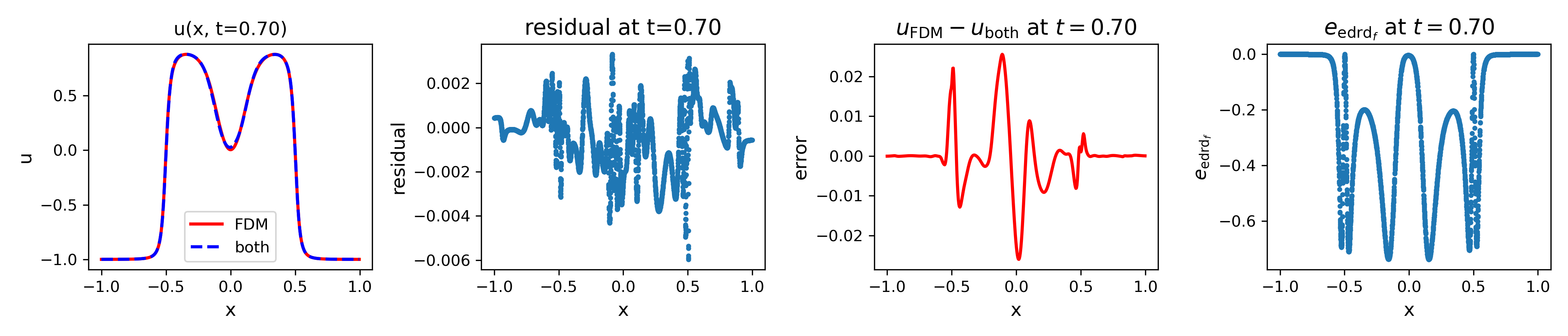}
     \end{minipage}
     \caption{Snapshots of  solutions obtained using PINN with RAR (top), EDRAS (middle) and RAR+EDRAS (bottom) at $t=0.7$.}
     \label{fig:ac1d-7}
 \end{figure}

Snapshots of the numerical solutions of the Allen-Cahn equation in 1D at $t=0.15, 0.3, 0.5$ and $t=0.7$ are presented in Figure \ref{fig:ac1d-015}, Figure \ref{fig:ac1d-03}, Figure \ref{fig:ac1d-5} and Figure  \ref{fig:ac1d-7} using EDRAS method, RAR method and the combination of both, respectively.

\begin{figure}[ht]
      \centering
         \includegraphics[scale=0.3]{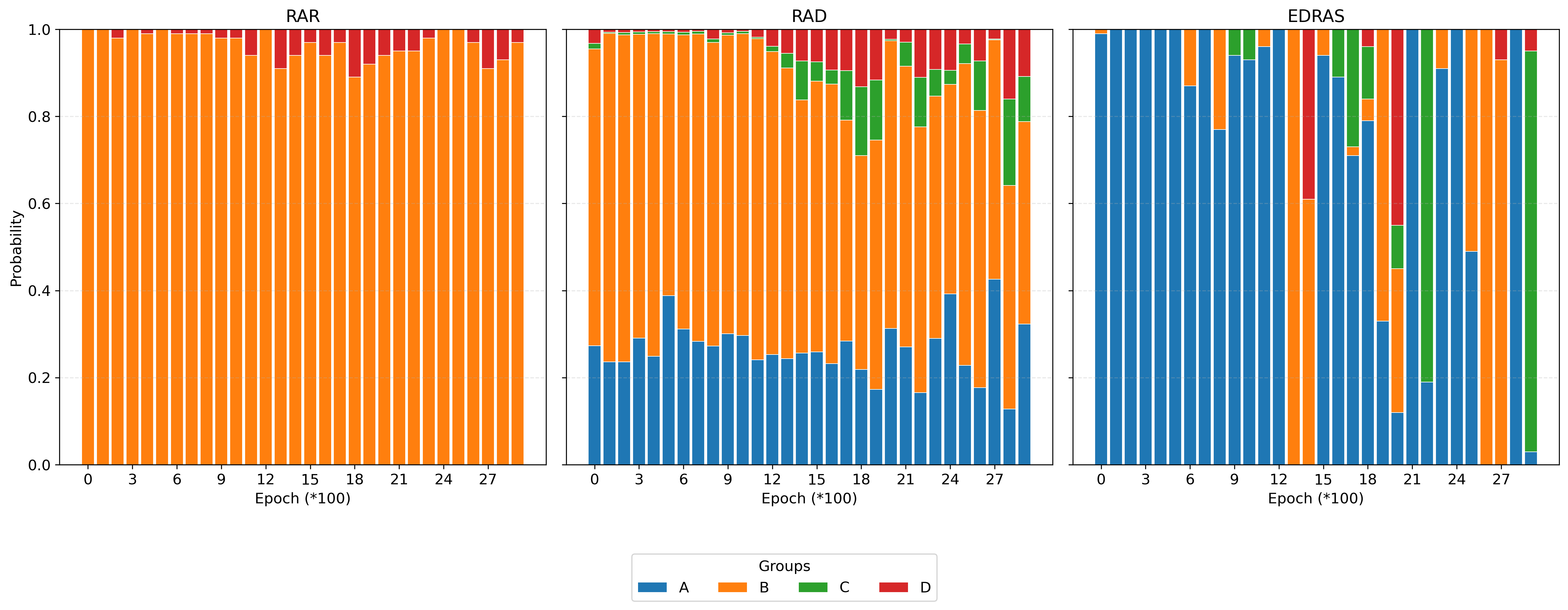}
      \caption{Probabilities of the 4 groups ($A,B,C,D$) for RAR, RAD and EDRAS in domain $[0, 0.1]\times [-1,1]$ plotted over training epochs. }
 \label{fig:prob_view01}
\end{figure}

\begin{figure}[ht]
      \centering
         \includegraphics[scale=0.3]{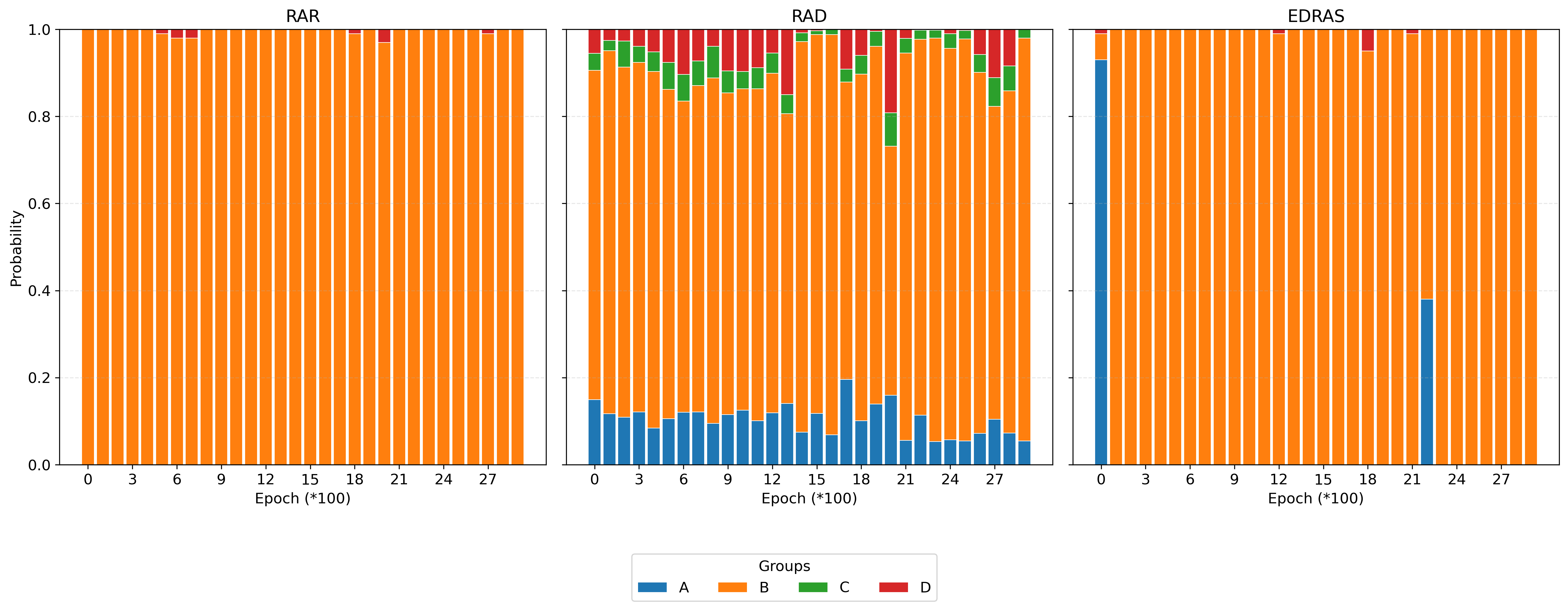}
      \caption{Probabilities of the 4 groups for RAR, RAD and EDRAS in domain $[0.1, 0.2]\times [-1,1]$ plotted over training epochs. }
 \label{fig:prob_view02}
\end{figure}

\begin{figure}[ht]
      \centering
         \includegraphics[scale=0.3]{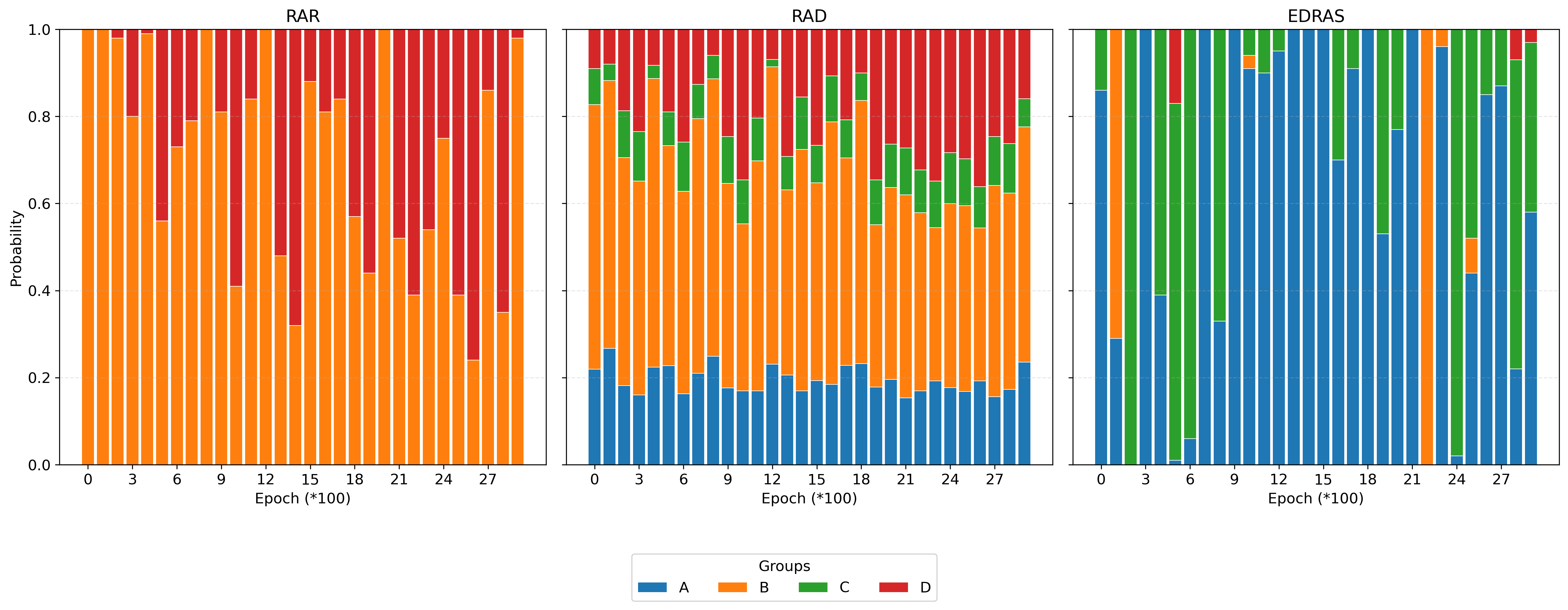}
      \caption{Probabilities of the 4 groups for RAR, RAD and EDRAS in domain $[0.2, 0.4]\times [-1,1]$ plotted over training epochs. }
 \label{fig:prob_view04}
\end{figure}

\begin{figure}[ht]
      \centering
         \includegraphics[scale=0.3]{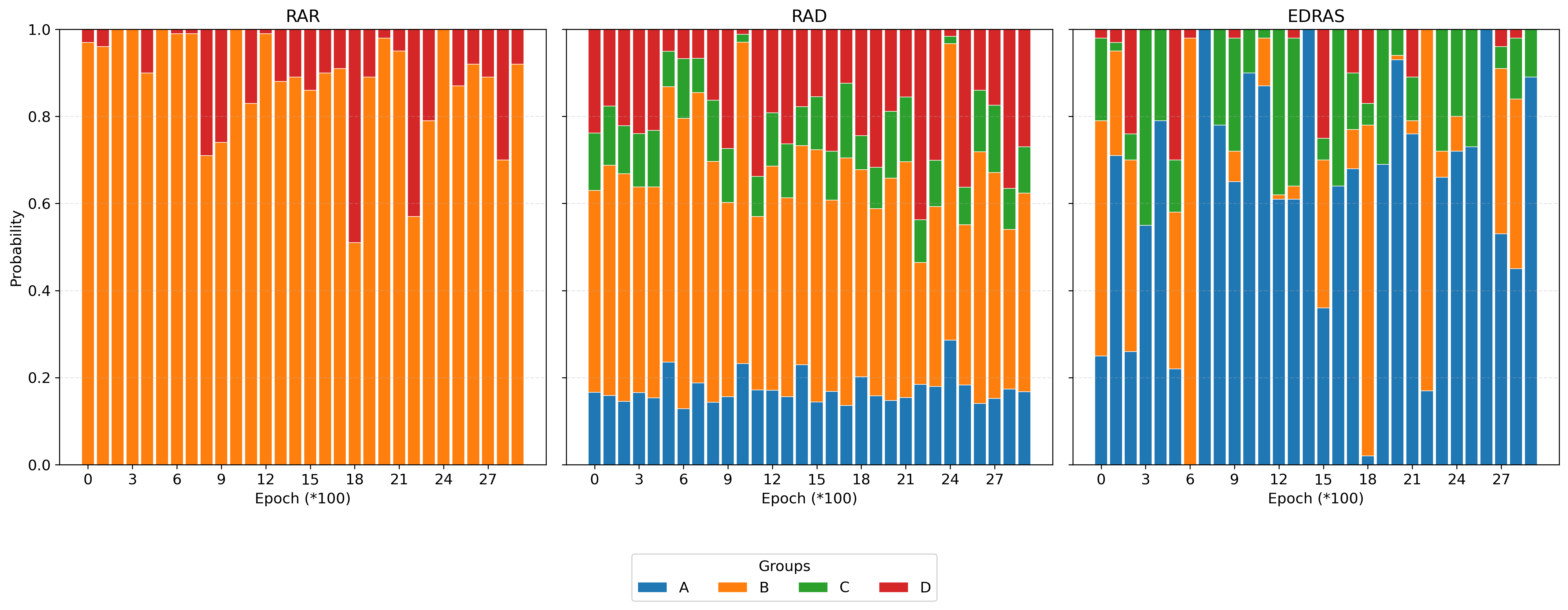}
      \caption{Probabilities of the 4 groups for RAR, RAD and EDRAS in domain $[0.4, 0.6]\times [-1,1]$ plotted over training epochs. }
 \label{fig:prob_view06}
\end{figure}

\begin{figure}[ht]
      \centering
         \includegraphics[scale=0.3]{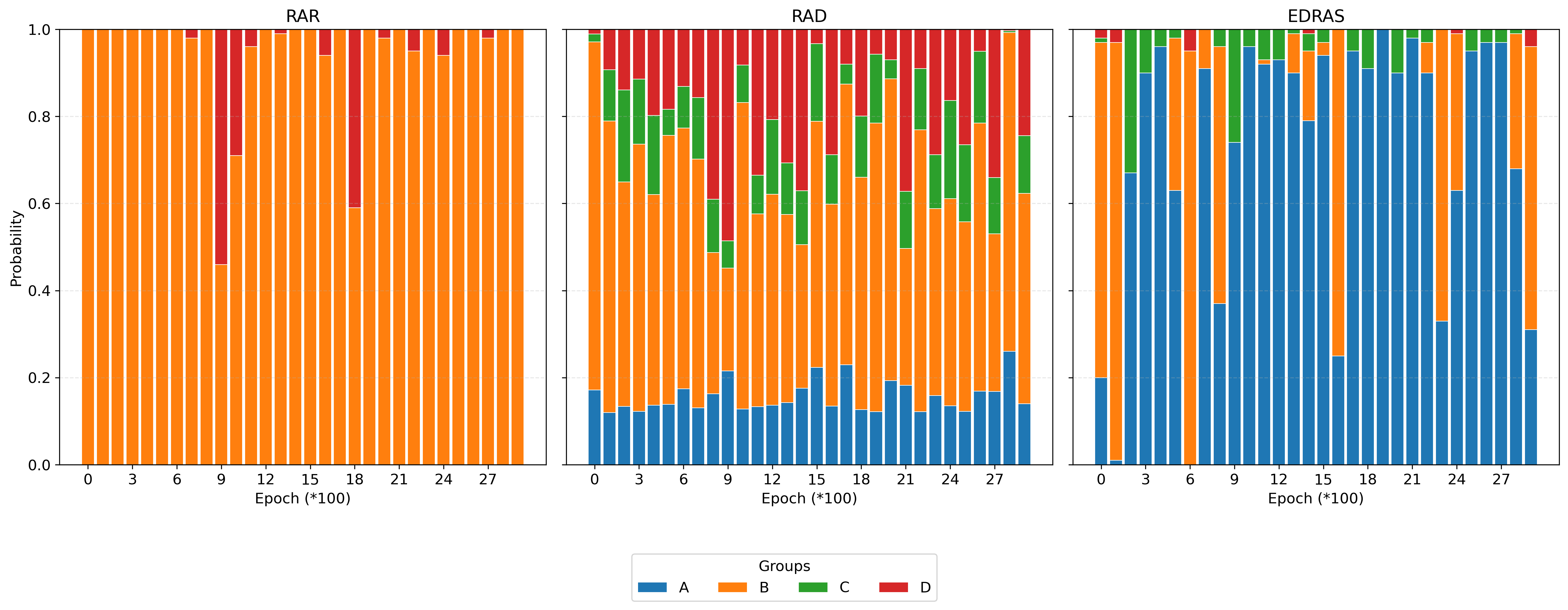}
      \caption{Probabilities of the 4 groups for RAR, RAD and EDRAS in domain $[0.6, 0.8]\times [-1,1]$ plotted over training epochs. }
 \label{fig:prob_view08}
\end{figure}

\end{document}